\documentclass[review]{elsarticle}

\usepackage[hidelinks,hyperfootnotes=false, bookmarks]{hyperref}
\usepackage{graphicx} 
\graphicspath{{../ps/}}
\usepackage{mathtools}
\usepackage{amsmath}
\usepackage{amssymb}
\usepackage{placeins}
\usepackage{algorithm}
\usepackage{algpseudocode}
\usepackage[utf8]{inputenc}
\usepackage{multirow}
\usepackage{color,soul}
\usepackage{array}
\biboptions{sort&compress}
\usepackage{import}
\usepackage{booktabs}
\usepackage{colortbl}
\usepackage{graphicx}
\usepackage[dvipsnames]{xcolor}
\usepackage{epstopdf}
\usepackage{mathrsfs}
\usepackage[figurename=Fig.]{caption}
\usepackage[numbers]{natbib}
\usepackage{array}
\usepackage{amsfonts}
\usepackage{scalerel}
\usepackage{bm}
\usepackage[labelformat=simple]{subcaption}
\usepackage[scientific-notation=true]{siunitx}

\newcommand{\grayrule}{\arrayrulecolor{gray!20} \midrule \arrayrulecolor{black}}
\newcommand{\mycmidrule}{\arrayrulecolor{gray!20} \cmidrule{2-4} \arrayrulecolor{black}}

\usepackage{blindtext, xcolor}
\usepackage{comment}

\DeclareMathOperator{\tr}{tr}
\newcommand{\arginf}{\operatornamewithlimits{arg\,inf}}
\newcolumntype{M}[1]{>{\centering\arraybackslash}m{#1}}

\newcolumntype{L}{>{\centering\arraybackslash}m{3cm}}
\usepackage[margin=1in]{geometry}
\usepackage{lineno,hyperref}
\algblock{Input}{EndInput}
\algnotext{EndInput}
\algblock{Output}{EndOutput}

\newcommand{\overbar}[1]{\mkern 1.5mu\overline{\mkern-1.5mu#1\mkern-1.5mu}\mkern 1.5mu}

\usepackage{bm}
\usepackage[color=gray!20, backgroundcolor=yellow, textwidth=2cm, textsize=footnotesize]{todonotes}
\renewcommand*{\today}{ }

\begin{document}

\begin{frontmatter}

%\title{Energetic formulation and numerical implementation of fatigue with cyclic plasticity and a phase-field approach to fracture}

\title{Phase-field modeling of fatigue coupled to \\ cyclic plasticity  in an energetic formulation\tnoteref{t1,t2}}

\author[add1]{Jacinto Ulloa}
\cortext[JU]{Corresponding author: Jacinto Ulloa}
\ead{jacintoisrael.ulloa@kuleuven.be}

\author[add1]{Jef Wambacq}
\ead{jef.wambacq@kuleuven.be}

\author[add2]{Roberto Alessi}
\ead{roberto.alessi@unipi.it}

\author[add1]{Geert Degrande}
\ead{geert.degrande@kuleuven.be}

\author[add1]{Stijn François}
\ead{stijn.francois@kuleuven.be}

\address[add1]{KU Leuven, Department of Civil Engineering, Kasteelpark Arenberg 40, B-3001 Leuven, Belgium}
\address[add2]{Università di Pisa, Department of Civil and Industrial Engineering, Largo Lucio Lazzarino 2, 56122 Pisa, Italy}

\tnotetext[t1]{{\itshape Postprint version.}}
\tnotetext[t2]{{\itshape Published version:}
       J.~Ulloa,  J.~Wambacq, R.~Alessi, G.~Degrande, and S.~François.
       \newblock Phase-field modeling of fatigue coupled to cyclic plasticity  in an energetic formulation.
       \newblock {\em Computer Methods in Applied Mechanics and Engineering}, 373:113473, 2021.\\
       {\itshape DOI: \tt\url{https://doi.org/10.1016/j.cma.2020.113473}} \vspace{3mm}}

\begin{abstract}

%The energetic formulation for rate-independent systems provides a foundation for the thermodynamically-consistent modeling of dissipative solids. The present study employs this framework in the development of  a variational model that describes cyclic failure in brittle and ductile materials. This is achieved by 

This paper presents a modeling framework to describe the driving mechanisms of cyclic failure in brittle and ductile materials, including cyclic plasticity and fatigue crack growth. A variational model is devised using the energetic formulation for rate-independent systems, coupling a phase-field description of fatigue fracture to a cyclic plasticity model that includes multi-surface kinematic hardening, gradient-enhanced isotropic hardening/softening and ratcheting. The coupled model embeds two distinctive fatigue effects. The first captures the characteristic features of low-cycle fatigue, driven by the accumulation of plastic strains, while the second  accounts for high-cycle fatigue, driven by free energy accumulation. The interplay between these mechanisms allows to describe a wide range of cyclic responses under both force loading and displacement loading, as shown in several numerical simulations. Moreover, the phase-field approach to fracture accounts for the initiation and propagation of fatigue-induced cracks. 

\end{abstract}


\end{frontmatter}

%\linenumbers

\section{Introduction}

%% Motivation and brief description of the objective
Solids and structures subjected to cyclic loading exhibit a progressive reduction in load-carrying capacity due to material degradation. This phenomenon is important in several branches of engineering, accounting for (up to) 90\% of all structural failures~\citep{stephens2000}. Despite its importance, the study of cyclic failure, including cyclic plasticity and various fatigue regimes, remains an open issue in computational mechanics. The objective of this study is to propose a model to describe cyclic failure in brittle and ductile materials using a mathematically and physically sound variational framework, and to demonstrate its capabilities in benchmark numerical~simulations. % embedded in continuum mechanics. In particular, models that consistently describe the dissipative behavior of cyclically deforming materials as well as the evolution of fatigue cracks are significantly limited. 

%% Describing fatigue
Depending on the loading conditions and material properties, different fatigue regimes can be distinguished~\cite{lemaitre1996}. At load amplitudes above a certain threshold, yet small enough to avoid plastic strains, high-cycle fatigue occurs, governed by a slow material degradation that leads to brittle fractures. Under higher load amplitudes, low-cycle fatigue occurs. This phenomenon is driven by a combination of damage and plastic strains, where the initial cyclic response can be rather complex. In particular, cyclic hardening or cyclic softening effects are generally observed, as well as asymmetrical loading effects, such as ratcheting under force loading and stress relaxation under displacement loading~\cite{ hassan1994a,hassan1994b,chaboche2008,paul2010}. These cyclically plastic responses may occur in a transient fashion, leading to stabilized hysteresis~loops. The interplay between these phenomena and material degradation strongly influences material behavior in the low-cycle fatigue regime. %, where the material may still be subject to high-cycle fatigue. 

%% Issues with current fatigue methods 
Classical fatigue analyses are based on (semi-) empirical methods, which require extensive data from experimental tests. Fatigue life is commonly assessed using  Wöhler curves, relating the applied stress amplitude to the number of cycles to failure in constant-amplitude cyclic loading. Statistical approaches such as the Basquin relation~\citep{basquin1910} fit empirical equations to Wöhler curves. The drawback of these techniques lies in their empirical nature, requiring tuning of problem-dependent parameters that do not account for the underlying fracture process. Concerning fatigue fracture analysis, the conventional approach is based on Paris' law \cite{paris1963} and its extensions \cite{forman1967,maierhofer2014}, which relate crack growth rate to stress intensity factors. Paris' law is rooted in classical fracture mechanics and is therefore unable to describe crack initiation and the final rupture stage. An alternative approach to fatigue consists of introducing fatigue effects in constitutive material models \citep{desmorat2006}, providing a more versatile framework \cite{Alessi2017fatigue}. However, most models include parameters with no clear physical interpretation, and often use Paris-type laws and Wöhler curves as input data~\cite{carrara2020}.% from specific structural configurations.%the cyclic increment of  rather than outputs 

%FWO: The main limitation of classical fracture mechanics is that the analysis is based on a predefined crack, which implies that the material evolution during the full loading process is not considered. As a result, the process of crack initiation cannot be consistently described. From a computational standpoint, classical fracture
%mechanics is hardly applicable to complex fracture processes such as crack branching and crack merging.
%%commonly characterized using stress-cycle (S/N) curves. S/N curves provide an empirical relationship between the number of cycles to failure for a corresponding maximum (normalized) stress level, but require a lot of time-consuming high-cycle fatigue tests to calibrate.

%% The material modeling approach 
%%%%% summarize and improve
This study adopts the framework of constitutive models with non-local internal variables. To properly describe failure modes, such as brittle and quasi-brittle fractures, or shear bands in ductile materials, a suitable representation of highly localized strains is required. Strain localization is captured by models with softening behavior, which, however, render the mathematical problem ill-posed and result in pathological mesh-dependence in numerical simulations. To alleviate these issues, regularization can be introduced by means of non-local effects that are governed by internal length scales. Gradient-enhanced models, as outlined in the formulations of \citet{maugin1990} and \citet{miehe2011}, are common examples of this approach.
% From a physical standpoint, this approach allows to account for interactions at the material's microstructure \citep{bazant1997}.
%From a physical standpoint, this approach allows to account for interactions at the material's microstructure \citep{bazant1997}. %In this setting, a damage variable is introduced to account for material degradation. Moreover, a suitable representation of low-cycle fatigue involves cyclic plasticity, where most theories rely on the Armstrong–Frederick model \cite{armstrong1966,fred2007} and its modifications \citep{chaboche1986,chaboche1989,chaboche1991,chaboche2008,houlsby2017}. 

In relation to gradient-enhanced models, developments in the modeling of crack nucleation and propagation have taken place in the last two decades due to the variational approach to fracture \citep{FrancMar1998,BourFrancMar2000,BourFrancMar2008}. This framework links classical fracture mechanics to gradient-damage models. Specifically, the $\Gamma$-converging regularization of \citet{BourFrancMar2000} towards the Griffith-based energy functional of~\citet{FrancMar1998} can be alternatively viewed as a phase-field description of fracture~\cite{MieHofWel2010,Bord-Hugh2012,ambati2015review} or as a gradient-damage model \cite{pham2011,pham2011b,marigo2016,deborst2016}. This framework is able to naturally predict crack initiation and propagation in complex crack topologies, overcoming limitations of classical fracture mechanics and discrete-crack approaches, such as XFEM-based methods \citep{MoeDolBely1999,MoeBely2002,SamBely2005}. While initially developed for brittle fracture,  several extended phase-field models have been developed in the literature, some of which are outlined in~\citet{wu2018phase}. 
%Examples of a rather extensive list of works for plasticity and damage include \citep{ muhlhaus1991,DeBorstMull1992,comi1996,fremond1996,peerlings1996,comi1999,gurtin2005,fleck2009a,fleck2009b}, while general formulations are outlined in \citep{maugin1990,miehe2011}. 

Two extensions of the phase-field approach to fracture are particularly relevant for the purposes of this work. The first is the extension to plasticity~\citep{AleMarVid2015,AleMarMauVid2017,Duda2015,AmbGerDeL2015,Bord-Hugh2016,kuhn2016,miehe2017,UllRodSam2016,rodriguez2018,fang2019,alessi2018comparison,yin2020}, which allows to describe the nucleation and propagation of ductile cracks. The second is the extension to fatigue that has gained attention in the recent literature~\citep{Alessi2017fatigue,carrara2020,mesgarnejad2018,seiler2019,lo2019,haveroth2020,schreiber2020,loew2020a,loew2020b}. This approach to fatigue has been shown to consistently recover experimental observations such as Wöhler curves and Paris' law. Nevertheless, extensions of phase-field-based fatigue to account for the main features of low-cycle fatigue are still missing in the literature. A suitable representation of low-cycle fatigue must include cyclic plasticity, for which most theories~\citep{chaboche1986,chaboche1989,chaboche1991,chaboche2008,houlsby2017} can be linked to the Armstrong–Frederick model \cite{armstrong1966,fred2007}. It is the objective of the present work to address this topic. % as post-processing results

% Modeling framework: classical vs. variational
%, a well-established framework for dissipative material behavior in the context of thermomechanics with internal variables.
The governing equations of a large class of material models can be derived from  the theory of generalized standard materials \cite{moreau1970,halphen1975,germain1983,ziegler1987}. In this setting, the evolution problem follows from an internal energy density and a dissipation potential, ensuring an a priori fulfilment of the second law of thermodynamics. This theory can be reformulated in a variational setting, where the governing equations emerge as the Euler-Lagrange equations of an energy minimization principle that can be solved using numerical optimization techniques. Pioneering works on this topic focused on rate-type variational principles for local elastoplasticity \citep{simo1989, comi1995, han1999, hackl1997, ortiz1999, miehe2002,carstensen2002}, and have been extended to, e.g., gradient-enhanced models \citep{comi1999, gurtin2009, miehe2011,fleck2009a}. In these works, the principle of virtual power is generally used to derive the governing equations~\citep{maugin1980,petryk2003,miehe2011}. 
%: the dissipation potential can be recovered as the support function of the set of admissible stresses
%Two equivalent approaches are possible: the primal form, which makes direct use of dissipation potentials, and  the dual form, which explicitly defines the set of admissible generalized stresses using yield functions. Both formulations are directly related to each other by standard arguments of convex analysis (see, e.g., \cite{han1999} and \cite{hackl2007} for details).

The energetic formulation \citep{mielke2006,Mielke2015} is a particularly attractive variational framework that furnishes a unified and rigorous mathematical setting for rate-independent dissipative processes. The theory handles non-smooth evolutions (discontinuities in space and time), since no derivatives appear in the most general setting, and provides tools for the analysis of structural and material stability \cite{Alessi2016}. Some applications of the energetic formulation in the modeling of dissipative solids include plasticity~\citep{lancioni2015, rokovs2016}, (quasi-)brittle \citep{BourFrancMar2008,pham2011,marigo2016,luege2018} and ductile \citep{AleMarVid2014,AleMarMauVid2017,UllRodSam2016,rodriguez2018} fracture, shape-memory alloys \citep{alessi2015shape,alessi2016shape} and fatigue \citep{Alessi2017fatigue,carrara2020}.

% Non-associativity: a problem for variational models
An important limitation of the theory of generalized standard materials and its variational formulation is the restriction imposed by the principle of maximum dissipation. This condition implies the normality law for the evolution equations and thus, in principle, excludes non-associative models. Nevertheless, extensions to non-associativity can be made by considering certain state-dependent dissipation potentials. References in this line of work include \citep{mosler2009,mosler2010,francfort2013,francfort2018,luege2018}. These developments are essential from a practical standpoint because a variational structure is no longer restricted to associative models, accounting for more realistic representations of, for instance, plasticity in geomaterials and in metals under cyclic loading.

In the framework of the energetic formulation, we propose a model that couples the phase-field approach to fatigue, as suggested by \citet{Alessi2017fatigue} and \citet{carrara2020}, to the main mechanisms of cyclic plasticity. The proposed model includes multi-surface kinematic hardening and a non-associative ratcheting variable \citep{houlsby2017}, as well as gradient-enhanced isotropic hardening/softening, governed by its own characteristic length scale. The coupling of the phase-field-based fatigue model to cyclic plasticity renders a general framework that encompasses the characteristic behavior of both high- and low-cycle fatigue in a  thermodynamically consistent setting and is able to objectively represent fatigue crack growth with plasticity. Several examples are presented to highlight the versatility of the proposed framework to capture representative responses under both force loading and displacement loading. 

This paper is organized as follows. To establish the modeling framework, section 2 presents the energetic formulation applied to a general class of dissipative solids with gradient-enhanced internal variables. These concepts are then used to construct the fatigue model with cyclic plasticity in section 3, with numerical experiments presented in section 4. 

%Using a compact tensor notation, a dot ($\cdot$) and a colon ($:$) denote products between vectors and tensors with simple and double contraction, respectively. The Euclidean norm of a vector- or tensor-valued quantity is denoted as $\Vert\Box\Vert$. The notation $\Box_\mathrm{dev}\coloneqq\Box-\tr\Box\bm{1}/3$ is used to denote the deviatoric projection of a tensor $\Box$.  Moreover, the notation $\Vert\bullet\Vert_ {\Box}\coloneqq\sqrt{\bullet:\Box:\bullet}$ for a second-order tensor $\bullet$ is used, where $\Box$ is a fourth-order tensor. Space and time dependence are denoted with vector and scalar arguments, respectively, e.g., $\Box(\bm{x})$ and $\Box(t)$. Functions and function values are denoted using the same symbol: $\Box$ is a space-time-dependent function evaluated as $\Box(\bm{x},t)$. The distinctions among functions and their values are  made when required for clarity purposes, and are otherwise omitted and inferred from context. An overdot, $\dot{\Box}$, is used to denote the derivative with respect to time, while the symbol $\nabla$ is used to denote the spatial gradient. In the context of convex analysis, $\partial\Box(\bullet)$ denotes the multivalued subdifferential of $\Box$ at $\bullet$. Further notational specifications are made as required throughout the~text.

We employ the following notations. A dot ($\cdot$) and a colon ($:$) denote inner products with simple and double contraction, respectively. The Euclidean norm of $\Box$ is denoted as $\Vert\Box\Vert$. The notation $\Box_\mathrm{dev}\coloneqq\Box-(1/3)\tr\Box\bm{1}$ is used to denote the deviatoric part of a tensor $\Box$.  The notation $\Vert\bullet\Vert_ {\Box}\coloneqq\sqrt{\bullet:\Box:\bullet}$ for a second-order tensor $\bullet$ is used, where $\Box$ is a fourth-order tensor.  Functions and function values are denoted using the same symbol: $\Box$ is a space-time-dependent function evaluated as $\Box(\bm{x},t)$. Time-dependent functions evaluated at a point in space are denoted with a vector argument, e.g. $\Box(\bm{x})$, while space-dependent functions parametrized by time are denoted with a scalar argument, e.g. $\Box(t)$. These distinctions are made when required for clarity purposes, and are otherwise omitted and inferred from context. The time derivative  is denoted as $\dot{\Box}$, while the spatial gradient reads $\nabla\Box$. In the context of convex analysis, $\partial_\Box$ is the multivalued subdifferential with respect to $\Box$, while $\partial\Box(\bullet)$ is the subdifferential of $\Box$ at~$\bullet$. Finally, $\delta_{\Box}\coloneqq\partial_{\Box}-\mathrm{div}[\partial_{\nabla{\Box}}]$ is the spatial Euler-Lagrange derivative with respect to $\Box$.

\section{Variational framework}\label{variat_general}

%This section presents the general framework for dissipative gradient-enhanced material models hereby considered, recalling concepts from the theory of generalized standard materials and the energetic formulation. 

\subsection{Problem outline}\label{referenceprob}

Consider an arbitrary solid of mass density $\rho$ occupying a domain $\Omega\subset\mathbb{R}^d$ of dimension $d\in\{1,2,3\}$, with boundary $\Gamma\subset\mathbb{R}^{d}$. The boundary consists of a Dirichlet part ${\Gamma_\mathrm{D}}$ with imposed displacements $\overbar{\bm{u}}(\bm{x},t)\in\mathbb{R}^d$ and a Neumann part $\Gamma_\mathrm{N}$ with imposed tractions $\bar{\bm{t}}(\bm{x},t)\in\mathbb{R}^d$, such that ${\Gamma_\mathrm{D} \cup \Gamma_\mathrm{N} = \Gamma}$ and ${\Gamma_\mathrm{D} \cap \Gamma_\mathrm{N} =\varnothing}$. The solid may be subjected to volume forces per unit mass, denoted as $\bm{b}(\bm{x},t)\in\mathbb{R}^d$.

The deformation process is assumed to be quasi-static,  occurring in a pseudo-time (loading) interval ${\mathrm{T}\coloneqq[0,t_\mathrm{max}]}$. The displacement field is given by $$\bm{u}\colon\begin{dcases}\Omega\times\mathrm{T}\to\mathbb{R}^d, \\ (\bm{x},t)\mapsto\bm{u}(\bm{x},t),\end{dcases}$$ which is kinematically admissible, satisfying  boundary conditions on $\Gamma_\mathrm{D}$.
%
%\begin{equation}
%\bm{u}(x\in\Gamma_\mathrm{D},t)=\bar{\bm{u}}(t),
%\end{equation}
%\\
%where $\bar{\bm{u}}(t)$ can be viewed as the trace on \Gamma_{\mathrm{D}} of a field $\bar{\bm{u}}(\bm{x},t)$ defined over $\Omega$. 
Assuming the small-strain hypothesis, the compatible strain tensor $\boldsymbol{\varepsilon}(\bm{x},t)\in\mathbb{R}^{d\times d}_\mathrm{sym}\coloneqq\{\bm{e}\in \mathbb{R}^{d\times d} \, | \,  \bm{e}=\bm{e}^\mathrm{T}\}$ is obtained from  $$\boldsymbol{\varepsilon}\coloneqq \nabla^{\mathrm{s}}\boldsymbol{u}, \quad \text{with} \quad \nabla^{\mathrm{s}}\Box\coloneqq\frac{1}{2}(\nabla\otimes\Box+\Box\otimes\nabla).$$ 
The Cauchy stress tensor $\boldsymbol{\sigma}(\bm{x},t)\in\mathbb{R}^{d\times d}_\mathrm{sym}$ is statically admissible, satisfying equilibrium for all $t\in \mathrm{T}$:
\begin{equation}
\mathrm{div}\bm{\sigma}+\rho\bm{b} = \boldsymbol{0} \quad \text{in} \quad \Omega \quad \text{and} \quad  \bm{\sigma}\cdot\boldsymbol{n}=\bar{\bm{t}} \quad \text{on} \quad \Gamma_{\mathrm{N}}, \quad	\text{with} \quad \bm{u}=\overbar{\bm{u}} \quad \text{on} \quad \Gamma_{\mathrm{D}}.
\label{st_adm}
\end{equation}
The dissipative mechanisms that lead to inelastic material behavior are characterized by a generic set of internal variables and their spatial gradients: $$\mathbf{a}\colon\begin{dcases}\Omega\times\mathrm{T}\to \mathbb{R}^{m},\\ (\bm{x},t)\mapsto\mathbf{a}(\bm{x},t),\end{dcases} \quad \nabla\mathbf{a}:\begin{dcases}\Omega\times\mathrm{T}\to \mathbb{R}^{md},\\ (\bm{x},t)\mapsto\nabla\mathbf{a}(\bm{x},t).\end{dcases}$$ We denote the set of primary fields by ${\mathbf{q}\coloneqq\{\bm{u},\mathbf{a}\}}$ and the constitutive state by ${\mathbf{c}\coloneqq\{\bm{\varepsilon},\mathbf{a},\nabla\mathbf{a}\}}$. Here, $\mathbf{a}$ is a vector arrangement of $m$ components associated to both scalar- and tensor-valued internal variables.

\subsection{Generalized standard materials}

%The present variational framework adopts the theory of generalized standard materials, [cite], which relies on the definition of two energy quantities: an internal energy density and a dissipation potential. 

Let $\psi\coloneqq\psi(\bm{\varepsilon},\mathbf{a},\nabla\mathbf{a})$ denote a Helmholtz-type internal energy density. To ensure physical soundness, the second law of thermodynamics is taken as an a priori restriction, given, for isothermal processes, by the Clausius-Planck inequality 
\begin{equation}
\bm{\sigma}:\dot{\bm{\varepsilon}} -\dot{\psi}(\bm{\varepsilon},\mathbf{a},\nabla\mathbf{a})\geq0.
\label{CD}
\end{equation}
The constitutive stress-strain relation
\begin{equation}
\bm{\sigma}=\frac{\partial\psi}{\partial{\bm{\varepsilon}}}(\bm{\varepsilon},\mathbf{a},\nabla\mathbf{a})
\label{constitutive}
\end{equation}
directly follows from equation~\eqref{CD}, along with the dissipation rate inequality
\begin{equation}
\phi=\mathbf{s}\cdot\dot{\mathbf{a}}\geq0, \quad \text{with} \quad \mathbf{s}=-\delta_{\mathbf{a}}\psi(\bm{\varepsilon},\mathbf{a},\nabla\mathbf{a}).
\label{dualforces} 
\end{equation}
%\quad \text{with} \quad \mathbf{s}\coloneqq-\frac{\partial}{\partial{\mathbf{a}}}\psi(\bm{\varepsilon},\mathbf{a}). (but assumed) 
The set $\mathbf{s}$ contains the thermodynamic forces, or generalized stresses, conjugate to $\mathbf{a}$. 

A thermodynamically admissible dissipation potential is defined as $\phi\coloneqq\phi(\dot{\mathbf{a}},\nabla{\dot{\mathbf{a}}};\mathbf{c})\geq0$, which is assumed to be convex with respect to $\{\dot{\mathbf{a}},\nabla{\dot{\mathbf{a}}}\}$ and to vanish for null rates. The dependence of the dissipation potential on the state $\mathbf{c}$ is a generality that accounts for a wide class of material models~\citep{germain1983}. The dissipation potential may also depend on history variables, not included in $\mathbf{c}$ but derived from its time history~\citep{Alessi2017fatigue,alessi2018c}, as well as on the generalized stresses $\mathbf{s}$, as for the case of non-associative models~\citep{houlsby2007,francfort2018,luege2018}. These more general cases will emerge in the model proposed in section~\ref{model}. For notational simplicity and without losing generality, only dependence on $\mathbf{c}$ is considered in the abstract formulation presented in this section. For rate-independent processes, $\phi$ is a homogeneous function of first degree  in $\{\dot{\mathbf{a}},\nabla{\dot{\mathbf{a}}}\}$, such that 
$$\phi(b\dot{\mathbf{a}},b\nabla{\dot{\mathbf{a}}};\mathbf{c})=b\phi(\dot{\mathbf{a}},\nabla{\dot{\mathbf{a}}};\mathbf{c}), \quad \forall \, b\geq0.$$ As a consequence, $\phi$ is not differentiable at null rates, and, from equation~\eqref{dualforces}, it follows that $$\mathbf{s}\in\partial_{\dot{\mathbf{a}}}\phi(\dot{\mathbf{a}},\nabla{\dot{\mathbf{a}}};\mathbf{c})-\mathrm{div}\big[\partial_{\nabla\dot{\mathbf{a}}}\phi(\dot{\mathbf{a}},\nabla{\dot{\mathbf{a}}};\mathbf{c})\big].$$ The evolution equation of $\mathbf{a}$ then takes the form
\begin{equation}
\delta_{{\mathbf{a}}}\psi(\bm{\varepsilon},\mathbf{a},\nabla\mathbf{a}) + \delta_{\dot{\mathbf{a}}}\phi(\dot{\mathbf{a}},\nabla\dot{\mathbf{a}};\mathbf{c})\ni \bm{0},
\label{biot}
\end{equation}% 
which is often referred to as Biot's equation~(cf. \citep{biot1965}). Equations \eqref{st_adm} and \eqref{biot} represent the strong form of the evolution problem of a general dissipative material model with gradient-enhanced internal variables.  

\subsection{Energetic formulation}\label{energetic_form}

In this section, the governing equations of the evolution problem are recovered in a variational setting. To this end, the basic energy quantities are first introduced in global form. Then, the evolution problem is defined in terms of the energetic formulation for rate-independent systems~\citep{Mielke2015}, based on notions of energy balance and stability. These principles naturally lead to an incremental energy minimization problem and constitute the building blocks of the variational fatigue model presented in section~\ref{model}. %The reader is referred to  for a rigorous survey of the energetic formulation. 

\subsubsection{Energy quantities}

To characterize how the system stores and dissipates energy in exchange with the external environment, the energetic formulation begins with the definition of global energy functionals. To this end, we denote a general function space of internal variables by ${\mathscr{A}}$ and the corresponding space of admissible variations that embeds evolution constraints (e.g., irreversibility conditions) by $\tilde{\mathscr{A}}$. The specific form of these function spaces depends on the material model, as will be clear in section~\ref{model}. The space of kinematically admissible displacement fields and the corresponding space of admissible variations are given~by
\begin{equation}
\mathscr{U}(t)\coloneqq\{\bm{w}\in \mathscr{F}\,|\,\bm{w}=\bar{\bm{u}}(t) \,\,\, \text{on} \,\,\, \Gamma_\mathrm{D}\} \quad \text{and} \quad \tilde{\mathscr{U}}\coloneqq\{\tilde{\bm{w}}\in \mathscr{F} \,|\,\tilde{\bm{w}}=\bm{0} \,\,\, \text{on} \,\,\, \Gamma_\mathrm{D}\},
\label{dispspace} 
\end{equation}
where the form  of~$\mathscr{F}$ also depends on the model and is specified in section~\ref{model}. The space of primary fields then reads ${\mathscr{Q}\coloneqq\mathscr{U}\times\mathscr{A}}$, such that ${\mathbf{q}(t)\in\mathscr{Q}}$ is a process with admissible variations~$\tilde{\mathbf{q}}\in{\tilde{\mathscr{Q}}\coloneqq\tilde{\mathscr{U}}\times\tilde{\mathscr{A}}}$. %, and the corresponding space of gradient fields by ${\mathscr{G}}(\Omega,G)$
%
%\begin{equation}
%\mathcal{P}(t,\mathbf{q}) \coloneqq \mathcal{E}(t,\mathbf{q}),
%\label{potfunc}
%\end{equation}

The stored internal energy functional $\mathcal{E}\colon\mathscr{Q}\to\mathbb{R}\cup\{+\infty\}$ is given by the state function
\begin{equation}
\mathcal{E}(\mathbf{q})\coloneqq \int_{\Omega}\psi\big(\bm{\varepsilon}(\bm{x}),\mathbf{a}(\bm{x}),\nabla\mathbf{a}(\bm{x})\big) \, \mathrm{d}\bm{x},
\label{stored_ext0}
\end{equation}
while the work of external actions is defined as the time-integral of the external power, namely:
\begin{equation}
\mathcal{L}\big(\bm{u};[0,t]\big)\coloneqq\int_{0}^{t}\bigg[\int_\Omega\rho\bm{b}(\bm{x},\tau)\cdot\dot{\bm{u}}(\bm{x},\tau)\,\mathrm{d}\bm{x} + \int_{\Gamma_\mathrm{N}}\bar{\bm{t}}(\bm{x},\tau)\cdot\dot{\bm{u}}(\bm{x},\tau)\,\mathrm{d}S + \int_{\Gamma_{\mathrm{D}}}\bm{t}_\mathrm{r}(\bm{x},\tau)\cdot\dot{\bar{\bm{u}}}(\bm{x},\tau)\, \mathrm{d}S\,\bigg]\mathrm{d}\tau,
\label{stored_ext1}
\end{equation}
where $\bm{t}_\mathrm{r}$ is the traction vector on $\Gamma_{\mathrm{D}}$. To formulate a stability condition, we introduce the external work distance~$\mathcal{T}\colon\mathrm{T}\times\mathscr{U}\times\mathscr{U}\to\mathbb{R}$, that is, the work done by the external forces at a given time between two admissible states $\{\bm{u}_0,\bm{u}_1\}$:
\begin{equation}
\mathcal{T}(t,\bm{u}_0,\bm{u}_1)=\int_{\Omega}\rho\bm{b}(\bm{x},t)\cdot\big({\bm{u}}_1(\bm{x})-{\bm{u}}_0(\bm{x})\big)\, \mathrm{d}\bm{x} +\int_{\Gamma_{\mathrm{N}}}\bar{\bm{t}}(\bm{x},t)\cdot\big({\bm{u}}_1(\bm{x})-{\bm{u}}_0(\bm{x})\big)\, \mathrm{d}S.
\label{extpower}
\end{equation}
On the other hand, the dissipative power functional ${\mathcal{R}\colon\tilde{\mathscr{A}}\times\mathscr{Q}\to[0,+\infty]}$ is defined as
\begin{equation}
\begin{aligned}
\mathcal{R}(\dot{\mathbf{a}};\mathbf{q})\coloneqq\int_{\Omega}\phi\big(\dot{\mathbf{a}}(\bm{x}),\nabla{\dot{\mathbf{a}}}(\bm{x});\mathbf{c}(\bm{x})\big)\, \mathrm{d}\bm{x},
\end{aligned}
\label{dissipated}
\end{equation}
while the dissipation distance ${\mathcal{D}\colon\mathscr{Q}\times\mathscr{Q}\to[0,+\infty]}$ between two states ${\{\mathbf{q}_0,\mathbf{q}_1\}}$ reads as follows~\cite{mielke2006}:
%
%\begin{equation}
%\mathcal{D}(\mathbf{q}_0,\mathbf{q}_1)\coloneqq\inf\bigg\{\mathrm{Diss}_\mathcal{D}(\mathbf{q};[0,1]) \,\, \big\vert \,\, \mathbf{q}\in \mathrm{C}^{1}([0,1],\mathrm{Q}), \, \mathbf{q}(0)=\mathbf{q}_0,  \ \mathbf{q}(1)=\mathbf{q}_1 \bigg\},
%\label{distance0}
%\end{equation} 
%\begin{equation}
%\mathcal{D}(\mathbf{q}_0,\mathbf{q}_1)\coloneqq\int_{\Omega}\inf\bigg\{\int_0^1\phi\big(\dot{\mathbf{a}}(\bm{x},s),\nabla{\dot{\mathbf{a}}}(\bm{x},s);\mathbf{c}(\bm{x},s)\big) \, \mathrm{d}s \, \vert \, \mathbf{q}\in \mathrm{C}^{1}(\Omega\times[0,1],\mathrm{Q}), \, \mathbf{q}(0)=\mathbf{q}_0,  \ \mathbf{q}(1)=\mathbf{q}_1 \bigg\}\, \mathrm{d}\bm{x},
%\label{distance0}
%\end{equation} 
%
\begin{equation}
\mathcal{D}(\mathbf{q}_0,\mathbf{q}_1)\coloneqq\inf\bigg\{\int_0^1\mathcal{R}\big(\dot{\mathbf{a}}(s);\mathbf{q}(s)\big) \, \mathrm{d}s \,\, \big\vert \,\, \mathbf{q}\in \mathrm{C}^{1}([0,1],\mathrm{Q}), \, \mathbf{q}(0)=\mathbf{q}_0,  \ \mathbf{q}(1)=\mathbf{q}_1 \bigg\},
\label{distance0}
\end{equation} 
with $\mathrm{Q}\subset\mathbb{R}^{d}\times\mathbb{R}^{m}$. This quantity allows to measure the energy dissipated along arbitrary minimizing paths $\mathbf{q}(s)\in\mathscr{Q}$. For smooth evolutions, the total energy dissipated in $[0,t]$ reads
\begin{equation}
\mathrm{Diss}_\mathcal{D}(\mathbf{q};[0,t]) \coloneqq \int_{0}^{t}\mathcal{R}\big(\dot{\mathbf{a}}(\tau);\mathbf{q}(\tau)\big)\, \mathrm{d}\tau.
\label{dissipenerg}
\end{equation}
For certain dissipation potentials obtained as the time-derivative of an energy function, the dissipated energy is a state function~\citep{AleMarMauVid2017,marigo2016}. However, as in the present study, this quantity is generally path-dependent.  

%The triangle inequality is assumed as a basic requirement for the dissipation distance, so that $\mathcal{D}(\mathbf{q}_0,\mathbf{q}_2)\leq \mathcal{D}(\mathbf{q}_0,\mathbf{q}_1)+\mathcal{D}(\mathbf{q}_1,\mathbf{q}_2)$ for all admissible states $\mathbf{q}_0,\mathbf{q}_1,\mathbf{q}_2$.

%The global  dissipation distance allows to evaluate the global dissipated energy during a time interval $[0,t]$ as
%
%\begin{equation}
%\mathrm{Diss}_\mathcal{D}(\mathbf{q};[0,t])\coloneqq\sup\bigg\{\sum_{n=1}^{m}\mathcal{D}(\mathbf{q}_{n-1},\mathbf{q}_{n}) \ \vert \ n\in\mathbb{N}, \ 0=t_0<t_1< \dots < t_m = t \bigg\}.
%\end{equation}
%\\

\subsubsection{Evolution problem}\label{energ_evo}
From the basic energy quantities, the evolution problem can be defined in variational form by means of the energetic formulation. This framework is based on a notion of energy balance and either a global or a local stability condition, and is directly related to the incremental energy minimization problem of~section~\ref{incremental}. 

\paragraph*{Global formulation}

A process $\mathbf{q}\colon\mathrm{T} \to\mathscr{Q}$ is an energetic solution if for all $t\in \mathrm{T}$:
\begin{align}
&\mathcal{E}\big(\mathbf{q}(t)\big) \leq \mathcal{E}(\tilde{\mathbf{q}}) -\mathcal{T}(t,{\bm{u}}(t),\tilde{\bm{u}}) + \mathcal{D}\big(\mathbf{q}(t),\tilde{\mathbf{q}}\big),\quad \forall \, \tilde{\mathbf{q}}=\{\tilde{\bm{u}},\tilde{\mathbf{a}}\}\in\mathscr{Q},
\label{GS}
%\tag{GS}
\\
%&\mathcal{E}\big(\mathbf{q}(t)\big) +  \mathrm{Diss}_\mathcal{D}(\mathbf{q};[0,t]) = \mathcal{E}\big(\mathbf{q}(0)\big)  +  \mathcal{L}\big(\bm{u};[0,t]\big)
%\label{EB}
%\tag{EB}
%\\
&\mathcal{E}\big(\mathbf{q}(t)\big) +  \mathrm{Diss}_\mathcal{D}(\mathbf{q};[0,t]) = \mathcal{E}\big(\mathbf{q}(0)\big)  +  \mathcal{L}\big(\bm{u};[0,t]\big).
\label{EB}
%\tag{EB}
\end{align}
Equations \eqref{GS} and \eqref{EB} represent, respectively, the global stability condition and global energy balance, and constitute the most general form of the energetic formulation, with regularity assumptions only required for the external loading functions. 

%Let $\check{\mathbf{a}}\subseteq\mathbf{a}$ denote a subset of scalar-valued internal variables with irreversible evolutions. Due to rate-independence, \eqref{DI} demands that 
%
%\begin{equation}
%\partial_{\dot{\check{\mathbf{a}}}}\phi(\dot{\mathbf{a}},\nabla{\dot{{\mathbf{a}}}};\mathbf{q})\cdot\dot{\check{\mathbf{a}}}+\partial_{\nabla\dot{\check{\mathbf{a}}}}\phi(\dot{\mathbf{a}},\nabla{\dot{\mathbf{a}}};\mathbf{q})\cdot\nabla\dot{\check{\mathbf{a}}}\geq0.
%\label{irr0}
%\end{equation}
%\\
%Assuming that $\partial_{\dot{\check{\mathbf{a}}}}\phi(\dot{{\mathbf{a}}},\nabla{\dot{{\mathbf{a}}}};\mathbf{q})\geq\bm{0}$ always holds and that $\partial_{\nabla\dot{\check{\mathbf{a}}}}\phi(\dot{{\mathbf{a}}},\nabla{\dot{{\mathbf{a}}}};\mathbf{q})$ vanishes at the boundary, equation~\eqref{irr0} yields the irreversibility condition
%
%\begin{equation}
%\dot{\check{\mathbf{a}}}(\bm{x},t) \geq \bm{0}  \quad \text{in} \quad \Omega\times\mathrm{T},
%\label{IR}
%\tag{IR}
%\end{equation}
%\\
%which enforces the dissipation inequality in irreversible processes.

\paragraph*{Local formulation}\label{gov_Eq_ener}\
To recover the local governing equations~\eqref{st_adm} and~\eqref{biot}, we depart from the notion of global stability and admit, as candidate solutions, those satisfying a local-directional stability condition. For this purpose, let~$h\in\mathbb{R}$ denote a variation parameter and ${\tilde{\mathbf{q}}\in\tilde{\mathscr{Q}}}$ be test directions in which admissible variations on $\mathbf{q}$ take place. Then, a process $\mathbf{q}\colon\mathrm{T}\to\mathscr{Q}$ is locally stable if for all $t\in \mathrm{T}$, there exists $\bar{h}>0$ such that %In particular, we account for irreversible variables, so that ${\tilde{\mathscr{A}}\coloneqq\{\mathbf{a}\in\mathscr{A} \, | \, \check{\mathbf{a}}\geq 0 \}}$.
\begin{equation}
\mathcal{E}\big(\mathbf{q}(t)\big) \leq \mathcal{E}\big(\mathbf{q}(t) + h\tilde{\mathbf{q}}\big) - \mathcal{T}\big(t,\bm{u}(t),\bm{u}(t)+h\tilde{\bm{u}}\big)+ \mathcal{D}\big(\mathbf{q}(t),\mathbf{q}(t)+h\tilde{\mathbf{q}}\big),\quad \forall \, \tilde{\mathbf{q}}=\{\tilde{\bm{u}},\tilde{\mathbf{a}}\}\in\tilde{\mathscr{Q}}, \quad \forall \, h\in [0,\bar{h}].
%\tag{DS}
\label{stab}
\end{equation}
In the context of local stability, it seems reasonable to consider variations in a small neighborhood of the current state, and to restrict the test directions to monotonic radial (straight) paths \citep{alessi2015shape}. In this case, the dissipation distance~\eqref{distance0} specializes to
%, due to fatigue conditions, 
%
\begin{equation}
\mathcal{D}(\mathbf{q}_0,\mathbf{q}_1)=\bigg\{\int_0^1\mathcal{R}\big(\dot{\mathbf{a}}(s);\mathbf{q}(s)\big) \, \mathrm{d}s \,\, \big\vert  \,\, \mathbf{q}(s)=\mathbf{q}_0+s(\mathbf{q}_1-\mathbf{q}_0) \bigg\}.
\label{distance1}
\end{equation} 
%\begin{equation}
%\mathcal{D}(\mathbf{q}_0,\mathbf{q}_1)=\int_{\Omega}\bigg\{\int_0^1\phi\big(\dot{\mathbf{a}}(\bm{x},s),\nabla{\dot{\mathbf{a}}}(\bm{x},s);\mathbf{c}(\bm{x},s)\big) \, \mathrm{d}s \, \vert  \, \mathbf{q}(s)=\mathbf{q}_0+s(\mathbf{q}_1-\mathbf{q}_0) \bigg\}\, \mathrm{d}\bm{x}.
%\label{distance1}
%\end{equation} 
%
Assuming that the functionals are Gâteaux-differentiable up to a certain order, a differential stability condition can be defined by replacing the right-hand side of the directional stability condition~\eqref{stab} by a Taylor expansion, and enforcing necessary and sufficient conditions of increasing order for inequality~\eqref{stab} to hold true \cite{pham2011,alessi2015shape}. In this sense,  a first-order condition is obtained by the functional derivatives
\begin{equation}
\frac{\mathrm{d}}{\mathrm{d}h}\mathcal{E}\big(\mathbf{q}+h\tilde{\mathbf{q}}\big)\bigg\vert_{h=0}- \frac{\mathrm{d}}{\mathrm{d}h}\mathcal{T}\big(t,\bm{u},\bm{u}+h\tilde{\bm{u}}\big)\bigg\vert_{h=0}+\frac{\mathrm{d}}{\mathrm{d}h}\mathcal{D}\big(\mathbf{q},\mathbf{q}+h\tilde{\mathbf{q}}\big)\bigg\vert_{h=0}\geq0, \quad \forall\,\tilde{\mathbf{q}}\in\tilde{\mathscr{Q}}.
\label{DS-1_0}
\end{equation}
As in equation~\eqref{DS-1_0}, for notational simplicity, implicit dependence on time (and space) will be omitted hereafter, unless required for clarity purposes. Using the 1-homogeneity of the dissipation potential, we evaluate the Gâteaux derivative of the dissipation distance~as 
\begin{equation}
\begin{aligned}
\frac{\mathrm{d}}{\mathrm{d}h}\mathcal{D}(\mathbf{q},\mathbf{q}+h\tilde{\mathbf{q}})\bigg\vert_{h=0}&=\lim_{h\to 0}\frac{\mathcal{D}(\mathbf{q},\mathbf{q}+h\tilde{\mathbf{q}})}{h}\\
&=\lim_{h\to 0}\frac{1}{h}\int_{\Omega}\bigg\{\int_0^1\phi\big(\dot{\mathbf{a}}(\bm{x},s),\nabla{\dot{\mathbf{a}}}(\bm{x},s);\mathbf{c}(\bm{x},s)\big) \, \mathrm{d}s \ \vert  \ \mathbf{q}(s)=\mathbf{q}+hs\tilde{\mathbf{q}}  \bigg\}\, \mathrm{d}\bm{x}\\
&=\lim_{h\to 0}\frac{1}{h}\int_{\Omega}\bigg\{\int_0^h\phi\big(\dot{\mathbf{a}}(\bm{x},r),\nabla{\dot{\mathbf{a}}}(\bm{x},r);\mathbf{c}(\bm{x},r)\big) \, \mathrm{d}r \ \vert  \ \mathbf{q}(r)=\mathbf{q}+r\tilde{\mathbf{q}}  \bigg\}\, \mathrm{d}\bm{x} \\ 
&=\lim_{h\to 0}\frac{1}{h}\int_{\Omega}\int_0^h\phi\big(\tilde{\mathbf{a}},\nabla{\tilde{\mathbf{a}}};\mathbf{c}+r\tilde{\mathbf{c}}\big) \, \mathrm{d}r\, \mathrm{d}\bm{x}\\
&=\int_{\Omega}\phi\big(\tilde{\mathbf{a}},\nabla{\tilde{\mathbf{a}}};\mathbf{c}\big) \, \mathrm{d}\bm{x},\end{aligned}
\end{equation}%\\&=\mathcal{R}(\tilde{\mathbf{a}};\mathbf{q}). 
where the change of variable $r=hs$ has been used. Using this result in equation~\eqref{DS-1_0} leads to the first-order stability condition 
\begin{equation}
\delta\mathcal{E}(\mathbf{q})(\tilde{\mathbf{q}})+\mathcal{R}(\tilde{\mathbf{a}};\mathbf{q})-\int_\Omega\rho\bm{b}\cdot\tilde{\bm{u}}\,\mathrm{d}\bm{x} - \int_{\Gamma_\mathrm{N}}\bar{\bm{t}}\cdot\tilde{\bm{u}}\,\mathrm{d}S \geq0, \quad \forall\,\tilde{\mathbf{q}}\in\tilde{\mathscr{Q}},
\label{DS-1}
%\tag{DS-1}
\end{equation}
where $\delta\mathcal{E}(\mathbf{q})(\tilde{\mathbf{q}})$ is the Gâteaux derivative of $\mathcal{E}$ in the direction $\tilde{\mathbf{q}}$. Equation~\eqref{DS-1} is a necessary condition for directional stability~\eqref{stab} that becomes sufficient if the inequality is strict. Otherwise, the study of higher-order conditions is required  \citep{pham2011b,alessi2015shape,marigo2016,alessi2016shape}, which is out of scope in the present study. Furthermore, if the energy quantities in the energy balance~\eqref{EB} are sufficiently regular in $[0,t]$, the time derivative of~\eqref{EB} yields the power balance equation
\begin{equation}
\frac{\mathrm{d}}{\mathrm{d}t}{\mathcal{E}}(\mathbf{q}) + \mathcal{R}(\dot{\mathbf{a}};\mathbf{q})
 - \int_\Omega\rho\bm{b}\cdot\dot{\bm{u}}\,\mathrm{d}\bm{x} - \int_{\Gamma_\mathrm{N}}\bar{\bm{t}}\cdot\dot{\bm{u}}\,\mathrm{d}S-\int_{\Gamma_{\mathrm{D}}}\bm{t}_\mathrm{r}\cdot\dot{\bar{\bm{u}}}\, \mathrm{d}S=0.
\label{EB-1}
%\tag{EB-1}
\end{equation}

\paragraph*{Dissipation inequality}
To ensure physical soundness, a dissipation inequality is included in the formulation \cite{Alessi2016}, ensuring the fulfilment of the second law of thermodynamics:
\begin{equation}
\phi(\dot{\mathbf{a}},\nabla\dot{\mathbf{a}};\mathbf{c})\geq0 \quad \text{in} \quad  \Omega\times\mathrm{T}.
\label{DI}
%\tag{DI}
\end{equation}

\paragraph*{Governing equations}  Using equations~\eqref{stored_ext0}\textendash\eqref{dissipated}, and in view of the stress-strain relation~\eqref{constitutive}, the first-order stability condition~\eqref{DS-1} yields 
\begin{equation}
\begin{aligned}
\int_{\Omega}\big\{\delta_{{\mathbf{a}}}\psi(\bm{\varepsilon},\mathbf{a},\nabla\mathbf{a}) + \delta_{\tilde{\mathbf{a}}}\phi(\tilde{\mathbf{a}},\nabla\tilde{\mathbf{a}};\mathbf{c})\big\}\cdot\tilde{\mathbf{a}}\, \mathrm{d}\bm{x} +\int_{\Gamma}\big\{\bm{n}\cdot\big(\partial_{\nabla\tilde{\mathbf{a}}}\phi(\tilde{\mathbf{a}},\nabla\tilde{\mathbf{a}};\mathbf{c})+\partial_{\nabla\mathbf{a}}\psi(\bm{\varepsilon},\mathbf{a},\nabla\mathbf{a})\big)\big\}\cdot\tilde{\mathbf{a}}\, \mathrm{d}S&\\
-\int_{\Omega}(\mathrm{div}\bm{\sigma}+\rho\bm{b})\cdot\tilde{\bm{u}}\, \mathrm{d}\bm{x}+\int_{\Gamma_\mathrm{N}}(\bm{\sigma}\cdot\bm{n}-\bar{\bm{t}}) \cdot\tilde{\bm{u}}\, \mathrm{d}S \geq0&,
\end{aligned}
\label{energ_stab}
\end{equation} 
which holds for all $\tilde{\mathbf{q}}\in\tilde{\mathscr{Q}}$. Likewise, the power balance~\eqref{EB-1} gives 
\begin{equation}
\begin{aligned}
\int_{\Omega}\big\{\delta_{{\mathbf{a}}}\psi(\bm{\varepsilon},\mathbf{a},\nabla\mathbf{a}) + \delta_{\dot{\mathbf{a}}}\phi(\dot{\mathbf{a}},\nabla\dot{\mathbf{a}};\mathbf{c})\big\}\cdot\dot{\mathbf{a}}\, \mathrm{d}\bm{x} +\int_{\Gamma}\big\{\bm{n}\cdot\big(\partial_{\nabla\dot{\mathbf{a}}}\phi(\dot{\mathbf{a}},\nabla\dot{\mathbf{a}};\mathbf{c})+\partial_{\nabla\mathbf{a}}\psi(\bm{\varepsilon},\mathbf{a},\nabla\mathbf{a})\big)\big\}\cdot\dot{\mathbf{a}}\, \mathrm{d}S\\
  -\int_{\Omega}(\mathrm{div}\bm{\sigma}+\rho\bm{b})\cdot\dot{\bm{u}}\, \mathrm{d}\bm{x}+\int_{\Gamma_{\mathrm{N}}}(\bm{\sigma}\cdot\bm{n}-\bm{\bar{t}}\,)\cdot\dot{\bm{u}}\, \mathrm{d}S+\int_{\Gamma_{\mathrm{D}}}(\bm{\sigma}\cdot\bm{n}-\bm{t}_\mathrm{r})\cdot\dot{\bar{\bm{u}}}\, \mathrm{d}S=0.
  \end{aligned}
\label{energ_energ}
\end{equation} 
Equations~\eqref{energ_stab} and~\eqref{energ_energ} yield:  
\begin{equation}
\begin{dcases}
\mathrm{div}\bm{\sigma}+\rho\bm{b} = \bm{0} & \text{in} \quad \Omega,\\
\bm{\sigma}\cdot\bm{n}=\bar{\bm{t}} & \text{on} \quad \Gamma_\mathrm{N},\\
\bm{\sigma}\cdot\bm{n}=\bm{t}_\mathrm{r} \quad &\text{on} \quad \Gamma_\mathrm{D},
\end{dcases} \quad \text{and} \quad \begin{dcases}
-\delta_{{\mathbf{a}}}\psi(\bm{\varepsilon},\mathbf{a},\nabla\mathbf{a}) - \delta_{\tilde{\mathbf{a}}}\phi(\tilde{\mathbf{a}},\nabla\tilde{\mathbf{a}};\mathbf{c})\leq \bm{0} & \text{in} \quad \Omega,\\
\big\{-\delta_{{\mathbf{a}}}\psi(\bm{\varepsilon},\mathbf{a},\nabla\mathbf{a}) - \delta_{\dot{\mathbf{a}}}\phi(\dot{\mathbf{a}},\nabla\dot{\mathbf{a}};\mathbf{c})\big\}\cdot\dot{\mathbf{a}}= 0 & \text{in} \quad \Omega,\\
\nabla\mathbf{a}\cdot\bm{n}=\bm{0}  \quad \text{on} \quad \Gamma,
\end{dcases}
\label{energ_eq0}
\end{equation}%\bm{u}=\bar{\bm{u}} \implies
encompassing the equilibrium equations~\eqref{st_adm} and the evolution equations~\eqref{biot} for the internal variables. Note that boundary conditions for internal variables are also recovered.

% , along with the dissipation inequality~\eqref{DI},

\subsection{Incremental minimization problem}\label{incremental}

%The energetic formulation naturally gives rise to an incremental variational problem, which yields a sequence of minimizers that approximate the solution at discrete pseudo-times. From a mathematical standpoint, this incremental minimization problem serves as the basis for existence and uniqueness results, and provides tools for error estimation. On the other hand, it allows for attractive numerical implementations based on optimization techniques. [move to introduction?]

An incremental minimization problem suitable for numerical implementation naturally follows from the global energetic formulation~\citep{Mielke2015}. Consider $n_{\mathrm{t}}+1$ discrete time instants $0=t_0<\dots<t_{n}<t_{n+1}<\dots<t_{n_\mathrm{t}}=t_\mathrm{max}$. Using equation~\eqref{distance1}, the incremental dissipated energy up to $t_{n+1}$ reads
\begin{equation}
\mathcal{D}(\mathbf{q}_0,\mathbf{q}_{n+1})=\mathcal{D}(\mathbf{q}_0,\mathbf{q}_{n})+\mathcal{D}(\mathbf{q}_n,\mathbf{q}_{n+1}).
\label{inc_diss00}
\end{equation}
The goal is to evaluate the dissipated energy increment from $t_n$ to $t_{n+1}$, which follows from~\eqref{dissipenerg} as
\begin{equation}
\mathcal{D}(\mathbf{q}_n,\mathbf{q}_{n+1})=\int_{\Omega}\int_{t_{n}}^{t_{n+1}}\phi\bigg(\frac{\mathbf{a}_{n+1}-\mathbf{a}_n}{t_{n+1}-t_n},\frac{\nabla\mathbf{a}_{n+1}-\nabla\mathbf{a}_n}{t_{n+1}-t_n};\mathbf{c}_n+\frac{t-t_n}{t_{n+1}-t_n}(\mathbf{c}_{n+1}-\mathbf{c}_n)\bigg)\,\mathrm{d}t\,\mathrm{d}\bm{x}.
\label{inc_diss0}
\end{equation}

In some cases, the dissipation potential is such that the dissipated energy is a state function. For this to occur, it is necessary (but not sufficient) that $\mathcal{D}(\mathbf{q}_n,\mathbf{q}_{n+1})\equiv \mathcal{D}(\mathbf{a}_n,\mathbf{a}_{n+1})$. In other words, the dissipation potential must not depend on the non-dissipative state variables. Then, the evaluation of~\eqref{inc_diss0} is straightforward and equivalent to the direct computation of the dissipated energy $\int_{t_n}^{t_{n+1}}\mathcal{R}\big(\dot{\mathbf{a}}(\tau);\mathbf{a}(\tau)\big)\mathrm{d}\tau$. In the model presented in section~\ref{model}, we shall see that only part of the dissipated energy is a state function.  
%\mathcal{R}(\dot{\mathbf{a}};\mathbf{q})\equiv\mathcal{R}(\dot{\mathbf{a}};\mathbf{a}) \implies 

In the more general case where the dissipated energy is not a state function, as is often the case for state-dependent dissipation potentials (e.g., in non-associative plasticity), an incremental approximation is required to evaluate the time integral in~\eqref{inc_diss0}. Taking the zeroth-order term of a Taylor expansion of the integral kernel (noting that all other terms imply rate-dependence), equation~\eqref{inc_diss0} is approximated as
%$\phi(\mathbf{a}_{n+1}-\mathbf{a}_n,\nabla\mathbf{a}_{n+1}-\nabla\mathbf{a}_n;\mathbf{c}_n)$
%% HERE: Not sure about ;q_n!
\begin{equation}
\mathcal{D}(\mathbf{q}_n,\mathbf{q}_{n+1})\approx  \mathcal{R}(\mathbf{a}_{n+1}-\mathbf{a}_n;\mathbf{q}_n)\eqqcolon\mathcal{D}(\mathbf{a}_n,\mathbf{a}_{n+1};\mathbf{q}_n).
\label{disspotapprox}
\end{equation} 
This expression represents an explicit approximation of the state-dependent dissipated energy, in agreement with the incremental framework of~\citet{miehe2011} (see also~\citet{mielke2007} and~\citet{luege2018}).\footnote{Other approximations involving implicit schemes are also possible; see~\citet{stainier2011} and~\citet{brassart2012}.}

Given all states up to $t_n$, the unknown state $\mathbf{q}_{n+1}$ is found from the variational principle
\begin{equation}
\inf_{\mathbf{q}_{n+1}\in \mathscr{Q}}\big\{\mathcal{E}(\mathbf{q}_{n+1})-\mathcal{T}(\bm{u}_n,\bm{u}_{n+1})+\mathcal{D}(\mathbf{a}_n,\mathbf{a}_{n+1};\mathbf{q}_n)\big\},
\label{discminprim}
\end{equation}
where the dissipation inequality~\eqref{DI} must be enforced in incremental form.

%In some cases, the dissipation potential is such that the dissipated energy is a state function. For this to occur, it is necessary (but not sufficient) that $\mathcal{R}(\dot{\mathbf{a}};\mathbf{q})\equiv\mathcal{R}(\dot{\mathbf{a}};\mathbf{a})$, i.e., the dissipation potential must not depend on the non-dissipative state variables. Then, the evaluation of~\eqref{inc_diss_indep} is straightforward, and we shall not resort to the incremental approximation~\eqref{disspotapprox}. Clearly, the increment $\int_{t_n}^{t_{n+1}}\mathcal{R}\big(\dot{\mathbf{a}}(\tau);\mathbf{a}(\tau)\big)\mathrm{d}\tau$ then replaces the term $\mathcal{R}(\mathbf{a}_{n+1}-\mathbf{a}_{n};\mathbf{q}_n)$ in the minimization problem~\eqref{discminprim}. In the model presented in section~\ref{model}, only part of the dissipated energy is a state function.

The incremental minimization problem~\eqref{discminprim} directly recovers the continuous evolution problem~\eqref{energ_eq0} when the dissipated energy is a state function~\citep{tanne2017,samaniego2020}. In particular, stationarity of the total energy at $t_{n+1}$ implies the first-order stability condition~\eqref{DS-1},  while taking the limit as $t_{n+1}-t_n\to 0$ allows to recover the energy balance~\eqref{EB-1}. Similarly, one can show the link between the continuous and incremental problems for the general case in which the dissipated energy is not a state function~\citep{luege2018}. As will become clear in section~\ref{model}, this link is not as straightforward as for the path-independent counterpart; nevertheless, it suffices to ensure that the incremental dissipation~\eqref{disspotapprox} is such that the Euler-Lagrange equations of~\eqref{discminprim} recover the continuous evolution problem~\eqref{energ_eq0} as~$t_{n+1}-t_n\to 0$ (see~\citet{stainier2011} for a deeper discussion on this topic). 

\subsection{Overview}

Table~\ref{overview1} presents an overview of the present variational formulation. In the following section, the abstract setting presented so far will take specific form after defining the state variables, the free energy density and the dissipation potential.
 %As shown in the model presented in the following section, the framework can be directly applied after defining the state variables and the basic energy quantities. 

\begin{table}[h!]
\renewcommand{\arraystretch}{1.1}
\small
\centering
\caption{Constitutive energy functionals, energetic formulation and incremental energy minimization.}
\begin{tabular}{ll}
\toprule
\multicolumn{2}{l}{\sffamily\textbf{Constitutive energy functionals}}
\\ 
\midrule
Stored energy            & $\mathcal{E}(\mathbf{q})=\int_\Omega \psi(\bm{\varepsilon},\mathbf{a},\nabla\mathbf{a})\,\mathrm{d}\bm{x}$ 
\\ \grayrule
Dissipative power        & $\mathcal{R}(\dot{\mathbf{a}};\mathbf{q})=\int_\Omega\phi(\dot{\mathbf{a}},\nabla\dot{\mathbf{a}};\mathbf{c})\,\mathrm{d}\bm{x}$,\; with \;$\phi(\dot{\mathbf{a}},\nabla\dot{\mathbf{a}};\mathbf{c})\geq0$ 
\\ \grayrule
Incremental dissipated energy  &                                                                                   
$\mathcal{D}(\mathbf{a}_{n+1},\mathbf{a}_n;\mathbf{q}_n)$, \; from~\eqref{inc_diss0} or~\eqref{disspotapprox}
\\
\midrule
\multicolumn{2}{l}{\sffamily\textbf{Continuous evolution problem}}                                                                                                                                                       
\\ \hline
Energy balance~\eqref{EB-1}       &   $\frac{\mathrm{d}}{\mathrm{d}t}{\mathcal{E}}(\mathbf{q}) + \mathcal{R}(\dot{\mathbf{a}};\mathbf{q})
 - \int_\Omega\rho\bm{b}\cdot\dot{\bm{u}}\,\mathrm{d}\bm{x} - \int_{\Gamma_\mathrm{N}}\bar{\bm{t}}\cdot\dot{\bm{u}}\,\mathrm{d}S-\int_{\Gamma_{\mathrm{D}}}\bm{t}_\mathrm{r}\cdot\dot{\bar{\bm{u}}}\, \mathrm{d}S=0$                                                                                
\\ \grayrule
Stability~\eqref{DS-1}            &    $\delta\mathcal{E}(\mathbf{q})(\tilde{\mathbf{q}})+\mathcal{R}(\tilde{\mathbf{a}};\mathbf{q})-\int_\Omega\rho\bm{b}\cdot\tilde{\bm{u}}\,\mathrm{d}\bm{x} - \int_{\Gamma_\mathrm{N}}\bar{\bm{t}}\cdot\tilde{\bm{u}}\,\mathrm{d}S \geq0$                                                                                    \\ \midrule
\multicolumn{2}{l}{\sffamily \textbf{Incremental evolution problem}}
\\ \hline
Incremental minimization~\eqref{discminprim}  & $\inf_{\mathbf{q}_{n+1}\in \mathscr{Q}}\big\{\mathcal{E}(\mathbf{q}_{n+1})-\mathcal{T}(\bm{u}_n,\bm{u}_{n+1})+\mathcal{D}(\mathbf{a}_n,\mathbf{a}_{n+1};\mathbf{q}_n)\big\}$ 
\\ \bottomrule
\end{tabular}
\label{overview1}
\end{table}

%\begin{table}[]
%\begin{tabular}{lll}
%\hline
%%\multicolumn{2}{c}{State variables}                                                          \\ \hline
%%Primary               & 
%%$\mathbf{q}\coloneqq\{\bm{u},\mathbf{a}\}$                            \\
%%Constitutive state    & 
%%$\mathbf{c}\coloneqq\{\bm{\varepsilon},\mathbf{a},\nabla\mathbf{a}\}$ \\ \hline
%%\multicolumn{2}{c}{Local energy quantities}                                                  \\ \hline
%%Free energy           &                                                                       $\psi(\bm{\varepsilon},\mathbf{a},\nabla\mathbf{a})$\\
%%Dissipation potential &                                                                      $\phi(\dot{\mathbf{a}},\nabla{\dot{\mathbf{a}}};\mathbf{c})$\\ \hline
%\multicolumn{3}{c}{Energy functionals}                                                       \\ \hline
%Continuous            &                                                                     
%&
%\\
%Incremental           &                                                                     
%&
%\\ \hline
%\multicolumn{3}{c}{Evolution problem}                                                        \\ \hline
%Continuous            &                                                                     
%&
%\\
%Incremental           &                                                                     
%&
%\\ \hline
%\end{tabular}
%\end{table}

\section{A phase-field model with fatigue effects coupled to cyclic plasticity}\label{model}

This section presents the proposed model that couples a phase-field approach to fatigue with cyclic plasticity, derived using the concepts presented in section~\ref{variat_general}. First, the modeling of cyclic plasticity is addressed, where a non-associative ratcheting variable is introduced in the variational framework. Then, the full model is elaborated by introducing the phase-field description of fatigue-induced fracture. 

\subsection{Modeling cyclic plasticity with ratcheting}\label{model:plast}% Discuss these concepts in an appendix?

To account for cyclic plasticity with ratcheting effects, we take the model of \citet{houlsby2017} as a point of departure, which was proposed to describe the behavior of cyclically loaded pile foundations \citep{abadie2015}. The original model is tightly linked to~\citet{armstrong1966}, but allows for multi-surface kinematic hardening with a single ratcheting variable. Herein, we also consider isotropic hardening/softening and gradient-extended plasticity~\citep{muhlhaus1991,DeBorstMull1992,fleck2009a,gurtin2009,nguyen2016}, accounting for a general class of elastoplastic materials.  %\citep{muhlhaus1991,DeBorstMull1992,gurtin2009,nguyen2016}

To present the model, we adopt the conventional setting of plasticity, where auxiliary internal variables are introduced to account for hardening effects~\citep{maugin1992,simo1998} and, in this case, ratcheting. As in classical elastoplasticity, the evolution of the internal variables is first presented in dual form, that is, in terms of yield functions in generalized-stress space, where the issue of associativity is addressed. Then, a state-dependent primal dissipation potential is obtained, which is used to construct the coupled fatigue model in section~\ref{fatiguemodel}.

%Herein, we elaborate the multi-dimensional model in a variational setting. 
%As previously mentioned, gradient-enhanced models offer several features that are key to model strain localization.

\subsubsection{Constitutive model}
We focus on a multi-surface representation of the plastic deformation process, with $n_\mathrm{y}\geq1$ yield surfaces of increasing yield strength. To this end, we introduce  the internal variables  
\begin{equation}
\{\bm{\varepsilon}^\mathrm{p},\bm{\kappa},\bm{\varepsilon}^\mathrm{r}\}, \quad \text{with} \quad \bm{\varepsilon}^\mathrm{p}\coloneqq\big\{\bm{\varepsilon}^{\mathrm{p}}_1,\dots,\bm{\varepsilon}^{\mathrm{p}}_s,\dots,\bm{\varepsilon}^{\mathrm{p}}_{n_\mathrm{y}}\big\} \quad \text{and} \quad \bm{\kappa}\coloneqq\big\{\kappa_1,\dots,\kappa_s,\dots,\kappa_{n_\mathrm{y}}\big\}.
\label{plast_intvar}
\end{equation}
The subscript $s$ denotes the $s^\text{th}$ yield surface with a corresponding plastic strain tensor $\bm{\varepsilon}^\mathrm{p}_{s}\colon \Omega\times\mathrm{T}\to\mathbb{R}^{d\times d}_{\mathrm{dev}}\coloneqq\{\bm{e}\in\mathbb{R}^{d\times d}_\mathrm{sym}\,|\,\tr(\bm{e})=0\}$ and an isotropic hardening/softening variable $\kappa_s\colon \Omega\times\mathrm{T}\to\mathbb{R}_+$, and ${\bm{\varepsilon}^{\mathrm{r}}\colon \Omega\times\mathrm{T}\to\mathbb{R}^{d\times d}_{\mathrm{dev}}}$ is a ratcheting strain tensor. As a modeling assumption, $\kappa_s$ evolves according to the hardening law
\begin{equation}
\dot{\kappa}_s\coloneqq\sqrt\frac{2}{3}\Vert\dot{\bm{\varepsilon}}^{\mathrm{p}}_s\Vert.
\label{hardrate} 
\end{equation}
 The ratcheting strain tensor $\bm{\varepsilon}^\mathrm{r}$ evolves according to the ratcheting law
\begin{equation}
\dot{\bm{\varepsilon}}^{\mathrm{r}}\coloneqq\beta\hat{\bm{r}}\sum_{s=1}^{n_{\mathrm{y}}}\Vert\dot{\bm{\varepsilon}}^{\mathrm{p}}_s\Vert,
\label{ratchrate} 
\end{equation}
where $\hat{\bm{r}}\coloneqq\partial\Vert\bm{\sigma}_\mathrm{dev}\Vert$ is the direction of the deviatoric stress~\citep{houlsby2017}. The parameter $\beta\in[0,1]$ defines the fraction of  plastic strains that contributes to ratcheting, assumed, for simplicity, to be equal for all yield~surfaces.

In light of the internal variables defined above, we shall present the plasticity model in the multi-surface setting, in agreement with the ratcheting model of~\citet{houlsby2017}. This choice has been made to favor generality and to endow the model with great flexibility to accommodate different material responses. Note that throughout the sequel, the single-surface case can always be  recovered, where~${n_\mathrm{y}=1}$. 

%, and $\bm{\sigma}_\mathrm{dev}$
%denotes the deviatoric part of the stress tensor.  
%From~\eqref{ratchrate} , it follows that the plastic strains can be taken as primary fields introduced in Section~\ref{variat_general}:
%
%\begin{equation}
%\mathbf{a}\coloneqq \big\{\mathbf{z},\bm{\kappa},\bm{\varepsilon}^{\mathrm{r}}\big\} \quad \text{and} \quad \mathbf{z}\coloneqq\bm{\varepsilon}^\mathrm{p}.
%\label{internal_model}
%\end{equation}
%\\

The free energy density is defined as
%\bigg(\bm{\varepsilon}-\sum_{s=1}^{n_\mathrm{y}}\bm{\varepsilon}^\mathrm{p}_{s}-\bm{\varepsilon}^r\big):\mathbf{C}:\big(\bm{\varepsilon}-\sum_{s=1}^{n_\mathrm{y}}\bm{\varepsilon}^\mathrm{p}_{s}-\bm{\varepsilon}^r\big)
\begin{equation}
\psi(\bm{\varepsilon},\bm{\varepsilon}^\mathrm{p},\bm{\kappa},\bm{\varepsilon}^\mathrm{r},\nabla\bm{\kappa})\coloneqq\underbrace{\frac{1}{2}\Vert\bm{\varepsilon}-\sum_{s=1}^{n_\mathrm{y}}\bm{\varepsilon}^\mathrm{p}_{s}-\bm{\varepsilon}^\mathrm{r}\big\Vert^2_{\textstyle\bm{\mathsf{C}}}}_{\textstyle \psi^{\mathrm{e}}(\bm{\varepsilon},\bm{\varepsilon}^\mathrm{p},\bm{\varepsilon}^\mathrm{r})}+\underbrace{\frac{1}{2}\sum\limits_{s=1}^{n_\mathrm{y}}\big(H^{\mathrm{kin}}_s\bm{\varepsilon}^\mathrm{p}_s:\bm{\varepsilon}^\mathrm{p}_s + H^{\mathrm{iso}}_s\kappa_s^2\big)+\frac{1}{2}\sum\limits_{s=1}^{n_\mathrm{y}}\eta_{\mathrm{p} s}^2\Vert\nabla\kappa_s\Vert^2}_{\textstyle \psi^{\mathrm{p}}(\bm{\varepsilon}^\mathrm{p},\bm{\kappa},\nabla\bm{\kappa})},
\label{free_ener_model}
\end{equation}
where $\bm{\mathsf{C}}$ is the fourth-order elastic tensor. $H^{\mathrm{kin}}_s$ and $H^{\mathrm{iso}}_s$ denote the kinematic and isotropic hardening moduli, respectively, for the $s^\mathrm{th}$ yield surface. The last term in equation~\eqref{free_ener_model} introduces non-local effects governed by the plastic internal length scale $\eta_{\mathrm{p}s}$, which is related to the plastic characteristic length $\ell_{\mathrm{p}s}$ by
\begin{equation}
\eta_{\mathrm{p}s}=\ell_{\mathrm{p} s}\sqrt{\sigma^\mathrm{p}_s},
\end{equation}
where $\sigma^{\mathrm{p}}_s$ is the $s^\text{th}$ plastic yield strength (cf.~\citep{MieTeiAld2016} for the single-surface case). Note that in contrast with $\eta_{\mathrm{p}s}$,  $\ell_{\mathrm{p}s}$ has the units of length.

In view of equation~\eqref{free_ener_model}, the stress tensor is obtained from equation~\eqref{constitutive} as
\begin{equation}
\bm{\sigma}(\bm{\varepsilon},\bm{\varepsilon}^\mathrm{p},\bm{\varepsilon}^\mathrm{r})=\frac{\partial \psi}{\partial \bm{\varepsilon}}=\bm{\mathsf{C}}:\bigg(\bm{\varepsilon}-\sum_{s=1}^{n_\mathrm{y}}\bm{\varepsilon}^\mathrm{p}_{s}-\bm{\varepsilon}^\mathrm{r}\bigg).
\label{constitutive_model}
\end{equation} 
The generalized stresses conjugate to~\eqref{plast_intvar} read 
\begin{equation}
\{\bm{s}^\mathrm{p},\bm{h},\bm{s}^\mathrm{r}\}, \quad \text{with} \quad \bm{s}^\mathrm{p}\coloneqq\{\bm{s}^\mathrm{p}_1,\dots,\bm{s}^\mathrm{p}_s,\dots,\bm{s}^\mathrm{p}_{n_\mathrm{y}}\} \quad \text{and} \quad \bm{h}\coloneqq\{h_1,\dots,h_s,\dots, h_{n_\mathrm{y}}\}.
\end{equation}
These dual variables follow from equation~\eqref{dualforces} as 
\begin{equation}
\bm{s}^\mathrm{p}_s =-\delta_{\bm{\varepsilon}^\mathrm{p}_s}\psi\equiv\bm{\sigma} - H^\mathrm{kin}_s\bm{\varepsilon}^\mathrm{p}_s, \quad h_s =-\delta_{\kappa_s}\psi\equiv- H^{\mathrm{iso}}_s\kappa_s + \eta_{\mathrm{p}s}^2\mathrm{div}[\nabla\kappa_s], \quad \bm{s}^\mathrm{r} =-\delta_{\bm{\varepsilon}^\mathrm{r}}\psi\equiv\bm{\sigma}.
\label{constitutive_model2}
\end{equation}

\subsubsection{Dissipation}

Focusing on $J_2$ plasticity, the yield functions are defined in generalized-stress space~as
\begin{equation}
f_{s}^\mathrm{p}(\bm{s}^\mathrm{p}_s ,h_s)\coloneqq\big\Vert\bm{s}^\mathrm{p}_{{s\,\mathrm{dev}}}\big\Vert - \sqrt\frac{2}{3}(\sigma^{\mathrm{p}}_s -  h_s) \leq0.
\label{yield_s}
\end{equation}
The set of admissible generalized stresses for each yield surface is the convex set $$\mathbb{K}_s\coloneqq\big\{\bm{s}^\mathrm{p}_s,h_s,\bm{s}^\mathrm{r} \ \vert \ f_s^\mathrm{p}(\bm{s}^\mathrm{p}_s ,h_s)\leq 0\big\},$$
such that the elastic domain is the intersection of all $\mathbb{K}_s$, and thus still a convex set:
\begin{equation}
\mathbb{K}\coloneqq\bigcap^{n_\mathrm{y}}_{s=1}\mathbb{K}_s=\big\{\bm{s}^\mathrm{p},\bm{h},\bm{s}^\mathrm{r} \ \vert \ f_s^\mathrm{p}(\bm{s}^\mathrm{p}_s ,h_s)\leq 0 \ \forall \, s\in\mathrm{Y}\big\},
\label{eldom}
\end{equation}
where $\mathrm{Y}\coloneqq\mathbb{Z}\cap[1,n_{\mathrm{y}}]$. To obtain the primal representation of the dissipation power in the case of associative models, the dissipation potential is obtained as the support function of $\mathbb{K}$ for given rates $\{\dot{\bm{\varepsilon}}^\mathrm{p},\dot{\bm{\kappa}},\dot{\bm{\varepsilon}}^\mathrm{r}\}$: 
\begin{equation}
\phi=\sup\{{\tilde{\bm{s}}^\mathrm{p}}:\dot{\bm{\varepsilon}}^\mathrm{p}+\tilde{\bm{h}}\cdot\dot{\bm{\kappa}} + {\tilde{\bm{s}}^\mathrm{r}}:\dot{\bm{\varepsilon}}^\mathrm{r} \ \vert  \ \{{\tilde{\bm{s}}^\mathrm{p}},\tilde{\bm{h}},{\tilde{\bm{s}}^\mathrm{r}}\}\in\mathbb{K}\}.
\label{maxdisp1}
\end{equation}
This optimization problem can be viewed as a particular statement of the principle of maximum dissipation. Equation~\eqref{maxdisp1} yields as necessary conditions the associative flow rule for the plastic strains, as well as the associative hardening law, consistent with equation~\eqref{hardrate} and the normality law~\citep{moreau1970,han1999}: 
\begin{equation}
\dot{{\bm{\varepsilon}}}^\mathrm{p}_s=\sum_{k=1}^{n_{\mathrm{y}}}\lambda_k\frac{\partial f_k^\mathrm{p}}{\partial {\bm{s}^\mathrm{p}_s}} \equiv \lambda_s \hat{\bm{n}}_s, \quad \dot{\kappa}_s= \sum_{k=1}^{n_{\mathrm{y}}}\lambda_k\frac{\partial f_k^\mathrm{p}}{\partial h_s}\equiv\sqrt\frac{2}{3}\lambda_s, 
\label{flowrules}
\end{equation}
where $\lambda_s\geq 0$, and $\hat{\bm{n}}_s\coloneqq\bm{s}^{\mathrm{p}}_{{s\,\mathrm{dev}}}/\Vert \bm{s}^{\mathrm{p}}_{{s\,\mathrm{dev}}}\Vert$ is the direction of the plastic flow. However, the fact that $\partial_{\bm{s}^\mathrm{r}}f_s^\mathrm{p}=\bm{0}$ implies  that $\dot{\bm{\varepsilon}}^\mathrm{r}\equiv\bm{0}$. Therefore, the ratcheting law~\eqref{ratchrate} is not obtained from~\eqref{maxdisp1} and is thus considered non-associative.
Equation~\eqref{maxdisp1} then leads to a dissipation potential that is independent of the ratcheting strain rate and is thus inconsistent with the ratcheting model.  
 
An approach to circumvent this issue and recover a variational structure in models with non-associative components is to let the set of admissible generalized stresses depend on the current state~\cite{francfort2018}. To this end, we define
the set %
\begin{equation}
\mathbb{L}(\bm{s}^\mathrm{r})\coloneqq\big\{\tilde{\bm{s}}^\mathrm{p},\tilde{\bm{h}},\tilde{\bm{s}}^\mathrm{r} \ \vert \  g_s^\mathrm{p}(\tilde{\bm{s}}^\mathrm{p}_s ,\tilde{h}_s,\tilde{\bm{s}}^\mathrm{r}) \leq  \beta\Vert\bm{s}^\mathrm{r}_\mathrm{dev}\Vert \ \forall s\in\mathrm{Y}\big\},
\label{eldom2}
\end{equation}
where $g_s^\mathrm{p}(\bm{s}^\mathrm{p}_s ,h_s,\bm{s}^\mathrm{r})\coloneqq f_s^\mathrm{p}(\bm{s}^\mathrm{p}_s ,h_s) +  \beta\Vert\bm{s}^\mathrm{r}_\mathrm{dev}\Vert$. The dissipation potential then follows as the stress-dependent support function of $\mathbb{L}(\bm{s}^\mathrm{r})$:
\begin{equation}
\phi=\sup\{\tilde{\bm{s}}^\mathrm{p}:\dot{\bm{\varepsilon}}^\mathrm{p}+\tilde{\bm{h}}\cdot\dot{\bm{\kappa}} + \tilde{\bm{s}}^\mathrm{r}:\dot{\bm{\varepsilon}}^\mathrm{r} \ \vert \ \{{\tilde{\bm{s}}^\mathrm{p}},\tilde{\bm{h}},{\tilde{\bm{s}}^\mathrm{r}}\}\in\mathbb{L}(\bm{s}^\mathrm{r})\}.
\label{maxdisp2}
\end{equation}
Equation~\eqref{maxdisp2} yields as necessary conditions the evolution equations~\eqref{flowrules} as well as the ratcheting law
\begin{equation}
\dot{\bm{\varepsilon}}^{\mathrm{r}} = \sum_{s=1}^{n_{\mathrm{y}}}\lambda_s\frac{\partial g_s}{\partial {\bm{s}^\mathrm{r}}} \equiv 
\beta\hat{\bm{r}}\sum_{s=1}^{n_{\mathrm{y}}}\lambda_s.
\label{maxdisp_new1}
\end{equation}
In equations~\eqref{flowrules} and~\eqref{maxdisp_new1}, we note that the direction of the plastic flow for each yield surface coincides with the direction of the relative deviatoric stress $\bm{s}^{\mathrm{p}}_{{s\,\mathrm{dev}}}$, while the direction of the ratcheting strain rate coincides with the direction of the deviatoric stress $\bm{\sigma}_\mathrm{dev}$. In view of equations~\eqref{hardrate} and \eqref{ratchrate}, equation~\eqref{maxdisp2} is evaluated for all $\{\dot{\bm{\varepsilon}}^\mathrm{p},\dot{\bm{\kappa}},\dot{\bm{\varepsilon}}^\mathrm{r}\} \in\mathbb{R}^{d\times d}_{\mathrm{dev}}\times\mathbb{R}_+\times \mathbb{R}^{d\times d}_{\mathrm{dev}}$ as
\begin{equation*}
\begin{aligned}
&\phi(\dot{\bm{\varepsilon}}^\mathrm{p},\dot{\bm{\kappa}},\dot{\bm{\varepsilon}}^\mathrm{r};\bm{s}^\mathrm{r})
=\sup\bigg\{\tilde{\bm{s}}^\mathrm{p}:\dot{\bm{\varepsilon}}^\mathrm{p}+\tilde{\bm{h}}\cdot\dot{\bm{\kappa}} + \tilde{\bm{s}}^\mathrm{r}:\dot{\bm{\varepsilon}}^\mathrm{r} \ \ \big\vert \ \  g_s(\tilde{\bm{s}}^\mathrm{p}_s ,\tilde{h}_s,\tilde{\bm{s}}^\mathrm{r}) \leq  \beta\Vert\bm{s}^\mathrm{r}_\mathrm{dev}\Vert \ \forall s\in\mathrm{Y}\bigg\}\\
&=\sup\bigg\{\sum_{s=1}^{n_{\mathrm{y}}} \big(\Vert\tilde{\bm{s}}^{\mathrm{p}}_{s\,\mathrm{dev}}\Vert \Vert\dot{{\bm{\varepsilon}}}^\mathrm{p}_s\Vert+\tilde{h}_s\dot{\kappa}_s\big) + \Vert\tilde{\bm{s}}^\mathrm{r}_\mathrm{dev}\Vert\Vert\dot{\bm{\varepsilon}}^\mathrm{r}\Vert   \ \ \big\vert \ \  \Vert\tilde{\bm{s}}^\mathrm{p}_{s\,\mathrm{dev}}\Vert  \leq  \beta\Vert\bm{s}^\mathrm{r}_\mathrm{dev}\Vert  - \sqrt{\tfrac{2}{3}}(\tilde{h}_s-\sigma^{\mathrm{p}}_s) -  \beta\Vert\tilde{\bm{s}}^\mathrm{r}_\mathrm{dev}\Vert \ \forall s\in\mathrm{Y}\bigg\}\\
&=\sup\bigg\{\sum_{s=1}^{n_{\mathrm{y}}} \bigg[\bigg(\beta\Vert\bm{s}^\mathrm{r}_\mathrm{dev}\Vert - \sqrt{\tfrac{2}{3}}(\tilde{h}_s-\sigma^{\mathrm{p}}_s) -  \beta\Vert\tilde{\bm{s}}^\mathrm{r}_\mathrm{dev}\Vert \bigg)\Vert\dot{{\bm{\varepsilon}}}^\mathrm{p}_s\Vert   + \tilde{h}_s\dot{\kappa}_s \bigg] + \Vert\tilde{\bm{s}}^\mathrm{r}_\mathrm{dev}\Vert\Vert\dot{{\bm{\varepsilon}}}^\mathrm{r}\Vert    \bigg\}\\
&=\sup\bigg\{\sum_{s=1}^{n_{\mathrm{y}}} \bigg[\bigg(\sqrt{\tfrac{2}{3}}\sigma^{\mathrm{p}}_s + \beta\Vert\bm{s}^\mathrm{r}_\mathrm{dev}\Vert\bigg)\Vert\dot{{\bm{\varepsilon}}}^\mathrm{p}_s\Vert   + \tilde{h}_s\bigg(\dot{\kappa}_s-\sqrt{\tfrac{2}{3}}\Vert\dot{{\bm{\varepsilon}}}^\mathrm{p}_s\Vert\bigg)\bigg]  +  \Vert\tilde{\bm{s}}^\mathrm{r}_\mathrm{dev}\Vert\bigg(\Vert\dot{{\bm{\varepsilon}}}^\mathrm{r}\Vert-  \beta  \sum_{s=1}^{n_{\mathrm{y}}}  \Vert\dot{{\bm{\varepsilon}}}^\mathrm{p}_s\Vert\bigg)    \bigg\}.
\end{aligned}
\end{equation*}%\equiv\sum_{s=1}^{n_{\mathrm{y}}}\sigma^{\mathrm{p}}_s\dot{\kappa}_s + \bm{\sigma}:\dot{\bm{\varepsilon}}^\mathrm{r} 
The expression inside the supremum becomes unbounded if the term multiplying $\Vert\tilde{\bm{s}}^\mathrm{r}_\mathrm{dev}\Vert$ is positive. Moreover, for isotropic hardening plasticity, $\tilde{h}_s\leq 0$, such that the expression inside the supremum becomes unbounded if the term multiplying $\tilde{h}_s$ is negative. The dissipation potential then takes the closed form
\begin{equation}
\phi(\dot{\bm{\varepsilon}}^\mathrm{p},\dot{\kappa},\dot{\bm{\varepsilon}}^\mathrm{r};\bm{s}^\mathrm{r})=\begin{dcases}\sum_{s=1}^{n_{\mathrm{y}}} \bigg(\sqrt\frac{2}{3}\sigma^{\mathrm{p}}_s + \beta\Vert\bm{s}^\mathrm{r}_\mathrm{dev}\Vert\bigg)\Vert\dot{{\bm{\varepsilon}}}^\mathrm{p}_s\Vert \quad \text{if} \quad \dot{\kappa}_s \geq \sqrt{\tfrac{2}{3}}\Vert\dot{\bm{\varepsilon}}^\mathrm{p}_s\Vert \quad \text{and} \quad  \Vert \dot{\bm{\varepsilon}}^\mathrm{r}\Vert\leq\beta\sum_{s=1}^{n_\mathrm{y}}\Vert \dot{\bm{\varepsilon}}^\mathrm{p}\Vert,\\
+\infty \quad \text{otherwise}.
\end{dcases} 
\label{diss_model0}
\end{equation}
A similar result follows in case of isotropic softening. Expression~\eqref{diss_model0} is consistent with the dissipation potential proposed by \citet{houlsby2017} (without isotropic hardening) in the context of hyperplasticity, a topic strongly linked to generalized standard materials~\citep{houlsby2007}.

%\subsubsection{The single-surface case}
%
%%%% HERE! Mention single-surface in 3.1.1, then again either here and/or after the damage model (at the end of the section), or in the examples?
%
%easily recovered from the multi-surface model elaborated in section~\eqref{model} by setting $n_\mathrm{y}=1$. Then, it suffices to define a single plastic strain tensor $\bm{\varepsilon}^\mathrm{p}\in\mathbb{R}_\mathrm{dev}$ and a single isotropic hardening/softening variable $\kappa\in\mathbb{R}_+$, such that the hardening and ratcheting laws~\eqref{hardrate} and~\eqref{ratchrate} become
%%
%\begin{equation}
%\dot{\kappa}=\sqrt{\frac{2}{3}}\Vert\dot{\bm{\varepsilon}}^\mathrm{p}\Vert,\quad\text{and}\quad\bm{\varepsilon}^\mathrm{r}=\beta\hat{\bm{r}}\Vert\dot{\bm{\varepsilon}}^\mathrm{p}\Vert.
%\end{equation}
%%
%The free energy and dissipation potential read
%%
%\begin{equation}
%\begin{aligned}
%\psi \\
%\phi
%\end{aligned}
%\end{equation}
%%
%while the plasticity and damage yield functions~\eqref{fp} and~\eqref{fd}) are rewritten as
%%
%\begin{equation}
%\begin{aligned}
%f^\mathrm{p}=\\
%f^\mathrm{d}=
%\end{aligned}
%\end{equation}
%%
%where $\eta_\mathrm{p}=\ell_\mathrm{p}\sqrt{\sigma^\mathrm{p}}$.

%Using \eqref{ratchrate}, we can rewrite the dissipation potential as
%%
%\begin{equation}
%\phi(\dot{\mathbf{a}},\nabla\dot{\mathbf{a}};\mathbf{q})=\sum_{s=1}^{n_{\mathrm{y}}}\big(\sigma^{\mathrm{p}}_s + \bm{\sigma}:R\hat{\bm{r}}_s\big)\dot{\kappa}_s \geq 0,
%\label{diss_model1}
%\end{equation}
%\\
%where irreversibility on $\kappa_s$ is assumed a priori. 
%%%% HEREEE!!

\subsection{Coupling cyclic plasticity to a phase-field approach to fatigue}\label{fatiguemodel}

\subsubsection{Constitutive model}

The formulation is now extended, coupling the cyclic plasticity model to a phase-field description of fatigue fracture. For this purpose, the usual phase-field/damage variable is defined as~${\alpha\colon \Omega\times\mathrm{T}\to[0,1]}$, with $\alpha=0$ and $\alpha=1$ corresponding to an undamaged and a fully degraded material state, respectively. Regularization is attained by the gradient of the phase field ${\nabla\alpha\colon \Omega\times\mathrm{T}\to\mathbb{R}^d}$. Then, along with the plastic variables introduced in the previous section, the constitutive state reads $$\mathbf{c}\coloneqq\{\bm{\varepsilon},\bm{\varepsilon}^\mathrm{p},\bm{\kappa},\bm{\varepsilon}^\mathrm{r},\alpha,\nabla\bm{\kappa},\nabla\alpha\},$$ with primary fields
$$\mathbf{q}\coloneqq\{\bm{u},\bm{\varepsilon}^\mathrm{p},\bm{\kappa},\bm{\varepsilon}^\mathrm{r},\alpha\}.$$ The damage variable is considered irreversible, thus excluding healing: 
\begin{equation}
\dot{\alpha}\geq 0  \quad \text{in} \quad \Omega\times\mathrm{T}.
\label{IR_model}
\end{equation}
Note that from the hardening law~\eqref{hardrate}, it follows that $\kappa_s\in\bm{\kappa}$ is also irreversible by definition. 

Material degradation is achieved by letting the free energy density decrease as a function of $\alpha$: 
\begin{equation}
\psi(\bm{\varepsilon},\bm{\varepsilon}^\mathrm{p},\bm{\kappa},\bm{\varepsilon}^\mathrm{r},\alpha,\nabla\bm{\kappa})\coloneqq g(\alpha)\big(\psi^{\mathrm{e}_+}(\bm{\varepsilon},\bm{\varepsilon}^\mathrm{p},\bm{\varepsilon}^\mathrm{r})+\psi^{\mathrm{p}}(\bm{\varepsilon}^\mathrm{p},\bm{\kappa},\nabla\bm{\kappa})\big)+\psi^{\mathrm{e}_-}(\bm{\varepsilon},\bm{\varepsilon}^\mathrm{p},\bm{\varepsilon}^\mathrm{r}),
\label{free_ener_model2}
\end{equation}
where $g(\alpha)$ is a damage degradation function endowed with the following properties:
\begin{equation}
g(0)=1, \quad g(1)=0, \quad g'(\alpha)\leq0 \quad \forall \ \alpha\in[0,1].
\label{degrad_prop}
\end{equation}
In this work, the widely used quadratic function is adopted:
\begin{equation}
g(\alpha)\coloneqq(1-\alpha)^2.
\label{quad_degr}
\end{equation}
In equation~\eqref{free_ener_model2}, $\psi^\mathrm{p}$ is the plastic free energy density given in equation~\eqref{free_ener_model}. On the other hand, the elastic free energy is split into damageable and undamageable parts, denoted as $\psi^{\mathrm{e}_+}$ and $\psi^{\mathrm{e}_-}$, respectively. This decomposition is introduced to consider asymmetric behavior in tension and compression, and can be performed in different ways. Herein, we consider the volumetric-deviatoric split \citep{AmorMarMau2009}:% and the spectral split \citep{MieHofWel2010}, which read respectively as
\begin{equation}
\begin{dcases}\psi^{\mathrm{e_+}}(\bm{\varepsilon},\bm{\varepsilon}^\mathrm{p},\bm{\varepsilon}^\mathrm{r})\coloneqq\frac{1}{2}K\langle\tr(\bm{\varepsilon}^\mathrm{e})\rangle_+^2+\mu(\bm{\varepsilon}^\mathrm{e}_\mathrm{dev}:\bm{\varepsilon}^\mathrm{e}_\mathrm{dev}),  \\ 
\psi^{\mathrm{e_-}}(\bm{\varepsilon},\bm{\varepsilon}^\mathrm{p},\bm{\varepsilon}^\mathrm{r})\coloneqq\frac{1}{2}K\langle\tr(\bm{\varepsilon}^\mathrm{e})\rangle_-^2,
\end{dcases} 
\label{split}
\end{equation} 
%\text{spectral} &\begin{dcases}\psi^{\mathrm{e_+}}(\bm{\varepsilon},\bm{\varepsilon}^\mathrm{p},\bm{\varepsilon}^\mathrm{r})\coloneqq\frac{1}{2}\lambda\langle\tr(\bm{\varepsilon}^\mathrm{e})\rangle_+^2+\mu(\bm{\varepsilon}^\mathrm{e}_+:\bm{\varepsilon}^\mathrm{e}_+), \\ 
%\psi^{\mathrm{e_-}}(\bm{\varepsilon},\bm{\varepsilon}^\mathrm{p},\bm{\varepsilon}^\mathrm{r})\coloneqq\frac{1}{2}\lambda\langle\tr(\bm{\varepsilon}^\mathrm{e})\rangle_-^2+\mu(\bm{\varepsilon}^\mathrm{e}_-:\bm{\varepsilon}^\mathrm{e}_-). 
%\end{dcases}
%\end{aligned}
%
where $K$ is the bulk modulus, $\mu$ is the shear modulus, and 
\begin{equation}
\bm{\varepsilon}^\mathrm{e}\coloneqq\bm{\varepsilon}-\sum_{s=1}^{n_\mathrm{y}}\bm{\varepsilon}^\mathrm{p}_s-\bm{\varepsilon}^\mathrm{r}
\label{elastic_strains}
\end{equation}
is the elastic strain tensor. From the decomposed energy density~\eqref{free_ener_model2}, the stress-strain relation~\eqref{constitutive}~reads
%Note that the same damage degradation function affects both $\psi^\mathrm{e}$ and $\psi^\mathrm{p}$, in agreement with most phase-field models for ductile fracture.
\begin{equation}
\bm{\sigma}(\bm{\varepsilon},\bm{\varepsilon}^\mathrm{p},\bm{\varepsilon}^\mathrm{r},\alpha)=\frac{\partial \psi}{\partial \bm{\varepsilon}}=g(\alpha)\frac{\partial \psi^{\mathrm{e_+}}}{\partial\bm{\varepsilon}}+\frac{\partial \psi^{\mathrm{e_-}}}{\partial\bm{\varepsilon}}.
\label{stress_model}
\end{equation}

\subsubsection{Dissipation}

The definition of the dissipation potential should entail coupling of the plastic strains to damage evolution, for which we follow the approaches in references  \citep{AleMarVid2014,AleMarVid2015,UllRodSam2016,alessi2018c,rodriguez2018}, and extend the formulation to multiple yield surfaces, kinematic hardening and ratcheting. Moreover, following \citet{Alessi2017fatigue} and \citet{carrara2020}, fatigue effects are considered by means of a local fatigue variable, defined here as ${\gamma\colon \Omega\times\mathrm{T}\to\mathbb{R}}_+$, and a fatigue degradation function  ${\gamma\mapsto d(\gamma)\in[0,1]}$. This function has the following properties:
\begin{equation}
d(\gamma\leq\gamma_0)=1, \quad d(\gamma>\gamma_0)\in[0,1], \quad d'(\gamma)\leq0,
\end{equation}
where $\gamma_0$ is a material threshold parameter. 

The different dissipative ingredients are introduced in the following dissipation potential:
\begin{equation}
\phi(\dot{\bm{\varepsilon}}^\mathrm{p},\dot{\bm{\kappa}},\dot{\bm{\varepsilon}}^\mathrm{r},\dot{\alpha},\nabla\dot{\alpha};\bm{\kappa},\alpha,\bm{\sigma},\gamma)\coloneqq\underbrace{\phi^\mathrm{p}(\dot{\bm{\varepsilon}}^\mathrm{p},\dot{\bm{\kappa}},\dot{\bm{\varepsilon}}^\mathrm{r};\alpha,\nabla\alpha,\bm{\sigma})}_{\text{plastic dissipation}}+\sum_{s=1}^{n_{\mathrm{y}}}\underbrace{g'(\alpha)\sigma^\mathrm{p}_s\kappa_s\dot{\alpha}}_{\text{coupling}}+\underbrace{\phi^\mathrm{d}(\dot{\alpha},\nabla\dot{\alpha};\alpha,\nabla\alpha,\gamma)}_{\text{damage dissipation}}.
\label{diss_pot_model}
\end{equation}
The coupling term $g'(\alpha)\kappa_s\dot{\alpha}$ was introduced by \citet{AleMarMauVid2017} for perfect plasticity to render the dissipated energy a state function. This idea is also applied here and generalized to the multi-surface hardening case. However, due to the dependence on $\bm{\sigma}$ and $\gamma$, only part of the dissipated energy  becomes a state function. The plastic dissipation corresponds to the damaged version of equation~\eqref{diss_model0}:
\begin{equation}
\begin{aligned}
\phi^\mathrm{p}(\dot{\bm{\varepsilon}}^\mathrm{p},\dot{\bm{\kappa}},\dot{\bm{\varepsilon}}^\mathrm{r};\alpha,\bm{\sigma})=\sum_{s=1}^{n_{\mathrm{y}}} \bigg(\sqrt\frac{2}{3}g(\alpha)\sigma^{\mathrm{p}}_s + \beta\Vert\bm{\sigma}_\mathrm{dev}(\bm{\varepsilon},\bm{\varepsilon}^\mathrm{p},\bm{\varepsilon}^\mathrm{r},\alpha)\Vert\bigg)\Vert\dot{{\bm{\varepsilon}}}^\mathrm{p}_s\Vert.
\end{aligned}
\end{equation}
%\\\equiv\sum_{s=1}^{n_{\mathrm{y}}} \bigg(g(\alpha)\sigma^{\mathrm{p}}_s + \sqrt\frac{3}{2}\beta\Vert\bm{\sigma}_\mathrm{dev}(\bm{\varepsilon},\bm{\varepsilon}^\mathrm{p},\bm{\varepsilon}^\mathrm{r},\alpha)\Vert\bigg)\dot{\kappa}_s
In this expression, we have enforced the constraints that are present in~\eqref{diss_model0} by assuming the hardening and ratcheting laws~\eqref{hardrate} and \eqref{ratchrate}, as will be done hereinafter to simplify the presentation. The degradation function $g(\alpha)$ is used to let the plastic yield strength decrease as a function of damage. More general choices of functions that comply with \eqref{degrad_prop} are possible. Herein, we only consider the quadratic function~\eqref{quad_degr}. The phase-field fracture dissipation potential with fatigue effects reads
\begin{equation}
\phi^\mathrm{d}(\dot{\alpha},\nabla\dot{\alpha};\alpha,\nabla\alpha,\gamma)\coloneqq d(\gamma)\big(w'(\alpha)\dot{\alpha}+\eta_\mathrm{d}^2\nabla\alpha\cdot\nabla\dot{\alpha}\big),
\end{equation}
where $\eta_\mathrm{d}$ is the damage internal length scale. The terms $w'(\alpha)\dot{\alpha}+\eta_\mathrm{d}^2\nabla\alpha\cdot\nabla\dot{\alpha}$ constitute the standard  phase-field dissipation power, whose time integral yields the regularized fracture energy density. This quantity is subject to degradation through $d(\gamma)$, which introduces a path-dependent fatigue effect. The function $w(\alpha)$ represents the local dissipated energy due to damage, for which two models (labeled AT-1 and AT-2 after~\citet{ambrosio1990}) are generally adopted~\citep{marigo2016}:
\begin{equation}
w(\alpha)\coloneqq\begin{dcases}w_0\alpha \quad &\text{AT-1}, \\ 
w_0\alpha^2 \quad &\text{AT-2}.\end{dcases}
\label{w0}
\end{equation} 
With these definitions, the damage characteristic length $\ell_\mathrm{d}$ can be recovered from the relation
\begin{equation}
\eta_\mathrm{d}=\ell_\mathrm{d}\sqrt{2w_0}.
\end{equation}
While the AT-2 model is used in most studies, the AT-1 model has the advantage of including an initial elastic stage before damage is triggered, where $w_0$ is a threshold parameter \citep{marigo2016,UllRodSam2016}. Regarding the fatigue degradation function, different options have been proposed~\citep{Alessi2017fatigue,carrara2020}, from which we adopt the following:
\begin{equation}
d(\gamma)\coloneqq \begin{dcases} 1 & \text{if}  \quad \gamma(\bm{x},t) \leq \gamma_0, \\ 
\bigg[1-k\log\bigg(\frac{\gamma(\bm{x},t)}{\gamma_0}\bigg)\bigg]^2 & \text{if} \quad \gamma_0\leq\gamma(\bm{x},t) \leq \gamma_010^{1/k}, \\ 
0 & \text{otherwise},
\end{dcases}
\end{equation}
where  $k$ is a material parameter that controls the rate of (logarithmic) decay of the fatigue degradation function. The fatigue variable $\gamma$ is defined as 
\begin{equation}
\gamma(\bm{x},t) \coloneqq \int_{0}^t\dot{\vartheta}(\bm{x},s)H\big(\dot{\vartheta}(\bm{x},s)\big)\, \mathrm{d}s, \quad \text{with} \quad \vartheta(\bm{x},t)\coloneqq g(\alpha)\big(\psi^{\mathrm{e}_+}(\bm{\varepsilon},\bm{\varepsilon}^\mathrm{p},\bm{\varepsilon}^\mathrm{r})+\psi^{\mathrm{p}}(\bm{\varepsilon}^\mathrm{p},\bm{\kappa},\nabla\bm{\kappa})\big).
\label{fatfield}
\end{equation}
While other definitions are possible for $\vartheta$, e.g., an accumulated strain measure, using the strain energy density ensures mesh-objectivity \citep{carrara2020}. Moreover, the Heaviside function $H$ precludes fatigue degradation in unloading stages. Note that in the present model, the plastic free energy also contributes to the evolution~of~$\gamma$.  

Table~\ref{overview2} presents a summary of the different mechanisms included in the present model. In the following subsections, these definitions will be used to derive the evolution problem according to the formulation elaborated in section~\ref{variat_general} (see table~\ref{overview1}).

\begin{table}[hb!]
\small
\centering
\caption{State variables and parameters included in the multifield phase-field-based fatigue model with cyclic plasticity. These variables enter the free energy density~\eqref{free_ener_model2} and the dissipation potential~\eqref{diss_pot_model}, which in turn constitute the basic ingredients of the energetic formulation summarized in table~\ref{overview1}.}
\begin{tabular}{llll}
\toprule
  & \multicolumn{2}{c}{\sffamily \textbf{Variable} \hspace{10mm}~}   
  & \sffamily \textbf{Type}
\\ \cmidrule{2-4} 
\multirow{2}{*}{\normalsize\color{MidnightBlue} \textsc{Elasticity}}   
  & Displacements 
  & $\bm{u}$           
  & Independent      
\\ \mycmidrule
  & Strains 
  & $\bm{\varepsilon}$          
  & Dependent        
\\ 
\midrule                      
\multirow{4}{*}{\normalsize \color{OliveGreen}\textsc{Plasticity}}   
  & Plastic strains 
  & $\bm{\varepsilon}^\mathrm{p}=\{\bm{\varepsilon}^\mathrm{p}_1,\dots,\bm{\varepsilon}^\mathrm{p}_s,\dots,\bm{\varepsilon}^\mathrm{p}_{n_\mathrm{y}}\}$              
  & Independent      \\ \mycmidrule
  & Hardening/softening variables 
  & $\bm{\kappa}=\{\kappa_1,\dots,\kappa_s,\dots\kappa_{n_\mathrm{y}}\}$              
  & Dependent, irreversible        \\ \mycmidrule
  & Ratcheting strains  
  & $\bm{\varepsilon}^\mathrm{r}$        
  & Dependent        \\ \mycmidrule
  & Plastic gradients
  & $\nabla\bm{\kappa}=\{\nabla\kappa_1,\dots,\nabla\kappa_s,\dots\nabla\kappa_{n_\mathrm{y}}\}$          
                              & Dependent      
\\ 
\midrule
\multirow{2}{*}{\normalsize\color{BrickRed} \textsc{Damage}}
  & Phase-field/damage variable
  & $\alpha$        
  & Independent, irreversible    \\ \mycmidrule
  & Damage gradient  
  & $\nabla\alpha$        
  & Dependent      
\\ 
\midrule
\textsc{\normalsize\textcolor{black!80}{Fatigue}}
  & Fatigue parameter  
  & $\gamma$
  & History-dependent, irreversible      
\\ 
\bottomrule
\end{tabular}
\label{overview2}
\end{table}
%\multicolumn{4}{c}{Basic energy  quantities}       \\ \hline
%Free energy                   &     $\psi(\bm{\varepsilon},\bm{\varepsilon}^\mathrm{p},\bm{\kappa},\bm{\varepsilon}^\mathrm{r},\alpha,\nabla\bm{\kappa})$    &      Equation~\eqref{free_ener_model2} & \\ 
%Dissipation                   &     $\phi(\dot{\bm{\varepsilon}}^\mathrm{p},\dot{\bm{\kappa}},\dot{\bm{\varepsilon}}^\mathrm{r},\dot{\alpha},\nabla\dot{\alpha};\bm{\kappa},\alpha,\bm{\sigma},\gamma)$    &      Equation~\eqref{diss_pot_model} & \\ \hline

\subsection{Energetic formulation} \label{energ_model}

With the previous definitions, the governing equations of the proposed model can be derived from the principles of the energetic formulation. For this purpose, we define the space of displacement fields in equation~\eqref{dispspace} by $\mathscr{F}\coloneqq \mathrm{H}^1(\Omega;\mathbb{R}^d)$, and the following function spaces for the internal variables: 
\begin{align} 
&\mathscr{B}\coloneqq \mathrm{L}^2(\Omega;\mathbb{R}^{d\times d}_{\mathrm{dev}}), &&\tilde{\mathscr{B}}\equiv\mathscr{B},  \\
&\mathscr{K}\coloneqq \mathrm{H}^1(\Omega;\mathbb{R}_+),  &&\tilde{\mathscr{K}}(\bm{b})\coloneqq \{z\in\mathscr{K} \  | \  z = \sqrt{2/3}\,\Vert\bm{b}\Vert,\,\bm{b}\in\tilde{\mathscr{B}}\},  \\
&\mathscr{R}\coloneqq \mathrm{L}^2(\Omega;\mathbb{R}^{d\times d}_{\mathrm{dev}}),  &&\tilde{\mathscr{R}}(\bm{b})\coloneqq \{\bm{z}\in\mathscr{R} \ | \ \Vert\bm{z}\Vert =  \textstyle{\beta\sum_{s=1}^{n_\mathrm{y}}}\Vert\bm{b}_s\Vert, \ {\bm{b}_s}\in\tilde{\mathscr{B}},\ \bm{b}_s\in\bm{b}\},\\
&\mathscr{D}\coloneqq \mathrm{H}^1(\Omega;[0,1]), &&\tilde{\mathscr{D}}\coloneqq \mathrm{H}^1(\Omega;\mathbb{R}_+).
\label{funscspaces}
\end{align}%\bm{m}\in\mathbb{R}^{d\times d}_{\mathrm{dev}}
The admissible primary fields are then given by
\begin{equation}
\mathbf{q}=\{\bm{u},\bm{\varepsilon}^\mathrm{p},\bm{\kappa},\bm{\varepsilon}^\mathrm{r},\alpha\}\in\mathscr{Q}\coloneqq \mathscr{U}\times\mathscr{P}\times\mathscr{D}, \quad \text{with} \quad \mathscr{P}\coloneqq \mathscr{B}\times\dots\times\mathscr{B}\times\mathscr{K}\times\dots\times\mathscr{K}\times\mathscr{R}.
\label{funcspace}
\end{equation} %\mathscr{B}\times\dots\times\mathscr{B}\times\mathscr{K}\times\dots\times\mathscr{K}\times\mathscr{R}\times\mathscr{D}
The space of admissible variations must ensure that $\phi(\tilde{\bm{\varepsilon}}^\mathrm{p},\tilde{\bm{\kappa}},\tilde{\bm{\varepsilon}}^\mathrm{r},\tilde{\alpha},\nabla\tilde{\alpha};\bm{\kappa},\alpha,\bm{\sigma},\gamma)<+\infty$. In this way, the global dissipative power entering the stability condition~\eqref{DS-1} remains finite, and we are able to consider the non-trivial conditions for stability. For the sake of conciseness, this requirement has been already embedded in the definition of the function spaces through the constraints. Thus, we set
%must account for evolution constraints of the internal fields, i.e., the plastic flow relations~\eqref{hardrate} and~\eqref{ratchrate}, and the irreversibility condition~\eqref{IR_model}. Therefore, we set
%
\begin{equation}
\tilde{\mathbf{q}}\in\tilde{\mathscr{Q}}\coloneqq \tilde{\mathscr{U}}\times\tilde{\mathscr{B}}\times\dots\times\tilde{\mathscr{B}}\times\tilde{\mathscr{K}}(\tilde{\bm{\varepsilon}}^\mathrm{p}_1)\times\dots\times\tilde{\mathscr{K}}(\tilde{\bm{\varepsilon}}^\mathrm{p}_{n_\mathrm{y}})\times\tilde{\mathscr{R}}\big(\tilde{\bm{\varepsilon}}^\mathrm{p}\big)\times\tilde{\mathscr{D}}.
\label{virtspace}
\end{equation}
The governing equations directly follow from enforcing the first-order stability condition~\eqref{DS-1}, the power balance~\eqref{EB-1} and the dissipation inequality~\eqref{DI}. Exploiting the generality of the formulation presented in section~\ref{energetic_form}, we directly replace the free energy density~\eqref{free_ener_model2} and the dissipation potential~\eqref{diss_pot_model} in the weak form of the first-order stability condition~\eqref{energ_stab} to obtain, for all $\tilde{\mathbf{q}}\in\tilde{\mathscr{Q}}$:
\begin{equation}
\begin{aligned}
\int_{\Omega}&\bigg[\sum_{s=1}^{n_\mathrm{y}}\bigg[\big(\bm{\sigma}-g(\alpha)H^\mathrm{kin}_s\bm{\varepsilon}^\mathrm{p}_s\big):\hat{\bm{n}}_s - \textstyle\sqrt\frac{2}{3}\bigg(g(\alpha)\big(\sigma^{\mathrm{p}}_s+H^\mathrm{iso}_s\kappa_s\big)- \eta_{\mathrm{p}s}^2\mathrm{div}\big[g(\alpha)\nabla\kappa_s\big]\bigg)\bigg]\Vert\tilde{\bm{\varepsilon}}^\mathrm{p}_s\Vert- \bigg(g'(\alpha)\big(\psi^{\mathrm{e_+}}+\psi^{\mathrm{p}}\big) \\
 & + g'(\alpha)\sum_{s=1}^{n_\mathrm{y}}\sigma^{\mathrm{p}}_s\kappa_s+d(\gamma)w'(\alpha) - \eta_\mathrm{d}^2\mathrm{div}\big[d(\gamma)\nabla\alpha\big]\bigg) \tilde{\alpha}\bigg]\, \mathrm{d}\bm{x} -\int_{\Gamma}\bigg(d(\gamma)\eta_\mathrm{d}^2\nabla\alpha\tilde{\alpha}+g(\alpha)\sum_{s=1}^{n_\mathrm{y}}\eta_{\mathrm{p}s}^2\nabla\kappa_s\tilde{\kappa}_s\bigg)\cdot\bm{n}\, \mathrm{d}S 
\\ & \ \hspace*{7.5cm} +\int_{\Omega}(\mathrm{div}\bm{\sigma}+\rho\bm{b})\cdot\tilde{\bm{u}}\, \mathrm{d}\bm{x}-\int_{\Gamma_\mathrm{N}}(\bm{\sigma}\cdot\bm{n}-\bar{\bm{t}}) \cdot\tilde{\bm{u}}\, \mathrm{d}S \leq0.
\end{aligned}
\label{model_stab}
\end{equation}
Likewise, using the plastic flow relations~\eqref{hardrate} and \eqref{ratchrate},  the power balance principle~\eqref{energ_energ} gives
\begin{equation}
\begin{aligned}
\int_{\Omega}&\bigg[\sum_{s=1}^{n_\mathrm{y}}\bigg[(\bm{\sigma}-g(\alpha)H^\mathrm{kin}_s\bm{\varepsilon}^\mathrm{p}_s):\hat{\bm{n}}_s - \textstyle\sqrt\frac{2}{3}\bigg(g(\alpha)\big(\sigma^{\mathrm{p}}_s+H^\mathrm{iso}_s\kappa_s\big)- \eta_{\mathrm{p}s}^2\mathrm{div}\big[g(\alpha)\nabla\kappa_s\big]\bigg)\bigg]\Vert\dot{\bm{\varepsilon}}^\mathrm{p}_s\Vert - \bigg(g'(\alpha)\big(\psi^{\mathrm{e_+}}+\psi^{\mathrm{p}}\big) \\
 & + g'(\alpha)\sum_{s=1}^{n_\mathrm{y}}\sigma^{\mathrm{p}}_s\kappa_s+ d(\gamma)w'(\alpha) - \eta_\mathrm{d}^2\mathrm{div}\big[d(\gamma)\nabla\alpha\big]\bigg) \dot{\alpha}\bigg]\, \mathrm{d}\bm{x} -\int_{\Gamma}\bigg(d(\gamma)\eta_\mathrm{d}^2\nabla\alpha\dot{\alpha}+g(\alpha)\sum_{s=1}^{n_\mathrm{y}}\eta_{\mathrm{p}s}^2\nabla\kappa_s\dot{\kappa}_s\bigg)\cdot\bm{n}\, \mathrm{d}S
\\ & \ \hspace*{3.7cm} +\int_{\Omega}(\mathrm{div}\bm{\sigma}+\rho\bm{b})\cdot\dot{\bm{u}}\, \mathrm{d}\bm{x} -\int_{\Gamma_{\mathrm{N}}}(\bm{\sigma}\cdot\bm{n}-\bm{\bar{t}}\,)\cdot\dot{\bm{u}}\, \mathrm{d}S-\int_{\Gamma_{\mathrm{D}}}(\bm{\sigma}\cdot\bm{n}-\bm{t}_\mathrm{r})\cdot\dot{\bar{\bm{u}}}\, \mathrm{d}S =0.
\end{aligned}
\label{model_energ}
\end{equation}
As in section~\ref{gov_Eq_ener}, equations~\eqref{model_stab} and \eqref{model_energ} yield the equilibrium equations and boundary conditions~\eqref{st_adm}, along with equilibrium at the Dirichlet boundary.
%
%\begin{equation}
%\begin{dcases}
%\nabla\cdot\bm{\sigma}+\rho\bm{b} = \bm{0} & \text{in} \quad \Omega,\\
%\bm{\sigma}\cdot\bm{n}=\bar{\bm{t}} & \text{on} \quad \Gamma_\mathrm{N},\\
%\bm{\sigma}\cdot\bm{n}=\bm{t}_\mathrm{r} & \text{on} \quad \Gamma_\mathrm{D},
%\end{dcases}
%\label{model_eq}
%\end{equation}
%
On the other hand, defining the $s^\text{th}$ plastic yield function~as
\begin{equation}
f^\mathrm{p}_s\coloneqq \big\Vert\bm{\sigma}_\mathrm{dev}-g(\alpha)H^\mathrm{kin}_s\bm{\varepsilon}^\mathrm{p}_s\big\Vert - \sqrt\frac{2}{3}\bigg(g(\alpha)\big(\sigma^{\mathrm{p}}_s+H^\mathrm{iso}_s\kappa_s\big)- \eta_{\mathrm{p}s}^2\mathrm{div}\big[g(\alpha)\nabla\kappa_s\big]\bigg),
\label{fp}
\end{equation}
and the damage yield function as
\begin{equation}
f^\mathrm{d}\coloneqq -g'(\alpha)\big(\psi^{\mathrm{e_+}}+\psi^{\mathrm{p}}\big) - g'(\alpha)\sum_{s=1}^{n_\mathrm{y}}\sigma^{\mathrm{p}}_s\kappa_s-d(\gamma)w'(\alpha) + \eta_\mathrm{d}^2\mathrm{div} \big[d(\gamma)\nabla\alpha\big],
\label{fd}
\end{equation} 
the evolution equations emerge from equations~\eqref{model_stab}, \eqref{model_energ} and the irreversibility condition \eqref{IR_model} as loading/unloading systems, along with boundary conditions of the gradient-enhanced variables: 
\begin{equation}
\begin{dcases}
f^\mathrm{p}_s\leq0, \quad f^\mathrm{p}_s\dot{\kappa}_s=0, \quad \dot{\kappa}_s\geq0,  \\
\nabla\kappa_s\cdot\bm{n}=0, \quad \forall \, s\in\mathrm{Y},
\end{dcases}
\quad \begin{dcases}
f^\mathrm{d}\leq0, \quad f^\mathrm{d}\dot{\alpha}=0, \quad \dot{\alpha}\geq0, \\
\nabla\alpha\cdot\bm{n}=0.
\end{dcases}
\label{load_unload_model}
\end{equation}

\subsection{An overview of the proposed model}

Equation~\eqref{fp} can be viewed as a damage-dependent version of the plastic yield function~\eqref{yield_s}. Specifically, the plastic driving force $\Vert\bm{\sigma}_\mathrm{dev}-g(\alpha)H^\mathrm{kin}_s\bm{\varepsilon}^\mathrm{p}_s\Vert$ is now a function of damage through the stress tensor~\eqref{stress_model} and the back stress $g(\alpha)H^\mathrm{kin}_s\bm{\varepsilon}^\mathrm{p}_s$. Moreover, the size of the yield surface is now given by $$\sqrt\frac{2}{3}\bigg(g(\alpha)\big(\sigma^{\mathrm{p}}_s+H^\mathrm{iso}_s\kappa_s\big) - \eta_{\mathrm{p}s}^2\mathrm{div}\big[g(\alpha)\nabla\kappa_s\big]\bigg),$$ which progressively decreases as $\alpha\to 1$. 
%, interpreted as the plastic resisting force,

In the damage yield function~\eqref{fd}, we identify 
\begin{equation}
R(\alpha,\nabla\alpha,\gamma)\coloneqq d(\gamma)w'(\alpha) - \eta_\mathrm{d}^2\mathrm{div} \big[d(\gamma)\nabla\alpha\big]
\label{damageresisting}
\end{equation} 
as the damage resisting force with fatigue effects. We denote the damage driving force by ${D\coloneqq D^{\mathrm{e}}+D^{\mathrm{p}}}$, where $D^{\mathrm{e}}$ and $D^{\mathrm{p}}$ are, respectively, elastic and plastic contributions given by
\begin{equation}
D^{\mathrm{e}}(\bm{\varepsilon}^\mathrm{e},\alpha)\coloneqq -g'(\alpha)\psi^{\mathrm{e_+}}\quad \text{and} \quad D^{\mathrm{p}}(\bm{\varepsilon}^\mathrm{p},\bm{\kappa},\nabla\bm{\kappa},\alpha)\coloneqq -g'(\alpha)\bigg(\psi^{\mathrm{p}} + \sum_{s=1}^{n_\mathrm{y}}\sigma^{\mathrm{p}}_s\kappa_s\bigg).
\label{damagedriving}
\end{equation}
As in standard phase-field models, the notion of \emph{driving} and \emph{resisting} forces refers to the energetic competition between the damage energy release rate and the critical fracture energy~\citep{MieHofWel2010}. The main difference is that in the present model, the driving force contains both elastic and plastic contributions, while the resisting force is degraded through the fatigue variable $\gamma$. These definitions are useful to describe the mechanical response in the numerical simulations presented in section~\ref{numsym}. Equation~\eqref{fd} can be written as 
\begin{equation}
f^{\mathrm{d}}=D(\bm{\varepsilon}^\mathrm{e},\bm{\varepsilon}^\mathrm{p},\bm{\kappa},\nabla\bm{\kappa},\alpha)-R(\alpha,\nabla\alpha,\gamma)=D^{\mathrm{e}}(\bm{\varepsilon}^\mathrm{e},\alpha)+D^{\mathrm{p}}(\bm{\varepsilon}^\mathrm{p},\bm{\kappa},\nabla\bm{\kappa},\alpha)-R(\alpha,\nabla\alpha,\gamma).
\label{fd_1}
\end{equation}
%%Here, we associate the plastic contribution to an increase in the damage driving force. A similar interpretation is made in the model of~\citet{MieTeiAld2016}, where, however, the plastic strain evolution is viewed as part of the resisting force. 
%
We observe in equation~\eqref{fd_1} two distinctive fatigue mechanisms: 

\begin{enumerate}
\item The accumulation of elastic free energy, which, in the absence of plastic strains, drives cracks in the high-cycle fatigue regime. This effect is attained by the multiplicative degradation of the damage resisting force through $d(\gamma)$. 

\item The accumulation of plastic energy (free and dissipated), which drives cracks in the low-cycle fatigue regime by additively increasing the damage driving force through~$D^{\mathrm{p}}$. This mechanism entails, on its own, a low-cycle fatigue process, that is accelerated when $d(\gamma)< 1$, where the plastic strains also contribute to the degradation of the damage resisting force (equation~\eqref{fatfield}). We associate the combination of $D^{\mathrm{p}}>0$ and $d(\gamma)< 1$ with very-low-cycle fatigue.
\end{enumerate}

%\begin{enumerate}
%\item The multiplicative degradation of the damage resisting force caused by the accumulation of free energy through $d(\gamma)$. In the absence of plastic strains, this mechanism describes the characteristic behavior of high-cycle fatigue by means of elastic free energy accumulation.
%\item The additive increase in the damage driving force caused by the accumulation of plastic strains through~$D^{\mathrm{p}}$. This mechanism entails, on its own, a low-cycle fatigue process, that can be accelerated when $d(\gamma)< 1$, which entails further degradation of the damage resisting force through a combination of elastic and plastic free energy. This type of response characterizes a very-low-cycle fatigue process.
%\end{enumerate}

%\noindent The interaction between these mechanisms will become clear in the numerical simulations presented in Section~\ref{numsym}.

\subsection{Incremental minimization and numerical implementation}

%The first step is the construction of the global potential energy functional~\eqref{stored_ext0} and the global dissipative power~\eqref{dissipated}. These quantities directly follow from replacing the free the free energy density~\eqref{free_ener_model2} in Equation~\eqref{stored_ext0} and the dissipation potential~\eqref{diss_pot_model} in Equation~\eqref{dissipated}. Therefore, we can immediately use the results of Section~\ref{gov_Eq_ener} to obtain the governing equations that fulfill~\eqref{DS-1}, \eqref{EB-1} and~\eqref{DI}. The main steps are briefly shown below for the sake of completeness.

As in section~\ref{incremental}, consider the time discretization ${0=t_0<\dots<t_{n}<t_{n+1}<\dots<t_{n_\mathrm{t}}=t_\mathrm{max}}$, where all quantities are known up to $t_n$, and the goal is to find the state at the current time step $t_{n+1}$. For convenience, we introduce the following  notations. A quantity $\Box$ evaluated at any previous time step $0\leq i\leq n$ is denoted as $\Box_i$, while a quantity evaluated at $t_{n+1}$ is written without a subscript, i.e., $\Box\coloneqq \Box_{n+1}$. Moreover, the operator $\Delta\Box\coloneqq \Box-\Box_n$ is used to denote an increment from $t_{n}$ to $t_{n+1}$.

Recall that the dissipation potential~\eqref{diss_pot_model} explicitly depends on the stress $\bm{\sigma}$ and the history variable $\gamma$. As a consequence, the dissipated energy is not a state function. In agreement with section~\ref{incremental}, the dissipated energy increment is evaluated as
\begin{equation}
\begin{aligned}
&\int_\Omega\int_{t_n}^{t_{n+1}}\bigg[\sum_{s=1}^{n_{\mathrm{y}}}\bigg(g(\alpha)\sigma^{\mathrm{p}}_s\dot{\kappa}_s+g'(\alpha)\sigma^\mathrm{p}_s\kappa_s\dot{\alpha}+\sqrt{\frac{3}{2}}\beta\Vert\bm{\sigma}_\mathrm{dev}\Vert\dot{\kappa}_s\bigg)+d(\gamma)\big(w'(\alpha)\dot{\alpha}+\eta_\mathrm{d}^2\nabla\alpha\cdot\nabla\dot{\alpha}\big)\bigg]\, \mathrm{d}t\,\mathrm{d}\bm{x}\\
=&\int_{\Omega}\int_{t_n}^{t_{n+1}}\bigg[\sum_{s=1}^{n_{\mathrm{y}}}\underbrace{\bigg(g(\alpha)\sigma^{\mathrm{p}}_s\dot{\kappa}_s+g'(\alpha)\sigma^\mathrm{p}_s\kappa_s\dot{\alpha}\bigg)}_{\text{path-independent}}+\underbrace{\vphantom{\bigg(}\sqrt{\frac{3}{2}}\beta\Vert\bm{\sigma}_\mathrm{dev}\Vert\sum_{s=1}^{n_{\mathrm{y}}}\dot{\kappa}_s+d(\gamma)\frac{\mathrm{d}}{\mathrm{d} t}\bigg(w(\alpha)+\frac{1}{2}\eta_\mathrm{d}^2\nabla\alpha\cdot\nabla\alpha\bigg)}_{\text{path-dependent}}\bigg]\, \mathrm{d}t\,\mathrm{d}\bm{x}.
\end{aligned}
\end{equation}
%\\
%&\approx\int_{\Omega}\int_{t_n}^{t_{n+1}}\bigg[\sum_{s=1}^{n_{\mathrm{y}}}\bigg(g(\alpha)\sigma^{\mathrm{p}}_s\dot{\kappa}_s+g'(\alpha)\sigma^\mathrm{p}_s\kappa_s\dot{\alpha}\bigg)+\vphantom{\bigg(}\sqrt{\frac{3}{2}}\beta\Vert\bm{\sigma}_{\mathrm{dev}\,n}\Vert\sum_{s=1}^{n_{\mathrm{y}}}\dot{\kappa}_s+d(\gamma_n)\frac{\mathrm{d}}{\mathrm{d} t}\bigg(w(\alpha)+\frac{1}{2}\eta_\mathrm{d}^2\nabla\alpha\cdot\nabla\alpha\bigg)\bigg]\, \mathrm{d}t\,\mathrm{d}\bm{x},
Evaluating the path-independent part directly and using the incremental approximation~\eqref{disspotapprox} for the path-dependent part, we obtain
\begin{equation}
\begin{aligned}
\mathcal{D}(\bm{\kappa}_n,\alpha_n,\bm{\kappa},\alpha;\bm{\sigma}_n,\gamma_n)
&=\int_{\Omega}\bigg[\sum_{s=1}^{n_{\mathrm{y}}}g(\alpha)\sigma^{\mathrm{p}}_s\kappa_{s}+\sqrt{\frac{3}{2}}\beta\Vert\bm{\sigma}_{\mathrm{dev}\,n}\Vert\sum_{s=1}^{n_{\mathrm{y}}}{\kappa}_s+d(\gamma_n)\big(w(\alpha)+\frac{1}{2}\eta_\mathrm{d}^2\nabla\alpha\cdot\nabla\alpha\big)\bigg]\bigg|_{n}^{n+1}\,\mathrm{d}\bm{x}.
\end{aligned}
\label{inc_diss_model0}
\end{equation}
In view of equation~\eqref{discminprim}, a time-discrete energy functional is defined as
\begin{equation}
\begin{aligned}
\Pi(\bm{u},\bm{\varepsilon}^\mathrm{p},\bm{\kappa},\bm{\varepsilon}^\mathrm{r},\alpha)
=\mathcal{E}(\bm{u},\bm{\varepsilon}^\mathrm{p},\bm{\kappa},\bm{\varepsilon}^\mathrm{r},\alpha)-\mathcal{T}(\bm{u}_n,\bm{u}) + \mathcal{D}(\bm{\kappa}_n,\alpha_n,\bm{\kappa},\alpha;\bm{\sigma}_n,\gamma_n),
\end{aligned}
\label{discfunctprim_model}
\end{equation}% - \int_{\Omega}\bigg(\sum_{s=1}^{n_\mathrm{y}}I_+(\Delta \kappa_s)+I_+(\Delta \alpha)\bigg)\, \mathrm{d}\bm{x}& where the indicator function $I_+\colon\mathbb{R}\to\{0,\infty\}$ is introduced to impose irreversibility. 
such that the incremental minimization problem~\eqref{discminprim} takes the form
\begin{equation}
\inf_{\{\bm{u},\bm{\varepsilon}^\mathrm{p},\bm{\kappa},\bm{\varepsilon}^\mathrm{r},\alpha\}\in \mathscr{Q}} \big\{\Pi(\bm{u},\bm{\varepsilon}^\mathrm{p},\bm{\kappa},\bm{\varepsilon}^\mathrm{r},\alpha)\,\, | \, \,\Delta \kappa_s=\textstyle\sqrt\frac{2}{3}\Vert\Delta\bm{\varepsilon}^\mathrm{p}_s\Vert \, \forall s\in\mathrm{Y},\,\, \Delta\bm{\varepsilon}^\mathrm{r}=\textstyle \beta\hat{\bm{r}}_n\sum_{s=1}^{n_\mathrm{y}}\Vert\Delta\bm{\varepsilon}^\mathrm{p}_s\Vert,\,  \,  \Delta\alpha\geq 0 \big\},
\label{discminprim_model}
\end{equation}
%
%\substack{\mathrm{inf}_{\bm{u}},\mathrm{inf}_{\mathbf{a}\ \\ \Vert\Delta\bm{\varepsilon}^\mathrm{p}_s\Vert=\Vert\Delta\bm{\varepsilon}^\mathrm{p}_s\Vert \, \forall s\in\mathrm{Y}\\ \Delta\bm{\varepsilon}^\mathrm{r}=\hat{\bm{r}}_n\sum_{s=1}^{n_\mathrm{y}}R_s\Vert\Delta\bm{\varepsilon}^\mathrm{p}_s\Vert}}
where we impose irreversibility and the incremental counterparts of the plastic flow relations~\eqref{hardrate} and~\eqref{ratchrate}. %The first variation of $\Pi$ reads, for all $\tilde{mathbf{q}}\in\tilde{\mathscr{B}}$:

The numerical solution of the variational problem~\eqref{discminprim_model} is based on an extension of the alternate minimization algorithm~\citep{bourdin2007} to ductile fracture~\citep{AleMarVid2015,UllRodSam2016}, where the solution of a sequence of convex optimization problems is sought, aiming to iteratively converge to a stationary state. The procedure is summarized in algorithm~\ref{alg1}, along with the equations presented below. 
%Algorithms of this kind have shown robustness in brittle fracture simulations with respect to monolithic schemes~\citep{gerasimov2016}. , which exploits the fact that while the functional~\eqref{discfunctprim_model} is non-convex, it is convex with respect to each primary field individually

\begin{algorithm}[!ht]
\small
 \hspace*{\algorithmicindent} \textbf{Input}: $\mathbf{q}_i\in\mathscr{Q}$ for all $i=0,\dots,n$. \\
 \hspace*{\algorithmicindent} \textbf{Output}: $\mathbf{q}_{n+1}\in\mathscr{Q}$. 
\begin{algorithmic}[1]
%\WHILE { \ are \\
\State Initialize iterations $j=0$ and set $\{\bm{\varepsilon}^{\mathrm{p}(0)},\bm{\kappa}^{(0)},\bm{\varepsilon}^{\mathrm{r}(0)},\alpha^{(0)}\}\coloneqq \{\bm{\varepsilon}^{\mathrm{p}}_n,\bm{\kappa}_n,\bm{\varepsilon}^{\mathrm{r}}_n,\alpha_n\}$.
\Repeat  
\State Set $j \leftarrow j+1$.
\State Find 
\begin{equation*}
\bm{u}^{(j)}\coloneqq \arginf_{\bm{u}\in\mathscr{U}} \Pi\big(\bm{u},{\bm{\varepsilon}^\mathrm{p}}^{(j-1)},\bm{\kappa}^{(j-1)},{\bm{\varepsilon}^\mathrm{r}}^{(j-1)},\alpha^{(j-1)}\big)
\end{equation*}
\hspace*{0.41cm} from equation~\eqref{weak_eq_model}.
 \State Find \begin{equation*}
\begin{aligned}
 \{{\bm{\varepsilon}^\mathrm{p}}^{(j)},\bm{\kappa}^{(j)},{\bm{\varepsilon}^\mathrm{r}}^{(j)}\}\coloneqq&\arginf_{\{\bm{\varepsilon}^\mathrm{p},\bm{\kappa},\bm{\varepsilon}^\mathrm{r}\}\in\mathscr{P}} \{\Pi(\bm{u}^{(j)},{\bm{\varepsilon}^\mathrm{p}},\bm{\kappa},{\bm{\varepsilon}^\mathrm{r}},\alpha^{(j-1)})\} \\ &\text{s.t.} \quad \Delta \kappa_s=\textstyle\sqrt\frac{2}{3}\Vert\Delta\bm{\varepsilon}^\mathrm{p}_s\Vert \, \forall s\in\mathrm{Y},\,\, \Delta\bm{\varepsilon}^\mathrm{r}=\textstyle \beta\hat{\bm{r}}_n\sum_{s=1}^{n_\mathrm{y}}\Vert\Delta\bm{\varepsilon}^\mathrm{p}_s\Vert
\end{aligned}
\end{equation*}
\hspace*{0.41cm} from equation~\eqref{weak_p_model}.
\State Find
\begin{equation*}
\alpha^{(j)}=\arginf_{\alpha\in \mathscr{D}} \Pi\big(\bm{u}^{(j)},{\bm{\varepsilon}^\mathrm{p}}^{(j)},\bm{\kappa}^{(j)},{\bm{\varepsilon}^\mathrm{r}}^{(j)},\alpha\big)
\end{equation*}
\hspace*{0.41cm} from equation~\eqref{weak_d_model}.
\Until $\Vert\bm{u}^{j}-\bm{u}^{j-1}\Vert$, $\Vert\bm{\kappa}^{j}-\bm{\kappa}^{j-1}\Vert$ and $\Vert {\alpha}^{j}-\alpha^{j-1}\Vert$ are sufficiently small.
\State Set $\mathbf{q}_{n+1}\coloneqq \mathbf{q}^{(j)}$.
%\ENDWHILE
\end{algorithmic}
\caption{Alternate minimization.}\label{alg1}
\end{algorithm} 

Given $\mathbf{q}_i$ for all $i\in\{0,\dots,n\}$, we solve, alternatively, the following coupled sub-problems. %\{{\bm{\varepsilon}^\mathrm{p}}^{(j)},\bm{\kappa}^{(j)},{\bm{\varepsilon}^\mathrm{r}}^{(j)},\alpha^{(j)}\}$

\paragraph*{Minimization with respect to the displacement field} Given $\{\bm{\varepsilon}^\mathrm{p},\bm{\kappa},\bm{\varepsilon}^\mathrm{r},\alpha\}$, 
find 
\begin{equation}
\bm{u}=\arginf_{\bm{u}\in\mathscr{U}} \{\Pi(\bm{u},\bm{\varepsilon}^\mathrm{p},\bm{\kappa},\bm{\varepsilon}^\mathrm{r},\alpha)\}
\end{equation}
from the necessary condition 
\begin{equation}
\int_{\Omega}\big(\bm{\sigma}:\nabla^{\mathrm{s}}\tilde{\bm{u}}-\rho\bm{b}\cdot\tilde{\bm{u}}\big)\,\mathrm{d}\bm{x}-\int_{\Gamma_\mathrm{N}}\bar{\bm{t}}\cdot\tilde{\bm{u}}\,\mathrm{d}S=0 \quad \forall \, \tilde{\bm{u}}\in\tilde{\mathscr{U}}.
\label{weak_eq_model}
\end{equation} %\delta\Pi(\bm{u},\bm{\varepsilon}^\mathrm{p},\bm{\kappa},\bm{\varepsilon}^\mathrm{r},\alpha)(\tilde{\bm{u}},\mathbf{0},\bm{0},\bm{0},0)=
Given that $\bm{u}\in\mathscr{U}$, this problem constitutes the weak form of the equilibrium equations~\eqref{st_adm}. This equation is non-linear due to the elastic energy decomposition~\eqref{split}.

\paragraph*{Minimization with respect to the plastic fields} Given $\{\bm{u},\alpha\}$, 
find 
\begin{equation}
\begin{aligned}
\{\bm{\varepsilon}^\mathrm{p},\bm{\kappa},\bm{\varepsilon}^\mathrm{r}\}=&\arginf_{\{\bm{\varepsilon}^\mathrm{p},\bm{\kappa},\bm{\varepsilon}^\mathrm{r}\}\in\mathscr{P}} \{\Pi(\bm{u},\bm{\varepsilon}^\mathrm{p},\bm{\kappa},\bm{\varepsilon}^\mathrm{r},\alpha)\} \\ &\text{s.t.} \quad \Delta \kappa_s=\textstyle\sqrt\frac{2}{3}\Vert\Delta\bm{\varepsilon}^\mathrm{p}_s\Vert \, \forall s\in\mathrm{Y},\,\, \Delta\bm{\varepsilon}^\mathrm{r}=\textstyle \beta\hat{\bm{r}}_n\sum_{s=1}^{n_\mathrm{y}}\Vert\Delta\bm{\varepsilon}^\mathrm{p}_s\Vert.
\end{aligned}
\label{minplast0}
\end{equation}
For admissible variations $\{\tilde{\bm{\varepsilon}}^\mathrm{p}_s,\tilde{\kappa}_s,\tilde{\bm{\varepsilon}}^\mathrm{r}\}\in\tilde{\mathscr{B}}\times\tilde{\mathscr{K}}(\tilde{\bm{\varepsilon}}^\mathrm{p}_s)\times\tilde{\mathscr{R}}(\tilde{\bm{\varepsilon}}^\mathrm{p})$, it can be shown that the functional derivative of~\eqref{discfunctprim_model} with respect to the plastic variables gives the weak form for the $s^\text{th}$ yield surface
% for the $s^\text{th}$ yield surface
\begin{equation}
\begin{aligned}
\int_{\Omega}\bigg[\bigg(-\textstyle\sqrt\frac{3}{2}\big\Vert\bm{\sigma}_\mathrm{dev}-&g(\alpha)H^\mathrm{kin}_s\bm{\varepsilon}^\mathrm{p}_s\big\Vert  + g(\alpha)(\sigma^{\mathrm{p}}_s+H^\mathrm{iso}_s\kappa_s)-\textstyle\sqrt\frac{3}{2}\beta\Delta\bm{\sigma}:\hat{\bm{r}}_n \bigg)\tilde{\kappa}_s + g(\alpha)\eta_{\mathrm{p}s}^2\nabla\kappa_s\cdot\nabla\tilde{\kappa}_s\bigg]\,\mathrm{d}\bm{x} \geq 0,\\ &\hspace*{0.5cm} \text{with} \quad  \kappa_s=\kappa_{sn}+\textstyle\sqrt\frac{2}{3}\Vert\Delta\bm{\varepsilon}^\mathrm{p}_s\Vert,\,\, \bm{\varepsilon}^\mathrm{r}=\bm{\varepsilon}^\mathrm{r}_n+\textstyle \beta\hat{\bm{r}}_n\sum_{j=1}^{n_\mathrm{y}}\Vert\Delta\bm{\varepsilon}^\mathrm{p}_j\Vert.
\end{aligned}
\label{weak_p_model0}
\end{equation}
%(\bm{\varepsilon},\bm{\varepsilon}^\mathrm{p},\bm{\varepsilon}^\mathrm{r},\alpha)
%\\  &\quad \text{with} \quad  \bm{\varepsilon}^\mathrm{p}_s=\bm{\varepsilon}^\mathrm{p}_{sn}+\textstyle\sqrt\frac{3}{2}\hat{\bm{n}}^\mathrm{tr}_s\Delta\kappa_s \quad \text{and} \quad \bm{\varepsilon}^\mathrm{r}=\bm{\varepsilon}^\mathrm{r}_n+ \textstyle\sqrt\frac{3}{2}\beta\hat{\bm{r}}_n\sum_{s=1}^{n_\mathrm{y}}\Delta\kappa_s.
This expression represents the incremental version of the plastic yield function in weak form, from which the continuous equation is recovered by letting the term $\sqrt{3/2}\beta\Delta\bm{\sigma}:\hat{\bm{r}}_n$ vanish for small-enough time steps. To solve equation~\eqref{weak_p_model0} in a convenient way, we formulate a reduced problem in terms of the scalar field $\kappa_s$. The goal is to express, for each yield surface, the tensor-valued quantities as a function of $\kappa_s$, which is readily achieved for the ratcheting strain tensor as 
\begin{equation}
{\varepsilon}^\mathrm{r}(\bm{\kappa})= \bm{\varepsilon}^\mathrm{r}_n+\sqrt{\frac{3}{2}}\beta\hat{\bm{r}}_n\sum_{s=1}^{n_\mathrm{y}}\Delta\kappa_s.
\label{ratchparam}
\end{equation}
For the plastic strain tensor, we make use of standard arguments of $J_2$ plasticity to show that
\begin{equation}
{\bm{\varepsilon}}^\mathrm{p}_s(\kappa_s)={\bm{\varepsilon}}^\mathrm{p}_{sn}+\sqrt{\frac{3}{2}}\hat{\bm{n}}_s^\mathrm{trial}\Delta\kappa_s \quad \text{with} \quad \hat{\bm{n}}_s^\mathrm{trial}=\frac{\bm{s}^{\mathrm{p\,trial}}_{{s\,\mathrm{dev}}}}{\Vert \bm{s}^{\mathrm{p\,trial}}_{{s\,\mathrm{dev}}}\Vert}\equiv\hat{\bm{n}}_s \quad \text{and} \quad  \bm{s}^{\mathrm{p\,trial}}_{s\,\mathrm{dev}}= \bm{\sigma}_{s\,{\mathrm{dev}}}^\mathrm
{trial}-H^\mathrm{kin}_s\bm{\varepsilon}^\mathrm{p}_{sn},
\label{plastparam}
\end{equation}
where, for the multi-surface ratcheting model, the trial deviatoric stress reads
\begin{equation*}
\bm{\sigma}^\mathrm{trial}_{s\,{\mathrm{dev}}}\big(\bm{\varepsilon},\bm{\varepsilon}^\mathrm{p}_{j\neq s},\bm{\kappa},\alpha\big)=2 g(\alpha)\mu\big(\bm{\varepsilon}_\mathrm{dev}-(\bm{\varepsilon}^\mathrm{p}_{sn}+\textstyle\sum_{j\neq s}\bm{\varepsilon}^\mathrm{p}_{j})-\bm{\varepsilon}^\mathrm{r}(\bm{\kappa})\big).
\end{equation*}
%=\bm{\sigma}_{{\mathrm{dev}}}(\bm{\varepsilon},\bm{\varepsilon}^\mathrm{p}_{sn}+\textstyle\sum_{j\neq s}\bm{\varepsilon}^\mathrm{p}_{j},\bm{\varepsilon}^\mathrm{r},\alpha)
Note that the parametrizations~\eqref{ratchparam} and~\eqref{plastparam} must be subject to the irreversibility condition $\Delta\kappa_s\geq0$. Introducing these results in~\eqref{weak_p_model0} yields the following non-linear PDE in terms of $\kappa_s$:
\begin{equation}
\begin{aligned}
\int_{\Omega}\bigg[\bigg(\textstyle-\sqrt\frac{3}{2}\big\Vert\bm{s}^{\mathrm{p\,trial}}_{s\,\mathrm{dev}}\big\Vert  + g(\alpha)\big(\sigma^{\mathrm{p}}_s+H^\mathrm{iso}_s\kappa_s+(2\mu+H^\mathrm{kin}_s)\Delta\kappa_s\big)+\partial_{\kappa_s}I_{+}(\Delta\kappa_s)\bigg)\tilde{\kappa}_s  + g(\alpha)\eta_{\mathrm{p}s}^2\nabla\kappa_s\cdot\nabla\tilde{\kappa}_s\bigg]\,\mathrm{d}\bm{x} \ni 0,
\end{aligned}
\label{weak_p_model}
\end{equation}
where the indicator function $I_{+}:\mathbb{R}\to\mathbb{R}\cup\{+\infty\}$ is used to impose irreversibility. %In view of the multivalued function $\partial_{\kappa_s}I_{+}(\Delta\kappa_s)$, equation~\eqref{weak_p_model} approximates, in weak form, the plastic loading/unloading conditions in equation~\eqref{load_unload_model} at $t_{n+1}$, where the exact equations are recovered by letting the term $\sqrt{3/2}\beta\Delta\bm{\sigma}:\hat{\bm{r}}_n$ vanish for small-enough time steps. 

Equation~\eqref{weak_p_model} evaluated for all $s\in\mathrm{Y}$ yields a system of $n_{\mathrm{y}}$ constrained, non-linear equations with coupled $\kappa_s$. This system can be solved iteratively, for instance, using a fixed-point iteration scheme. Then, in each iteration, equation~\eqref{weak_p_model} is solved independently for the $s^\text{th}$ yield surface given the current estimates of $\{\bm{\varepsilon}_{j\neq s}^\mathrm{p},\kappa_{j\neq s}\}$. Clearly, for the single-surface case, a single non-linear equation must be~solved.

\paragraph*{Minimization with respect to the damage field} Given $\{\bm{u},\bm{\varepsilon}^\mathrm{p},\bm{\kappa},\bm{\varepsilon}^\mathrm{r}\}$, 
find 
\begin{equation}
\alpha=\arginf_{\alpha\in\mathscr{D}} \{\Pi(\bm{u},\bm{\varepsilon}^\mathrm{p},\bm{\kappa},\bm{\varepsilon}^\mathrm{r},\alpha)+\textstyle\int_{\Omega}I_{+}(\Delta\alpha)\,\mathrm{d}\bm{x}\}.
\label{mindam}
\end{equation}
The indicator function is used to impose irreversibility, while the box constraint $\alpha(\bm{x})\in[0,1]$ in~$\Omega$ must also be enforced for $\alpha\in\mathscr{D}$ to hold. Equation~\eqref{mindam} yields, for all $\tilde{\alpha}\in\tilde{\mathscr{D}}$, the necessary condition 
\begin{equation}
\begin{aligned}
\int_{\Omega}\bigg[\bigg(g'(\alpha)(\psi^{\mathrm{e_+}}+\psi^{\mathrm{p}}) + g'(\alpha)\sum_{s=1}^{n_\mathrm{y}}\sigma^{\mathrm{p}}_s\kappa_s + d(\gamma_n) w'(\alpha)+\partial_{\alpha}I_{+}(\Delta\alpha)\bigg)\tilde{\alpha} + d(\gamma_n)\eta_\mathrm{d}^2\nabla\alpha\cdot\nabla\tilde{\alpha}\bigg]\,\mathrm{d}\bm{x} \ni 0,
\end{aligned}
\label{weak_d_model}
\end{equation} %\delta\Pi(\bm{u},\bm{\varepsilon}^\mathrm{p},\bm{\kappa},\bm{\varepsilon}^\mathrm{r},\alpha)(\tilde{\bm{u}},\mathbf{0},\bm{0},\bm{0},0)=
which recovers the damage loading/unloading conditions in equation~\eqref{load_unload_model}.

Equations~\eqref{weak_eq_model}, \eqref{weak_p_model} and~\eqref{weak_d_model} are suitable for spatial discretization using standard finite elements.  The non-linearities in the mechanical balance and plasticity equations~\eqref{weak_eq_model} and~\eqref{weak_p_model} are tackled at each alternate minimization iteration with a standard Newton scheme, for which the maximum number of iterations is generally much lower than the required number of alternate minimization iterations at the corresponding time step. On the other hand, the constraints present in the plasticity and damage equations~\eqref{weak_p_model} and~\eqref{weak_d_model} can be tackled using techniques for PDE-constrained optimization. Examples of such techniques in the context of phase-field fracture modeling are outlined by~\citet{gerasimov2019}. Herein, we apply a simple algorithmic procedure adopted in previous works~\citep{lancioni2009, UllRodSam2016}, where unconstrained equations are first solved and an a posteriori correction is applied to the solution.% require a delicate treatment. Most studies on this topic focus on the damage irreversibility and box constraints for brittle fracture~\citep{MieHofWel2010,gerasimov2019}. Herein, we apply the simple algorithmic procedure adopted in references~\citep{UllRodSam2016} and \citep{rodriguez2018}.

%{\color{blue} Concerning the spatial discretization, sufficiently fine meshes are  generally required when strain localization is expected. In particular, similar to~\cite{MieTeiAld2016}, the plastic and damage characteristic lengths, denoted respectively as $\ell_\mathrm{p\,}$ and $\ell_\mathrm{d}$, can be computed as}

\section{Numerical simulations}\label{numsym}

%%% HERE. Run 1D simulations without damage and include those graphs as well.

This section presents numerical simulations that highlight the main features of the model described in section~\ref{model}. In order to highlight the various dissipative mechanisms and their interplay, the homogeneous 1D uniaxial response is first studied. Then, 2D finite element simulations are performed under plane strain conditions. As done for the 1D case, the first 2D simulation addresses the response of the cyclic plasticity model without damage. This example allows to describe a variety of plastic responses with cyclic effects. Then, examples involving fatigue crack growth coupled to cyclic plasticity are presented, capturing the initiation and propagation of ductile fatigue cracks.

\subsection{Homogeneous uniaxial response}\label{sec:hom}

The aim of this example is to provide an interpretation of the failure mechanisms that result from coupling cyclic plasticity to damage with fatigue effects. For illustrative purposes, we study the response of a single 1D element under either force loading or displacement loading. To this end, consider a straightforward reformulation of the multidimensional model presented in section~\ref{model} to the 1D case, where all vector- and tensor-valued quantities are replaced by scalar quantities, and a homogeneous response is assumed. With an obvious change of notation, the plasticity and damage yield functions~\eqref{fp} and~\eqref{fd} become
\begin{equation}
\begin{aligned}
f^\mathrm{p}_s &= \big\vert\sigma-g(\alpha)H^\mathrm{kin}_s\varepsilon^\mathrm{p}_s\big\vert - g(\alpha)\big(\sigma^{\mathrm{p}}_s+H^\mathrm{iso}_s\kappa_s\big),\\
f^\mathrm{d}&= -g'(\alpha)\big(\psi^{\mathrm{e}}+\psi^{\mathrm{p}}\big) - g'(\alpha)\sum_{s=1}^{n_\mathrm{y}}\sigma^{\mathrm{p}}_s\kappa_s-d(\gamma)w'(\alpha),
\end{aligned}
\end{equation} 
where, denoting the Young's modulus by $E$:
\begin{equation}
\begin{aligned}
\sigma &= g(\alpha)E\bigg(\varepsilon-\sum_{s=1}^{n_\mathrm{y}}\varepsilon^\mathrm{p}_s-\varepsilon^\mathrm{r}\bigg),\quad
\psi^\mathrm{e} = \frac{1}{2}E\bigg(\varepsilon-\sum_{s=1}^{n_\mathrm{y}}\varepsilon^\mathrm{p}_s-\varepsilon^\mathrm{r}\bigg)^2, \quad
\psi^\mathrm{p} = \frac{1}{2}\sum_{s=1}^{n_\mathrm{y}}\big(H^\mathrm{iso}_s\kappa_s^2+H^\mathrm{kin}_s{\varepsilon^\mathrm{p}_s}^2\big).
\end{aligned}
\end{equation}

\begin{table}[!ht]
\centering
\caption{Fixed parameters for the homogeneous uniaxial responses, with varying  parameters shown in table~\ref{tab:fp2}. $\sigma_s^\mathrm{p}$ and $H^\mathrm{kin}_s$ vary linearly from $s=1$ to $s=n_\mathrm{y}$.}
\small \begin{tabular}{lccccccccccc}
\hline
\multirow{2}{*}{Load type}                     & \multirow{2}{*}{Model} & $E$ & $w_0$ & $\gamma_0$ & $k$ & $n_{\mathrm{y}}$ & $\sigma^\mathrm{p}_1$ & $\sigma^\mathrm{p}_{n_\mathrm{y}}$ & $H^\mathrm{kin}_1$ & $H^\mathrm{kin}_{n_\mathrm{y}}$  \\
                                                  &                             & {[}MPa{]}          & {[}MPa{]}          & {[}MPa{]} & {[}-{]} & {[}-{]} & {[}MPa{]} & {[}MPa{]} & {[}MPa{]} & {[}MPa{]} \\ \hline
Force                                                    & AT-2  & 10  & $\infty$ / 260    & $\infty$ / 10         & 0.4 & 20               & 0.6                                                                     & 1.4                                & 100                & 9.09                            \\ \hline
Displacement                                             & AT-1  & 1   & $\infty$ / 750 / 75     & $\infty$ / 10     / 1  & 0.4 & 10               & 0.4                                                                     & 0.7                                & 8                  & 0.73    
\\ \hline                       
\end{tabular}
\label{tab:fp1}
\end{table}

\begin{table}[!ht]
\centering
\caption{Varying parameters for the homogeneous uniaxial responses, with fixed parameters shown in table~\ref{tab:fp1}. $H^\mathrm{kin}_s$ varies linearly from $s=1$ to $s=n_\mathrm{y}$, with $n_{\mathrm{y}}=20$ for force loading and $n_{\mathrm{y}}=10$ for displacement loading.}
\small \begin{tabular}{lcccc} \hline
\multirow{2}{*}{Load type}                     & \multirow{2}{*}{Response}  & $H^\mathrm{iso}_1$ & $H^\mathrm{iso}_{n_\mathrm{y}}$ & $\beta$ \\ &                             & {[}MPa{]}          & {[}MPa{]}   & {[}-{]} \\ \hline
\multirow{6}{*}{Force}                            & KH       & 0                  & 0                               & 0   \\
                                                  & KH-IH    & 0.2                & 0.0182                          & 0   \\
                                                  & KH-IS    & -0.08              & -0.0073                         & 0   \\
                                                  & KH-R     & 0                  & 0                               & 0.5 \\
                                                  & KH-IH-R  & 0.2                & 0.0182                          & 0.5 \\
                                                  & KH-IS-R  & -0.08              & -0.0073                         & 0.5 \\ \hline
\multicolumn{1}{c}{\multirow{6}{*}{Displacement}} & KH       & 0                  & 0                               & 0   \\
\multicolumn{1}{c}{}                              & KH-IH    & 0.02               & 0.0018                          & 0   \\
\multicolumn{1}{c}{}                              & KH-IS    & -0.015             & -0.0014                         & 0   \\
\multicolumn{1}{c}{}                              & KH-R     & 0                  & 0                               & 0.2 \\
\multicolumn{1}{c}{}                              & KH-IH-R  & 0.02               & 0.0018                          & 0.2 \\
                                                  & KH-IS-R  & -0.015 / -0.075             & -0.0014 / -0.0068 & 0.2 \\ \hline
\label{tab:fp2}
\end{tabular}
\end{table}

\begin{figure}[!h]
\centering
\includegraphics[scale=0.69, trim={0cm 0cm 0cm 0cm},clip]{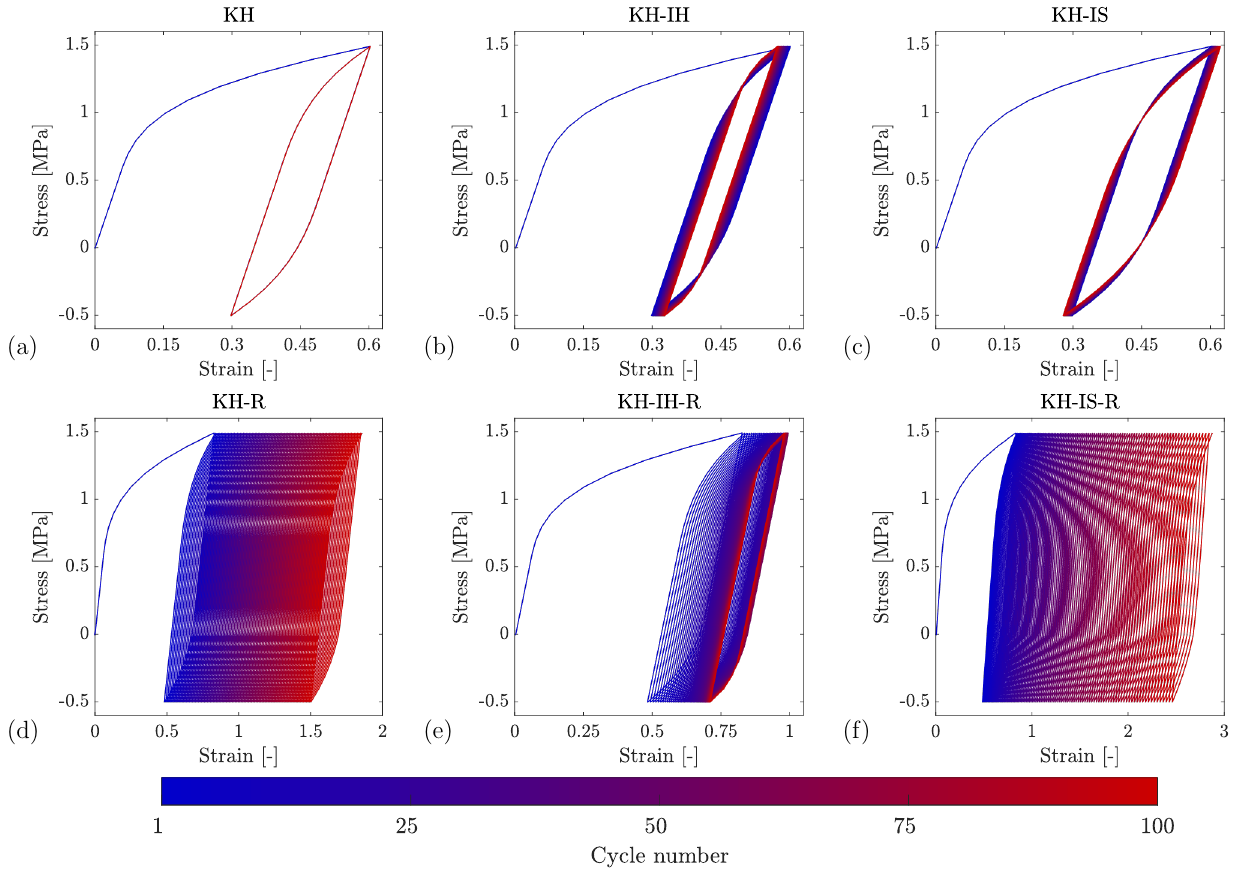}
\caption{Stress-strain curves of the homogeneous uniaxial tests under force loading, showing the undamaged response for (a) KH, (b) KH-IH, (c) KH-IS, (d) KH-R, (e) KH-IH-R and (f) KH-IS-R.}
\label{fpd01}
\end{figure}

\begin{figure}[!h]
\centering
\includegraphics[scale=0.69, trim={0cm 0cm 0cm 0cm},clip]{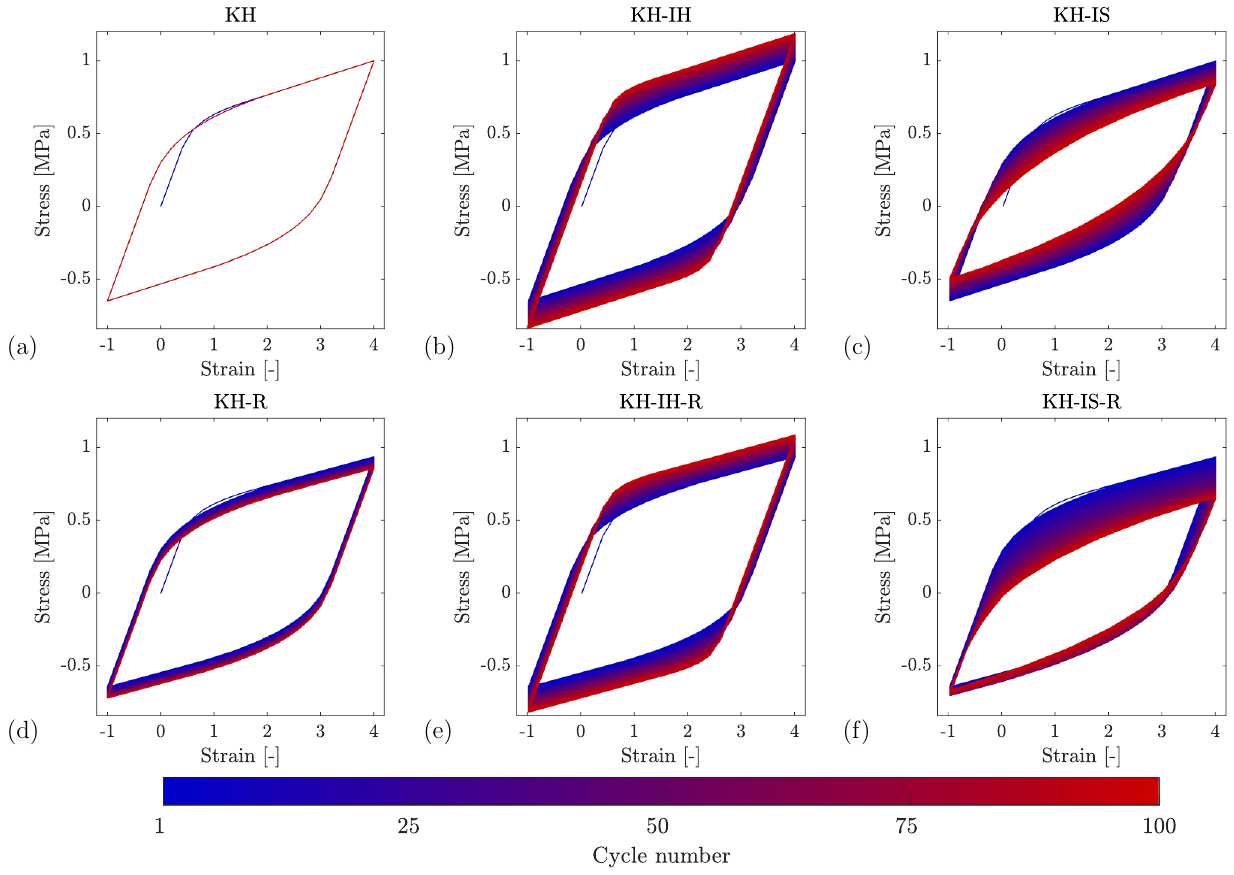}
\caption{Stress-strain curves of the homogeneous uniaxial tests under displacement loading, showing the undamaged response for (a) KH, (b) KH-IH, (c) KH-IS, (d) KH-R, (e) KH-IH-R and (f) KH-IS-R.}
\label{fpd02}
\end{figure}

\begin{figure}[!h]
\centering
\includegraphics[scale=0.69, trim={0cm 0cm 0cm 0cm},clip]{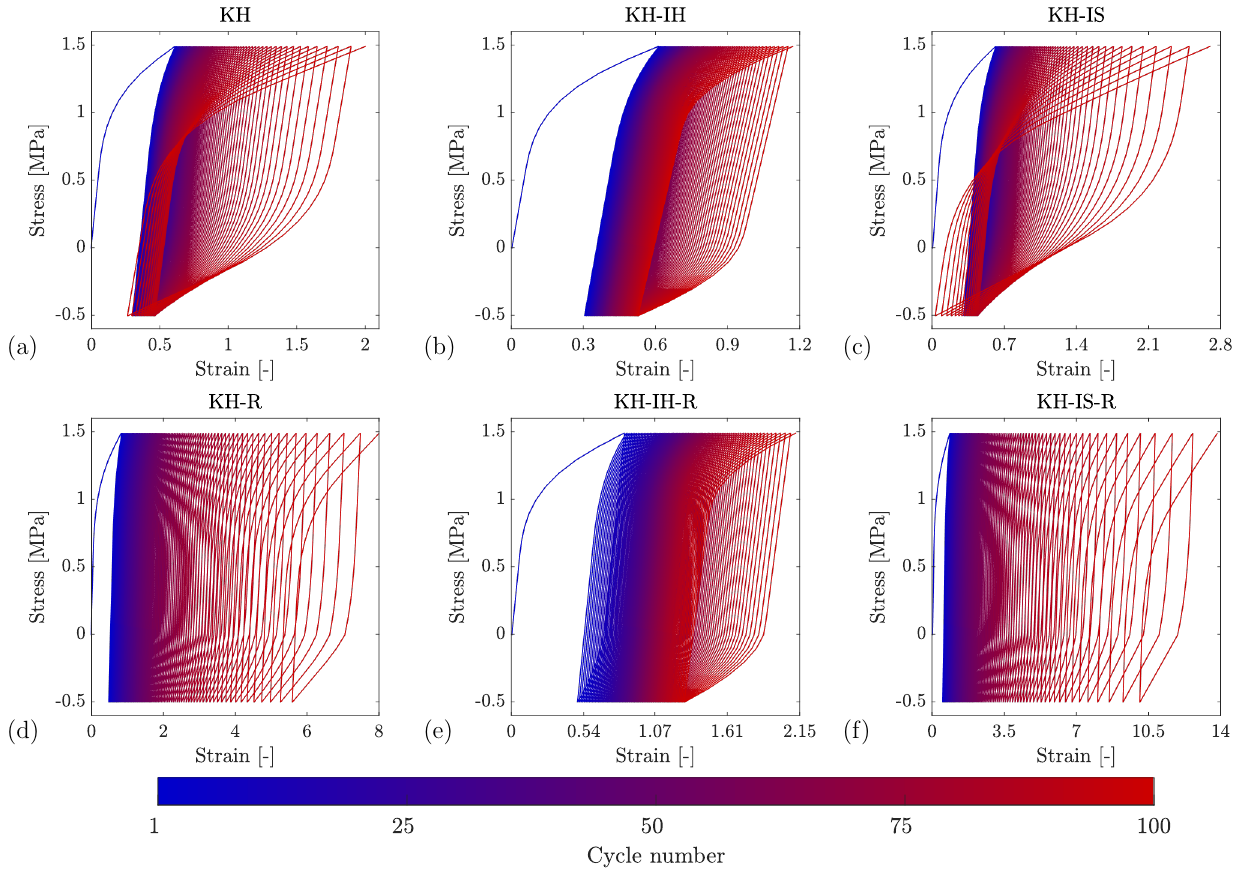}
\caption{Stress-strain curves of the homogeneous uniaxial tests under force loading, showing (a) KH, (b) KH-IH, (c) KH-IS, (d) KH-R, (e) KH-IH-R and (f) KH-IS-R coupled to damage with fatigue effects.}
\label{fpd1}
\end{figure}

\begin{figure}[!h]
\centering
\includegraphics[scale=0.69, trim={0cm 0cm 0cm 0cm},clip]{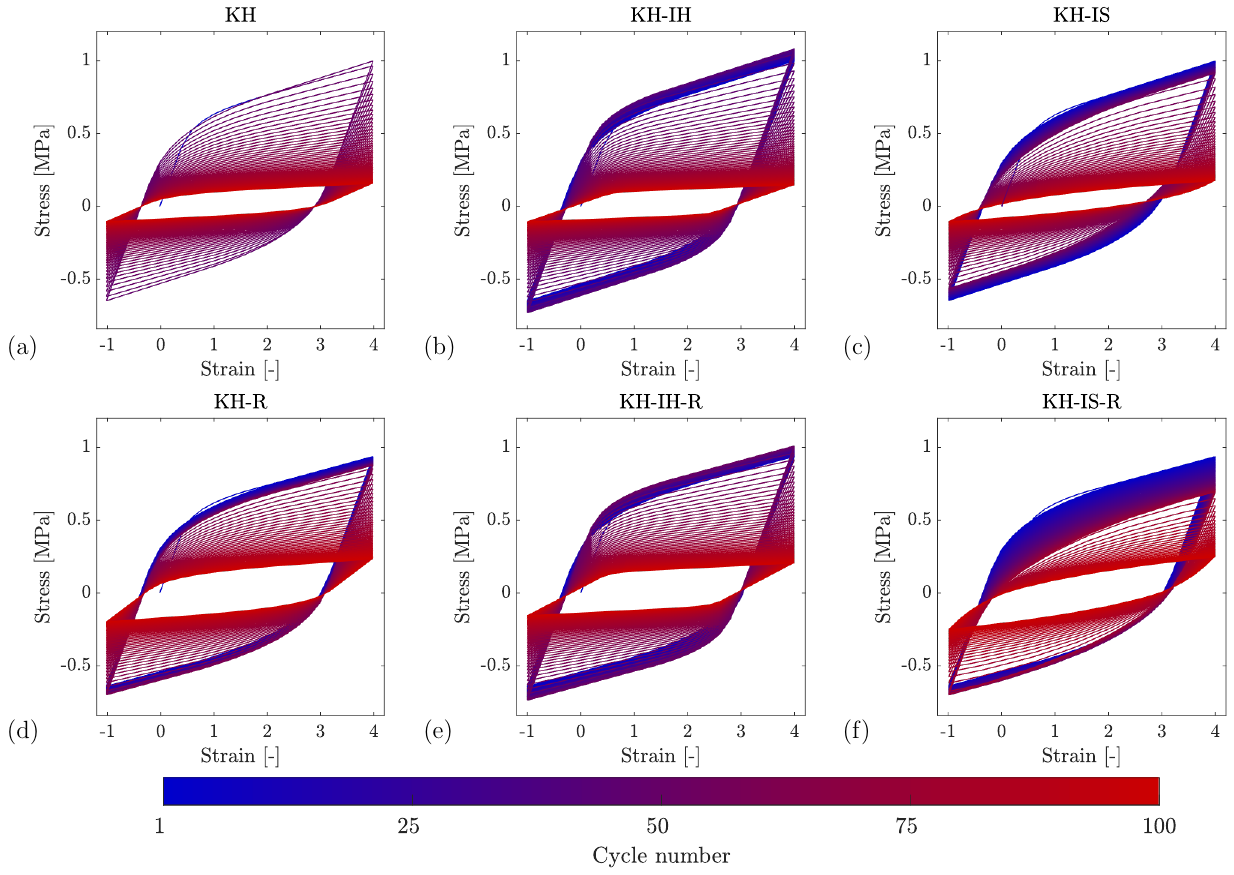}
\caption{Stress-strain curves of the homogeneous uniaxial tests under displacement loading, showing (a) KH, (b) KH-IH, (c) KH-IS, (d) KH-R, (e) KH-IH-R and (f) KH-IS-R coupled to damage with fatigue effects.}
\label{fpd2}
\end{figure}

To capture a smooth stress-strain response, the 1D element is modeled with multi-surface plasticity ($n_\mathrm{y}>1$). The parameters shown in tables~\ref{tab:fp1} and~\ref{tab:fp2} are considered to trigger a variety of behaviors, combining the purely  kinematic response (KH) with isotropic hardening (KH-IH); isotropic softening (KH-IS); ratcheting (KH-R); isotropic hardening and ratcheting (KH-IH-R); and isotropic softening and ratcheting (KH-IS-R). The loading consists of imposed stress cycles between $-0.5$~MPa and $1.5$~MPa for force loading and imposed strain cycles between -1 and 4 for displacement loading.

Figure~\ref{fpd01} shows the undamaged responses under force loading, ensured by setting $w_0=\infty$ and $\gamma_0=\infty$. The KH response exhibits closed hysteresis loops, while including isotropic hardening (KH-IH) results in a progressive decrease in cyclic strain amplitude. The opposite occurs with isotropic softening (KH-IS), where the cyclic strain amplitude progressively increases as the size of the yield surface decreases. The combination of these responses with ratcheting effects results in more complex cyclic evolutions. KH-R exhibits mean cyclic strain increments at a constant rate, i.e., a purely ratcheting response. Combining ratcheting with isotropic hardening (KH-IH-R) results in a competition between both mechanisms, where the ratcheting effect tends to vanish as the equivalent plastic strains increase. As a result of the loading pattern, this response occurs asymmetrically, with more pronounced cyclic hardening during tensile loading. On the other hand, KH-IS-R  leads to an asymmetric accelerated ratcheting response.

Figure~\ref{fpd02} shows the undamaged responses under displacement loading. As under force loading conditions, a closed cycle is observed for KH, while the yield surface progressively grows for KH-IH, leading to an elastic response after a sufficiently large number of cycles. The opposite occurs for KH-IS, where the size of the yield surface progressively decreases. In both KH-IH and KH-IS, the response is symmetric in tension and compression. On the other hand, for $\beta> 0$ (non-zero ratcheting strains), the response becomes asymmetric. This interesting effect of the ratcheting model allows to capture stress relaxation, as shown in figure~\ref{fpd02}(d) for KH-R, where the plastic cycles shift downwards. Moreover, an asymmetrical response is observed for both KH-IH-R and KH-IS-R, which combine stress relaxation with cyclic hardening and cyclic softening, respectively. For KH-IS-R, we have set $H^\mathrm{iso}_1=-0.015$~MPa and $H^\mathrm{iso}_{n_\mathrm{y}}=-0.0014$~MPa (table~\ref{tab:fp2}).

The result of coupling these responses to damage evolution is shown in figure~\ref{fpd1} for force loading, obtained by setting $w_0=260$~MPa and $\gamma_0=10$~MPa (table~\ref{tab:fp1}). In this case,  the initial plastic responses resemble the results in figure~\ref{fpd01}. However, as damage evolves, cyclic softening is triggered in all cases. Therefore, the softening responses are accelerated, while the (initially) hardening responses shift to a cyclic softening regime. An analogous result is observed in figure~\ref{fpd2} for the case of displacement loading with $w_0=750$~MPa and $\gamma_0=10$~MPa.

As previously discussed, the coupled plastic-damage model includes a fatigue mechanism that degrades the damage resisting force through $d(\gamma)$ as a function of free energy accumulation. Thus, in the absence of plastic strains, damage is accelerated for the AT-2 model and triggered after an initial elastic response for the AT-1 model, leading to a high-cycle fatigue process. Responses of this type are thoroughly studied in references~\citep{Alessi2017fatigue} and \citep{carrara2020}. On the other hand, the accumulation of plastic strains leads, on its own, to a plastic fatigue mechanism that promotes damage evolution through the plastic driving force $D^{\mathrm{p}}$. We associate this process with low-cycle fatigue, and the combination of $D^\mathrm{p}$ and $d(\gamma)$ with very-low-cycle fatigue. 

\begin{figure}[!ht]
\centering
\includegraphics[scale=0.68, trim={0cm 0cm 0cm 0cm},clip]{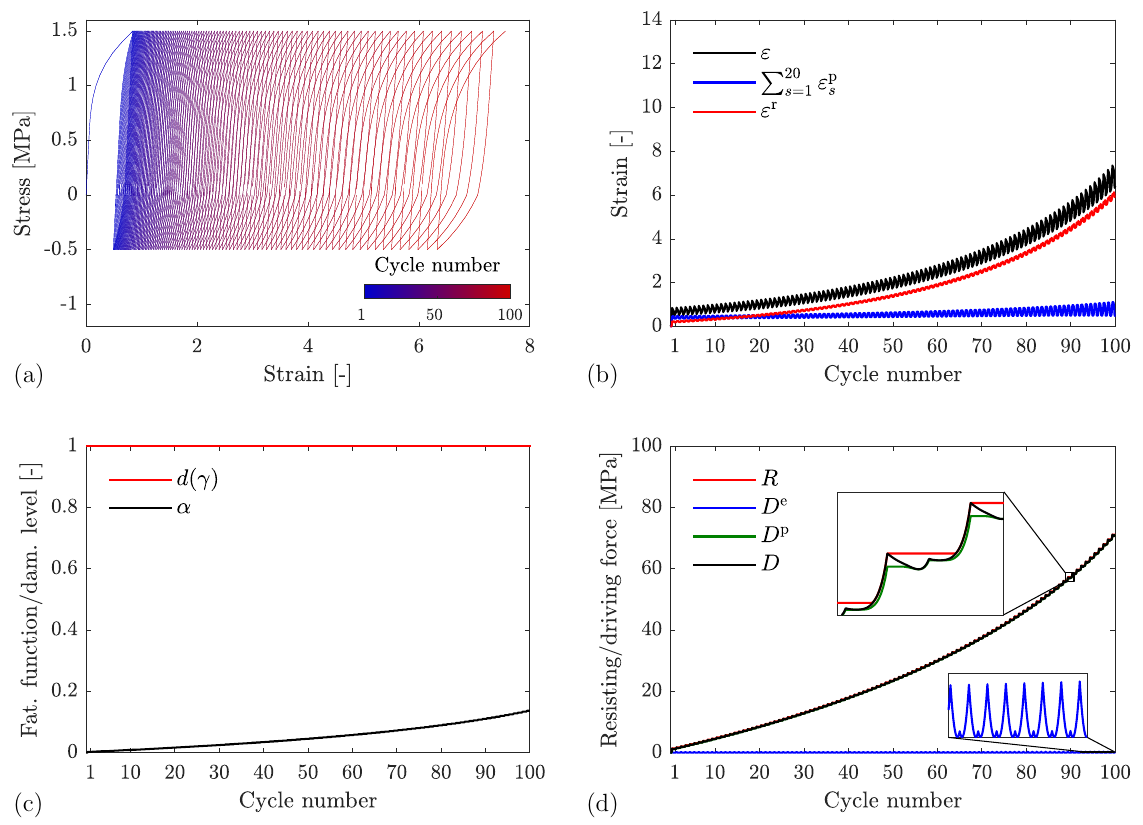}
\caption{KH-IS-R homogeneous response with damage ($\gamma_0=\infty$) under force loading: (a) stress-strain curve and corresponding time histories for the (b) strains, (c) damage and fatigue degradation and (d) damage driving and resisting forces.}
\label{fnogamma}
\end{figure}

\begin{figure}[!ht]
\centering
\includegraphics[scale=0.68, trim={0cm 0cm 0cm 0cm},clip]{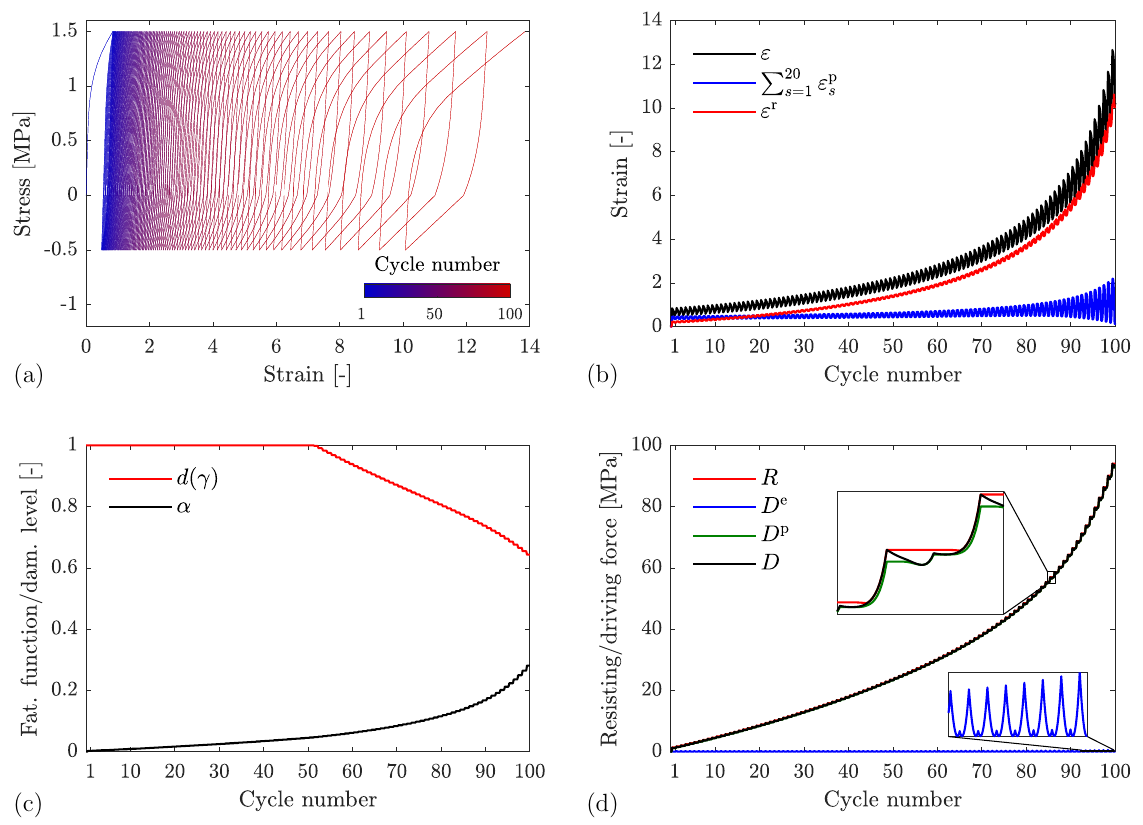}
\caption{KH-IS-R homogeneous response with damage ($\gamma_0=10$ MPa) under force loading: (a) stress-strain curve and corresponding time histories for the (b) strains, (c) damage and fatigue degradation and (d) damage driving and resisting forces.}
\label{fgamma}
\end{figure}

\begin{figure}[!ht]
\centering
\includegraphics[scale=0.68, trim={0cm 0cm 0cm 0cm},clip]{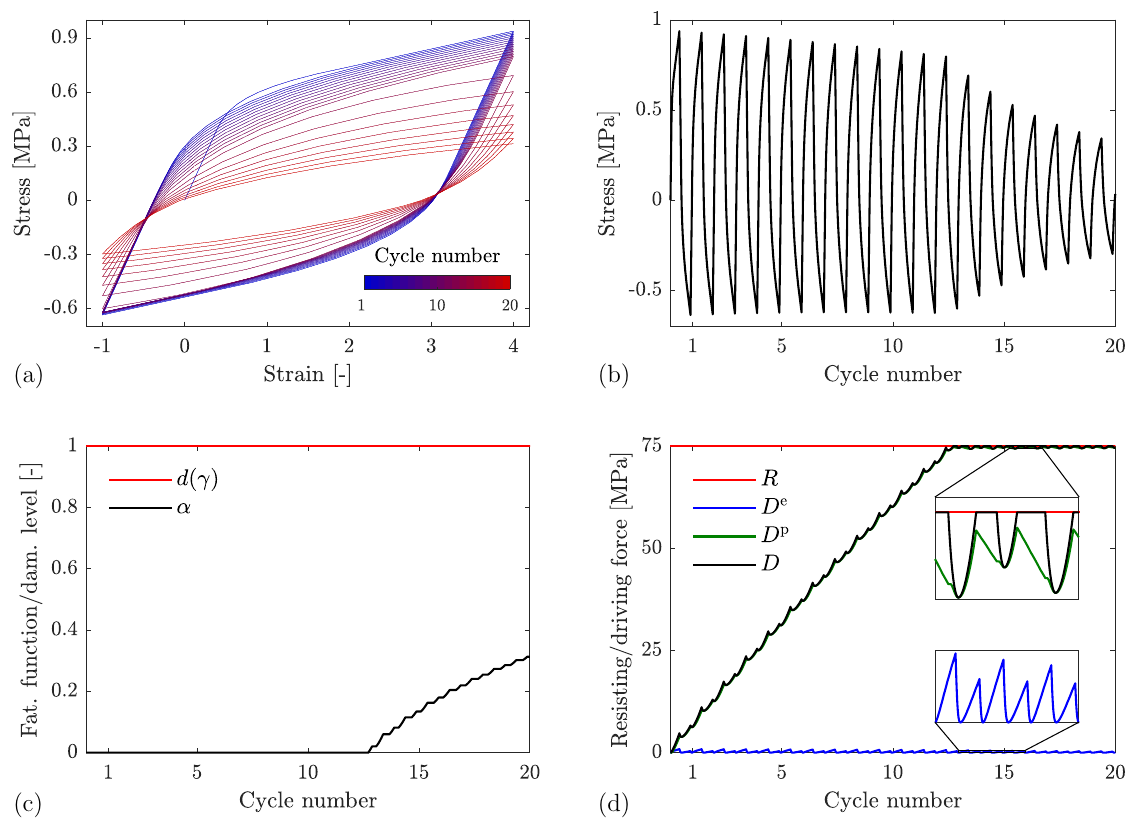}
\caption{KH-IS-R homogeneous response with damage ($\gamma_0=\infty$) under displacement loading: (a) stress-strain curve and corresponding time histories for the (b) stress, (c) damage and fatigue degradation and (d) damage driving and resisting forces.}
\label{dnogamma}
\end{figure}

\begin{figure}[!ht]
\centering
\includegraphics[scale=0.68, trim={0cm 0cm 0cm 0cm},clip]{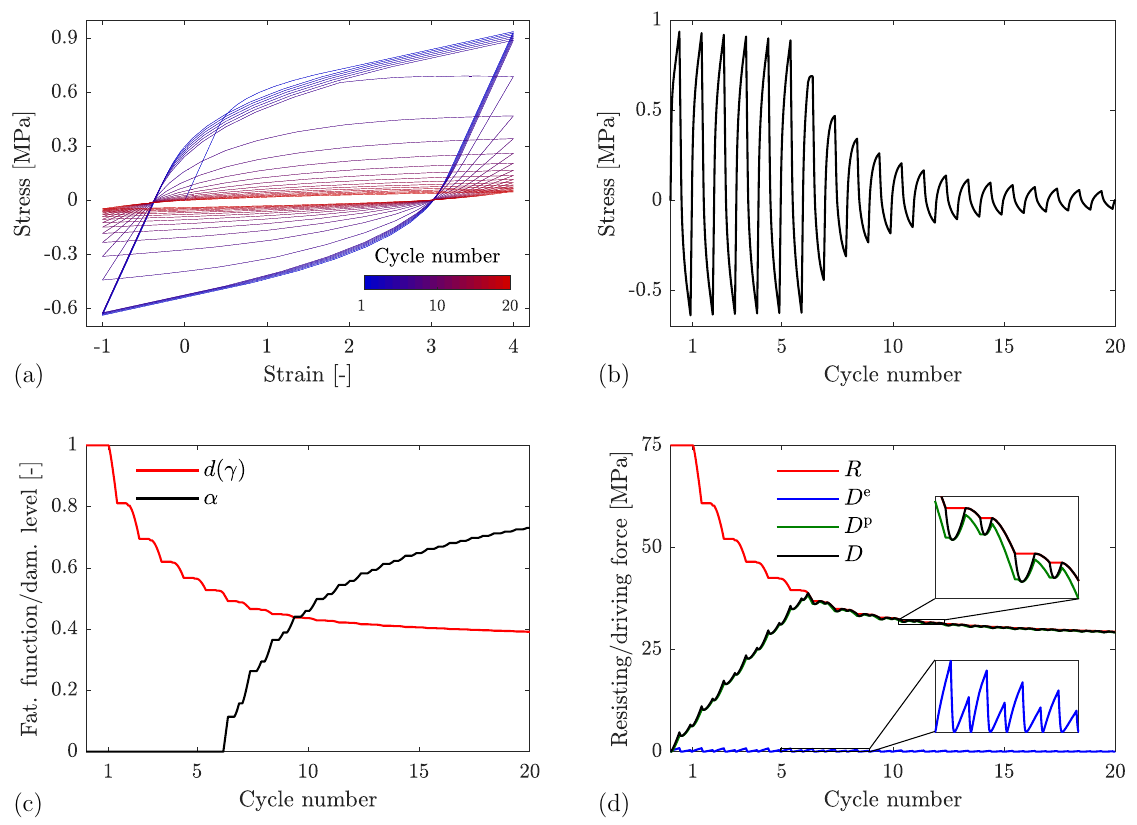}
\caption{KH-IS-R homogeneous response with damage ($\gamma_0=1$ MPa) under displacement loading: (a) stress-strain curve and corresponding time histories for the (b) stress, (c) damage and fatigue degradation and (d) damage driving and resisting forces.}
\label{dgamma}
\end{figure}

The fatigue mechanisms and their interplay are analyzed in figures~\ref{fnogamma} to~\ref{dgamma} for KH-IS-R. Figure~\ref{fnogamma} shows the mechanical response for the force loading test with the parameters from tables~\ref{tab:fp1} and~\ref{tab:fp2}, this time setting $\gamma_0=\infty$, such that damage is only driven by the accumulation of plastic strains. The ratcheting strain $\varepsilon^\mathrm{r}$ grows in an accelerated manner, where the mean cyclic value progressively increases, while the plastic strains $\varepsilon^\mathrm{p}_s$ grow in amplitude. These responses are reflected in the cyclic evolution of the total strain $\varepsilon$.  As a result of the AT-2 model, the damage resisting force $R=2w_0\alpha$ shown in figure~\ref{fnogamma}(d) grows as damage evolves. This occurs at every loading stage, where the resisting force intersects the total driving force $D$. Figure~\ref{fgamma} shows the counterpart simulation with $\gamma_0=10$~MPa. The ratcheting strain, the total strain and the plastic strains present the same evolution pattern, but grow at a notably higher rate.  This is a consequence of the accelerated growth of the damage variable that occurs as $d(\gamma)$~decreases. This response results from the combined effect of $D^\mathrm{p}$ and $d(\gamma)$, where the fatigue variable $\gamma$ is driven by the sum of elastic free energy and plastic free~energy.

For the case of displacement loading (figures~\ref{dnogamma} and~\ref{dgamma}), we take the values from tables~\ref{tab:fp1} and~\ref{tab:fp2} as $w_0=75$~MPa, $H^\mathrm{iso}_1=-0.075$~MPa and $H^\mathrm{iso}_{n_\mathrm{y}}=-0.0068$~MPa, in order to accelerate the softening response for illustrative purposes. We then set, alternatively, $\gamma_0=\infty$ and $\gamma_0=1$~MPa. Figure~\ref{dnogamma} shows the mechanical response of the displacement loading test with $\gamma_0=\infty$,  exhibiting a relatively slow stress decay. As a result of the AT-1 model, damage is triggered after 12 cycles. This occurs when the total driving force $D$, modulated by the plastic driving force $D^\mathrm{p}$, intersects the constant-valued resisting force $R$. Figure~\ref{dgamma} shows the counterpart with $\gamma_0=1$~MPa, where the resisting force begins to decrease after 2 cycles, causing damage to be triggered after only 6 cycles (figures~\ref{dgamma}(c) and~\ref{dgamma}(d)). Consequently, figure~\ref{dgamma}(b) shows a notably faster stress decay than figure~\ref{dnogamma}(b), as expected in very-low-cycle~fatigue. 

\FloatBarrier

\subsection{Non-homogeneous finite element simulations}\label{sec:coupled}

In this subsection, we present a series of 2D simulations under plane strain conditions that highlight the versatility of the proposed model. For simplicity, and to alleviate computational cost, we consider in all examples the case of single-surface plasticity, which is easily recovered from the model presented in section~\ref{model} by letting $n_\mathrm{y}=1$. For notational simplicity, the subscripts are thus dropped from the plasticity parameters, such that the isotropic hardening modulus, the kinematic hardening modulus, the yield strength and the characteristic length read $H^\mathrm{iso}$, $H^\mathrm{kin}$, $\sigma^\mathrm{p}$ and $\ell_\mathrm{p}$ (alternatively $\eta_\mathrm{p}$), respectively.

In all computations, bilinear quadrilateral elements are employed. The use of $\mathrm{C}^0$ continuous elements is feasible in the present model due to the gradient regularization employed for both plasticity and damage.

\subsubsection{Undamaged perforated specimen}\label{sec:plate}

\begin{figure}[!ht]
\centering
\includegraphics[width=0.55\textwidth]{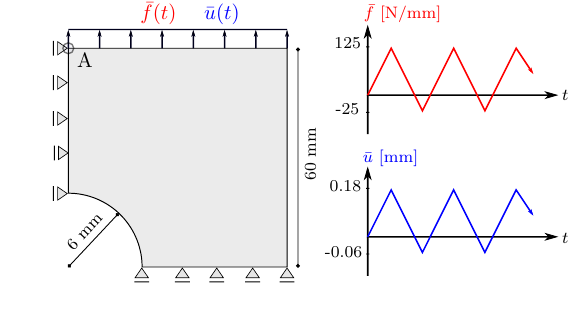}
%\def\svgwidth{0.55\textwidth}
%\import{}{holeinplate_scheme.pdf_tex}
\caption{Schematic representation of the perforated specimen under force loading and displacement loading.}
\label{fig:plast_setup}
\end{figure}

%The aim of this experiment is to present the range of possible cyclic responses embedded in the ratcheting plasticity model. Therefore, d
This example consists of a square specimen with a central hole under plane strain conditions and cyclic loading, alternatively subjected to force loading and displacement loading. Due to symmetry conditions, only the top-right quarter of the specimen is analyzed (figure~\ref{fig:plast_setup}), and uniform loading is applied on the top~border. Unless stated otherwise, a uniform mesh of 800 elements is considered in the simulations, with 120 (160) time steps per cycle for force (displacement) loading.

Damage is disallowed by setting, in the general model, a sufficiently large~$w_0$. We consider a Young's modulus $E=205\times10^3$ MPa, a Poisson's ratio $\nu=0.3$, a plastic yield strength $\sigma^\mathrm{p}=100$ MPa and an internal length scale $\eta_\mathrm{p}=0.6$ N$^{1/2}$. It is worth mentioning that the effect of the non-local term modulated by $\eta_\mathrm{p}$ has little effect for the stress hardening responses of this example. The remaining parameters vary according to table~\ref{tab:plast}, where the kinematic hardening modulus is combined with a positive (negative) isotropic hardening modulus to describe cyclic hardening (softening). Moreover, ratcheting and stress relaxation effects are attained by setting $\beta>0$. 
% $\eta_\mathrm{p}=0.6$ N$^{1/2}$

\begin{table}[!ht]
\centering
\caption{Varying parameters for the cyclically plastic responses of the perforated specimen.}
\small \begin{tabular}{lllll}
\hline
\multirow{2}{*}{Load type}                     & \multirow{2}{*}{Response}   & $H^\mathrm{kin}$ & $H^\mathrm{iso}$ & $\beta$     \\
                                                  &                             & {[}MPa{]}          & {[}MPa{]}          & {[}-{]} \\ \hline
\multirow{6}{*}{Force}                            & KH                          & 16513.89           & 0                  & 0       \\
                                                  & KH-IH                       & 16348.75           & 165.14             & 0       \\
                                                  & KH-IS                       & 16447.83           & -66.06             & 0       \\
                                                  & KH-R                        & 16513.89           & 0                  & 0.075     \\
                                                  & KH-IH-R                     & 16348.75           & 165.14             & 0.1     \\
                                                  & KH-IS-R                     & 16447.83           & -66.06             & 0.1     \\ \hline
\multicolumn{1}{c}{\multirow{5}{*}{Displacement}} & KH                          & 16513.89           & 0                  & 0       \\
\multicolumn{1}{c}{}                              & KH-IH                       & 16500.68           & 13.21              & 0       \\
\multicolumn{1}{c}{}                              & KH-IS                       & 16500.68           & -13.21             & 0       \\
\multicolumn{1}{c}{}                              & KH-R                        & 16513.89           & 0                  & 0.4     \\
\multicolumn{1}{c}{}                              & KH-IH-R                     & 16500.68           & 13.21              & 0.8     \\
                                                  & KH-IS-R                     & 16500.68           & -13.21             & 0.8     \\ \hline
\end{tabular}
\label{tab:plast}
\end{table}

\begin{figure}[!ht]
\centering
\includegraphics[scale=0.75, trim={0cm 0cm 0cm 0cm},clip]{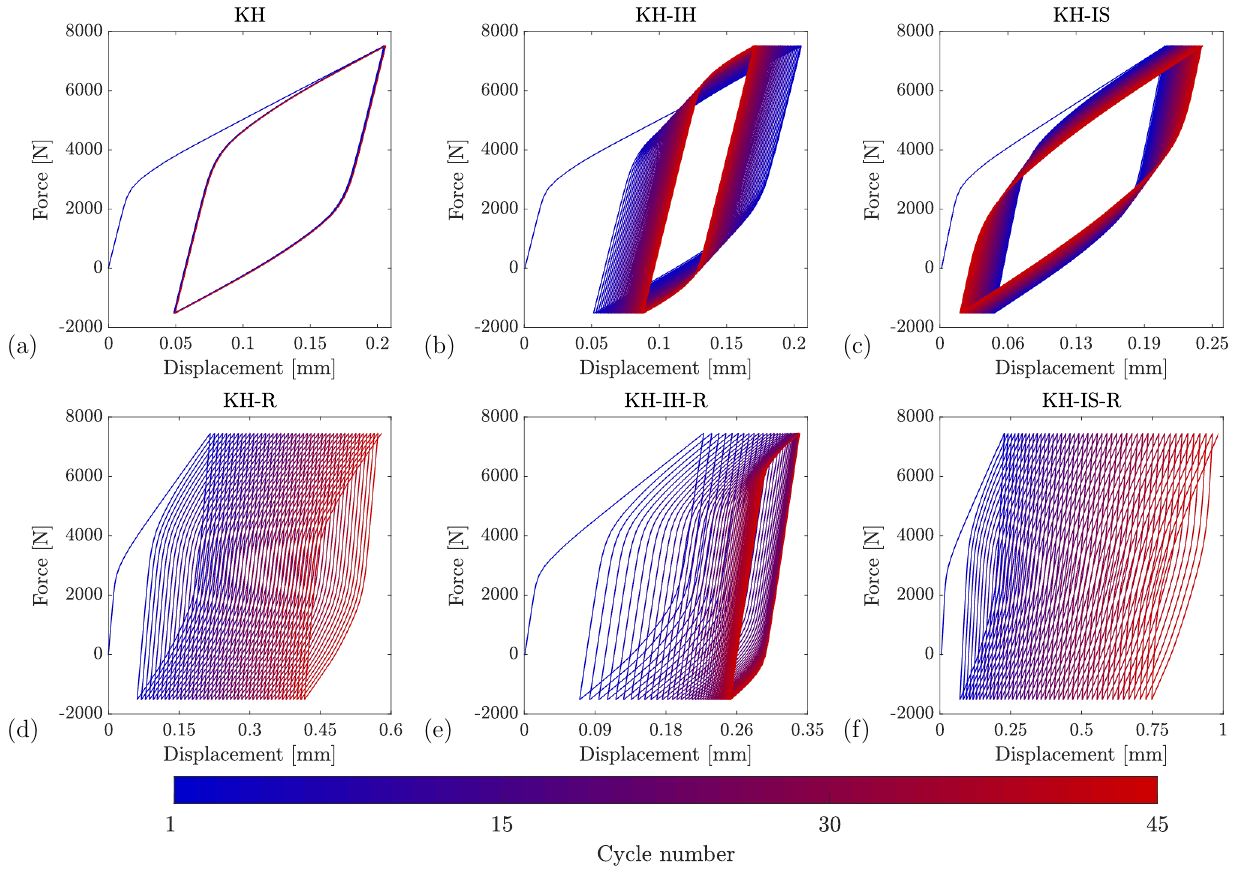}
\caption{Force-displacement curves measured at location A of the perforated specimen under force loading (figure~\ref{fig:plast_setup}), showing the cyclic response with  (a) KH, (b) KI-IH, (c) KH-IS, (d) KH-R, (e) KH-IH-R and (f) KH-IS-R.}
\label{pf1}
\end{figure}

\begin{figure}[!ht]
\centering
\includegraphics[scale=0.75, trim={0cm 0cm 0cm 0cm}]{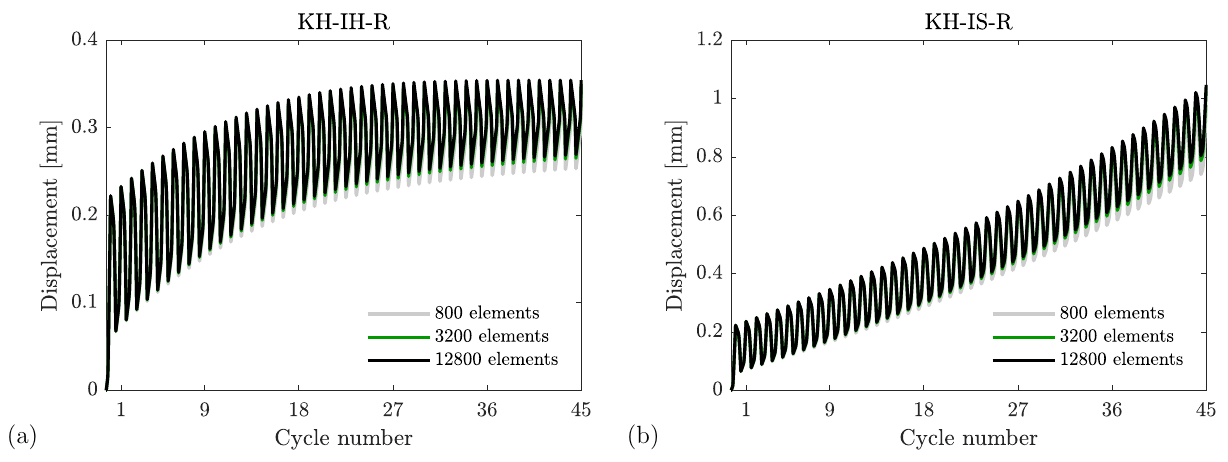}
\caption{Time history of the displacement at location A for the perforated specimen under force loading (figure~\ref{fig:plast_setup}) computed with different mesh densities, showing (a) the KH-IH-R response and (b) the KH-IS-R response.}
\label{pf3}
\end{figure}

Figure~\ref{pf1} shows the cyclic responses under force loading, corresponding to the force-displacement curves in the vertical direction (with displacements measured at location A) due to a uniform distributed force of magnitude $\bar{f}(t)$ (figure~\ref{fig:plast_setup}). The results resemble and further highlight the behaviors presented for the homogeneous uniaxial case: the KH  response exhibits closed hysteresis cycles, capturing the Bauschinger effect; KH-IH results in a progressive decrease in cyclic displacements, leading to a closed  plastic loop; and KH-IS leads to a progressive increase in cyclic displacements as the size of the yield surface decreases. Likewise, KH-R exhibits ratcheting at a constant rate;  KH-IH-R results in a vanishing ratcheting response (i.e., a \emph{shakedown} response) occurring asymmetrically in tension and compression; and KH-IS-R results in an asymmetric accelerated ratcheting response. The last two responses are further observed in figure~\ref{pf3}, which shows the time history of the displacements at location A for KH-IH-R and KH-IS-R. The results are shown for different mesh densities, indicating mesh-objectivity.

\begin{figure}[!ht]
\centering
\includegraphics[scale=0.75, trim={0cm 0cm 0cm 0cm},clip]{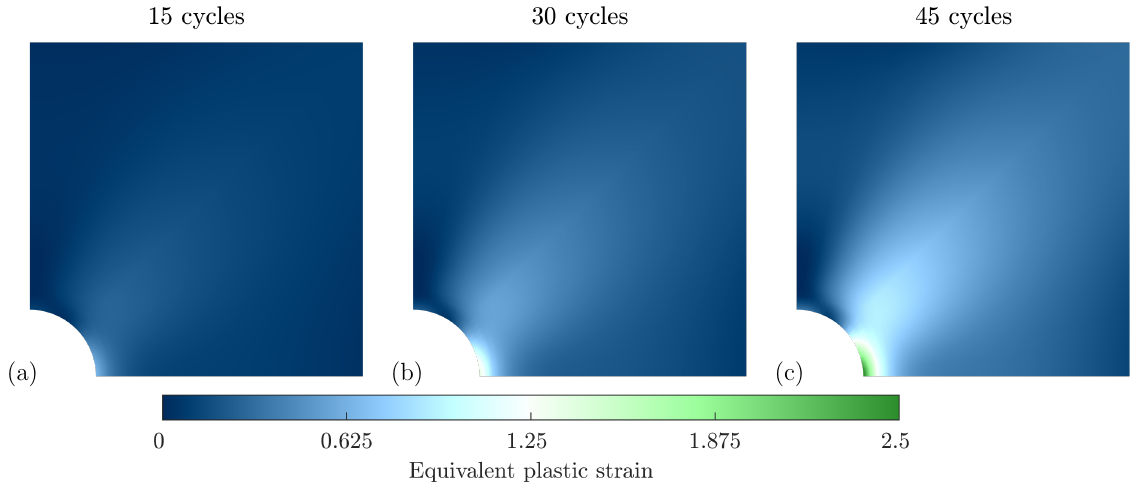}
\caption{Equivalent plastic strain for the KH-IS-R response of the perforated specimen under force loading after (a) 15, (b) 30 and (c) 45 cycles.}
\label{pf2}
\end{figure}

Figure~\ref{pf2} shows the equivalent plastic strains at different load cycles for KH-IS-R, while the contour plots for the other responses are not shown due to their qualitative similarity. The plastic strains concentrate near the hole and propagate in an inclined pattern, as expected for deviatoric-driven plasticity. %where the non-local effect governed by the plastic length scale allows to control the diffuseness of the spatial distribution.    

\begin{figure}[!ht]
\centering
\includegraphics[scale=0.8, trim={0cm 0cm 0cm 0cm},clip]{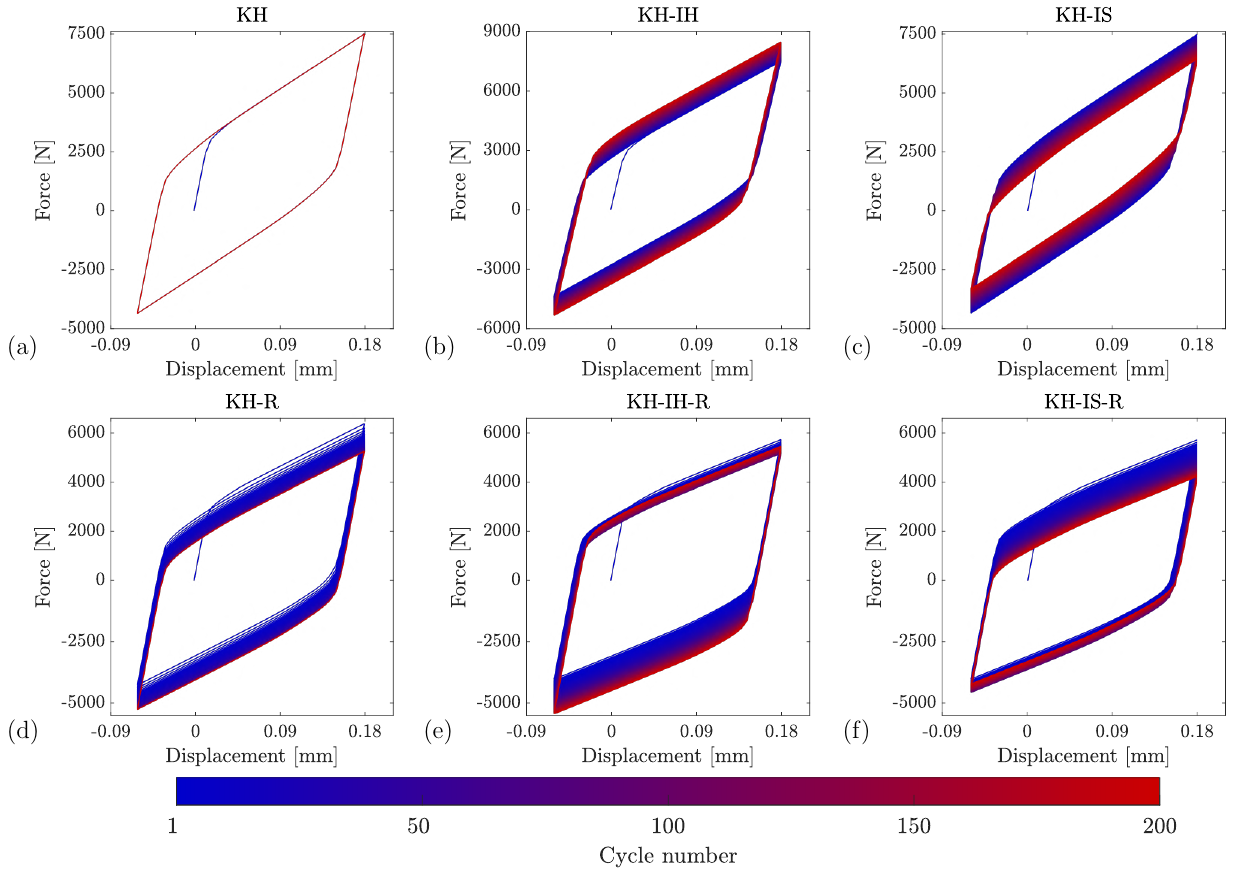}
\caption{Force-displacement curves measured at the top border of the perforated specimen under displacement loading (figure~\ref{fig:plast_setup}), showing the cyclic response with (a) KH, (b) KI-IH, (c) KH-IS, (d) KH-R, (e) KH-IH-R and (f) KH-IS-R.}
\label{pd1}
\end{figure}

\begin{figure}[!ht]
\centering
\includegraphics[scale=0.75, trim={0cm 0cm 0cm 0cm},clip]{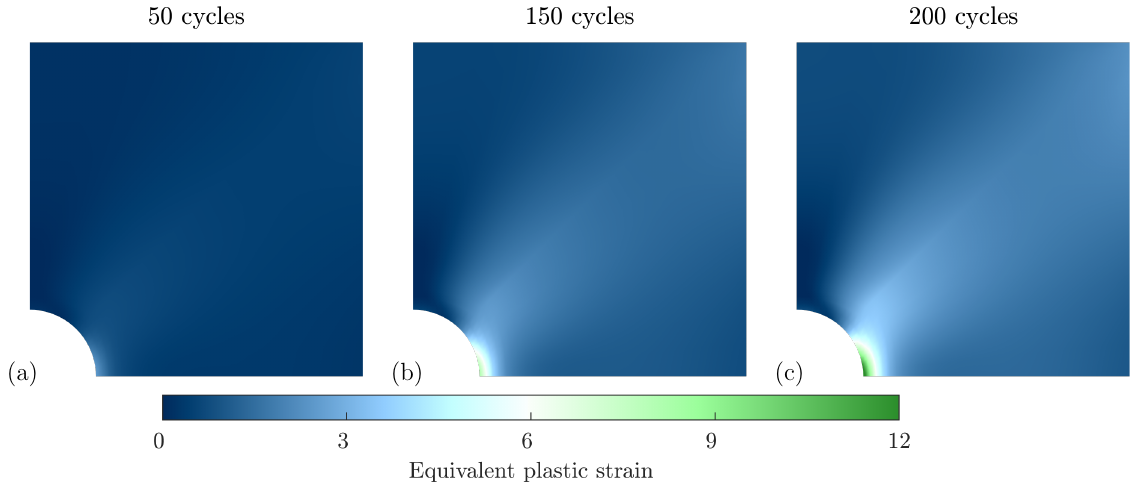}
\caption{Equivalent plastic strain for the KH-IS-R response of the perforated specimen under displacement loading after (a) 50, (b) 150 and (c) 200 cycles.}
\label{pd2}
\end{figure}

Figure~\ref{pd1} shows the cyclic responses under displacement loading of magnitude $\bar{u}(t)$ applied on the top border (figure~\ref{fig:plast_setup}). As for the case of force loading, the results resemble and further highlight the behaviors presented for the homogeneous uniaxial case, with a closed cycle observed for KH; a  progressively growing yield surface for KH-IH; and a progressively vanishing yield surface for KH-IS. Likewise, including ratcheting effects leads to stress relaxation for KH-R, and pronounced asymmetric responses for KH-IH-R and KH-IS-R with cyclic hardening and cyclic softening, respectively. Figure~\ref{pd2} shows the equivalent plastic strains at different load cycles for KH-IS-R, which qualitatively resemble the results obtained under force loading. The simulation was computed with 800, 3200 and 12800 elements. The results indicate similar mesh convergence as for the case of force loading and are thus not shown for brevity.

%{\color{blue} This behavior is further observed in figure~\ref{meshdisp}, where the time history of the reaction force at the top border is shown for KH-IH-R and KH-IS-R. Moreover, the results upon mesh refinement indicate mesh convergence.}
%This behavior is further observed in figure~\ref{meshdisp}, where the time history of the reaction force at the top border is shown for KH-IH-R and KH-IS-R. Moreover, the results upon mesh refinement indicate mesh convergence.
 
From this study, we conclude that the plastic responses observed in figures~\ref{pf1} and~\ref{pd1} encompass a wide range of material responses in a specimen with geometrical effects or imperfections. Of particular interest are the responses exhibiting ratcheting and stress relaxation, which are further enriched by including cyclic hardening or cyclic softening effects. This allows to recover complex cyclic behaviors often observed in experimental works~\cite{hassan1994a,hassan1994b}.

\subsubsection{Asymmetrically notched specimen}

\begin{figure}[!ht]
\centering
\includegraphics[width=0.55\textwidth]{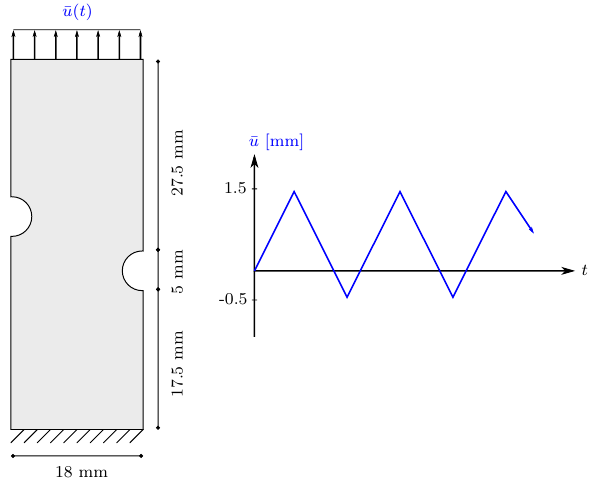}
%\def\svgwidth{0.55\textwidth}
%\import{}{asn_scheme.pdf_tex}
\caption{Schematic representation of the asymmetrically notched specimen under displacement loading.}
\label{plate2}
\end{figure}

This example aims to describe the initiation, growth and merging of ductile cracks in a low-cycle fatigue process, driven by the accumulation and localization of isochoric plastic deformations. For this purpose, we subject an asymmetrically notched specimen in plane strain conditions to displacement cycles, as schematically shown in figure~\ref{plate2}, where the bottom border is fixed and vertical displacements of magnitude $\bar{u}(t)$ are imposed on the top border, with 300 time steps per cycle. %We describe the plastic response with a single yield surface ($n_\mathrm{y}=1$) and kinematic hardening, where failure is triggered at a low number of load cycles. 

%%% Remove ny and put another row for the length scale study (test 1; test2), and add columns for umax, umin and AT model 
\begin{table}[!h]
\centering
\caption{Material parameters for the asymmetrically notched specimen under cyclic loading.}
\small \begin{tabular}{cccccccccccc}
\hline
$K$       & $\nu$   & $w_0$     & $\ell_\mathrm{d}$ & $\gamma_0$ & $k$  & $\sigma^\mathrm{p}$ & $H^\mathrm{iso}$ & $H^\mathrm{kin}$ & $\ell_{\mathrm{p}}$ & $\beta$ \\ 
{[}MPa{]} & {[}-{]} & {[}MPa{]} & {[}mm{]}   & {[}MPa{]}  & {[}-{]}         & {[}MPa{]}             & {[}MPa{]}          & {[}MPa{]}          & {[}mm{]}     & {[}-{]} \\ \hline
71659.46  & 0.331   &  1428.3       &         0.2       & 300         & 0.4   & 345                   & 112.5                  & 1591.67                & 0.25                   & 0.4    \\ \hline  
\end{tabular}
\label{tab:asn}
\end{table}

\begin{figure}[!ht]
\centering
\includegraphics[scale=0.80, trim={0cm 0cm 0cm 0cm},clip]{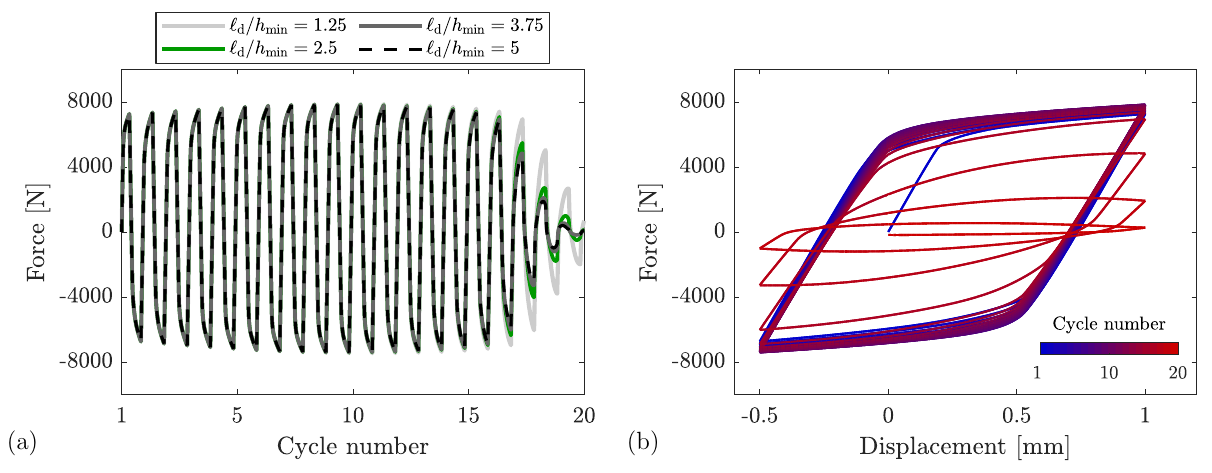}
\caption{Results of the asymmetrically notched specimen showing (a) the time history of the reaction force for different mesh sizes and (b) the force-displacement curve for $h_\mathrm{min}=0.053$ mm.}
\label{asn_ss}
\end{figure}

\begin{figure}[!ht]
\centering
\includegraphics[scale=0.70, trim={0cm 0cm 0cm 0cm},clip]{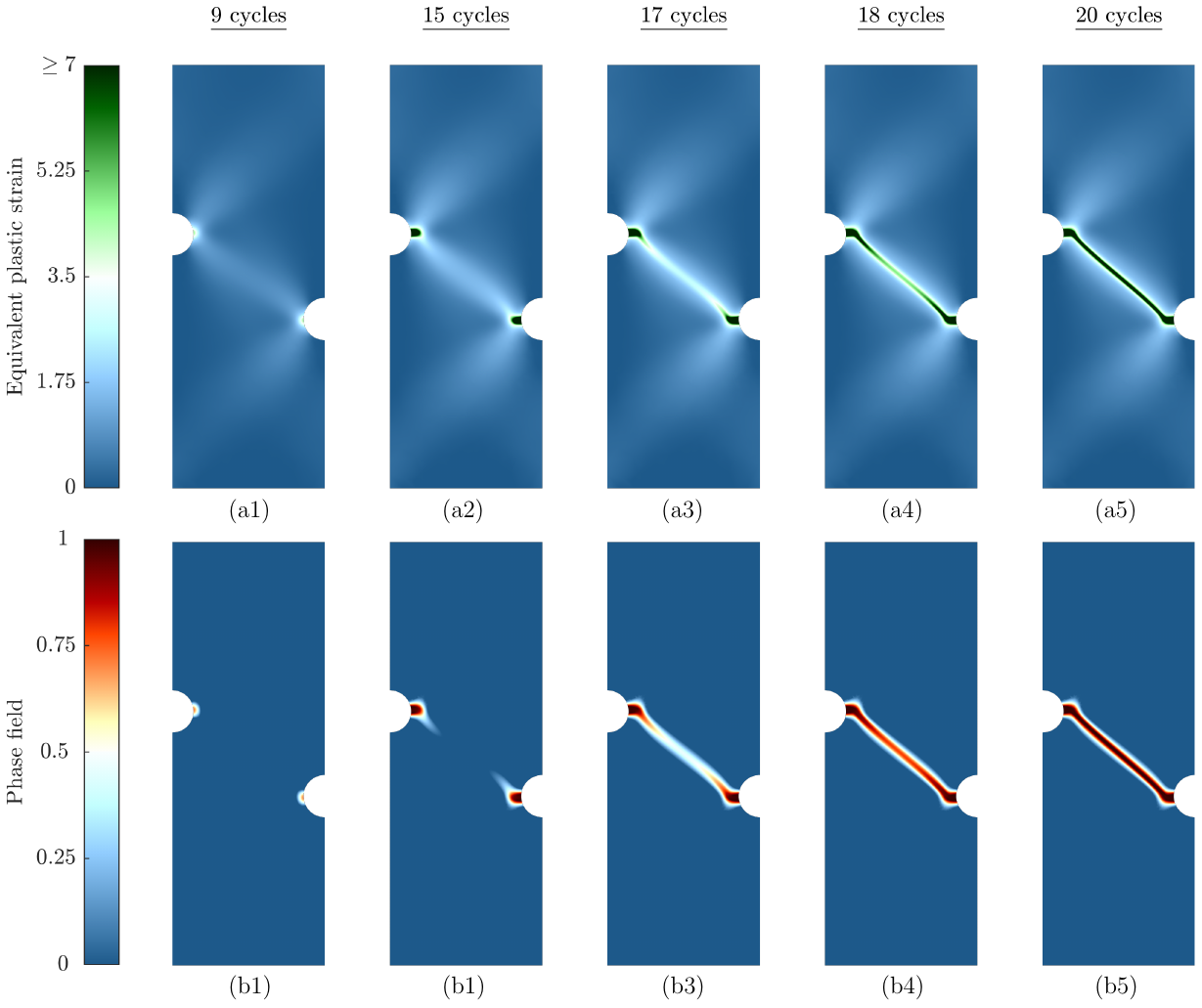}
\caption{Equivalent plastic strains (a1)-(a5) and phase field (b1)-(b5) during damage initiation (9 cycles), crack initiation (15 cycles), crack propagation (17-18 cycles), and final merging and fracture stage (20 cycles).}
\label{asn}
\end{figure}

%%% HERE!!!!!
% In case there's a Test 2: Consider first the parameters shown in table~\ref{tab:asn} corresponding to Test 1.
Consider the parameters shown in table~\ref{tab:asn} and the AT-1 damage model (equation~\eqref{w0}). For the numerical simulations, an unstructured mesh of bilinear quadrilateral elements was employed, with local refinement in the region where the crack is expected to develop.  In this region, a characteristic element size $h_\mathrm{min}\in\{0.16,0.08,0.053,0.04\}$ mm  was chosen, aiming for a compromise between (i) a sufficiently large ratio $\ell_\mathrm{d}/h_\mathrm{min}$ required to properly describe the phase-field fracture process~\citep{BourFrancMar2008,freddi2019} and (ii) the computational cost associated with the present cyclic simulations. In addition, note from table~\ref{tab:asn} that, in agreement with~\citet{miehe2017}, the characteristic length scales have been chosen such that $\ell_\mathrm{p} > \ell_\mathrm{d}>h_\mathrm{min}$, aiming for an objective and physically meaningful description of the coupled plasticity-damage evolution.  

Figure~\ref{asn_ss}(a) shows the time history of the reaction force for the different mesh sizes, where mesh-objective results are observed as $h_\mathrm{min}$ decreases with fixed characteristic lengths $\ell_\mathrm{p}$ and $\ell_\mathrm{d}$. Figure~\ref{asn_ss}(b) shows the force-displacement response measured at the top border of the specimen. Similar to the homogeneous uniaxial response, we observe that the initial hardening loops progressively  decrease in amplitude as damage evolves, with the global stiffness vanishing after 20 cycles.

Figure~\ref{asn} shows contour plots for the equivalent plastic strains and the phase-field/damage variable at different loading stages. After 9 cycles, the plastic strains are localized in relatively wide bands, governed by the plastic length scale. Within these regions, damage is triggered and begins to evolve. After about 15 cycles, individual cracks initiate and propagate from the notches, resulting in a strong localization of plastic strains. At this stage, the diffuse nature of the plastic strains in the pre-cracked states shifts to a strongly localized evolution, thus capturing the behavior of ductile cracks. At subsequent stages, the cracks propagate and finally merge along a central shear band, leading to a slip-like failure mode after 20 cycles.

In the context of the phase-field approach to ductile fracture, similar examples have been studied under monotonic loading~\citep{Duda2015,AmbGerDeL2015,rodriguez2018}. In particular, \citet{rodriguez2018} present different responses that depend on the degree of ductility: a mixed Mode I/II failure associated with an elastoplastic brittle response, and a Mode II-dominated failure associated with a ductile response. In the present study, the cyclic response is strongly driven by plastic strains and therefore resembles the ductile failure mode obtained for monotonic loading by~\citet{rodriguez2018}.

\begin{figure}[!h]
\centering
\includegraphics[scale=0.75, trim={0cm 0cm 0cm 0cm},clip]{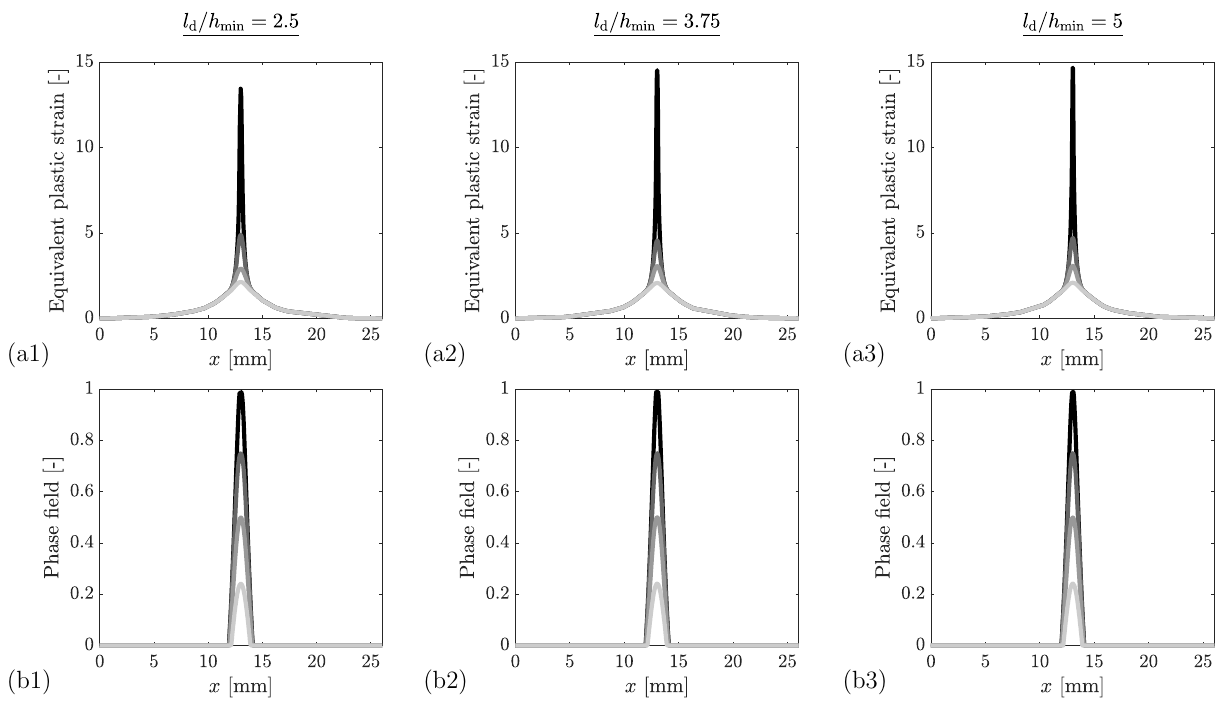}
\caption{Equivalent plastic strain (a1)-(a3) and damage (b1)-(b3) profiles  for different mesh sizes, showing mesh-objective results for both plasticity and damage. The results are shown when the maximum damage level reaches 0.25 (light gray), 0.5 (medium gray), 0.75 (dark gray) and 0.99 (black). The curves are plotted along a direction orthogonal to the crack in the center of the specimen, where $x$ denotes the distance along the cross-section.}
\label{asn_prof}
\end{figure}

An important feature of the proposed model is the non-local treatment of the localized responses for both plasticity and damage, governed by two length scale parameters $\ell_\mathrm{p}$ and $\ell_\mathrm{d}$. As discussed by~\citet{miehe2017} for monotonic fracture, the use of gradient plasticity in phase-field models overcomes mesh-sensitivity in the post-critical stage and avoids unrealistic localization patterns observed in the local plasticity counterpart. This behavior is embedded in the present model, as shown in figure~\ref{asn_prof} for varying mesh sizes, where the profiles of the equivalent plastic strains and the phase-field variable are plotted along a direction orthogonal to the crack in the center of the specimen.

% In particular, the control over the plastic localization region offered by the plastic length scale allows to avoid unphysical damage growth in regions outside the plastic zone. 

%It has already been shown that for a sufficiently fine discretization, mesh-objective results can be achieved (figure~\ref{asn_ss}). Thus, it is left to assess the ability ... assess the influence of these mechanisms in the present fatigue model, we consider in the present simulation the parameters \dots. Figure X shows the results of \dots. 
%Recall that the proposed model  gradient-extended plasticity two length scale parameters that govern the plasticity and damage localized responses. To highlight the influence of these mechanisms and their interplay ability of the coupled gradient-extended plasticity-damage model to, \dots

\subsubsection{Compact-tension test}

\begin{figure}[!ht]
\centering
\includegraphics[width=0.75\textwidth]{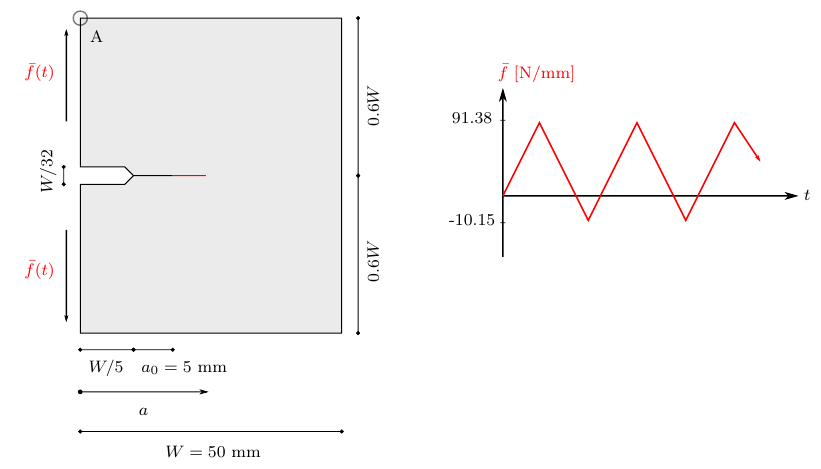}
%\def\svgwidth{0.55\textwidth}
%\import{}{asn_scheme.pdf_tex}
\caption{Schematic representation of the simplified compact-tension specimen with an initial crack of 5 mm under force loading.}
\label{plate3}
\end{figure}

The final example considers a typical fatigue test consisting of a compact-tension specimen subjected to  force loading. Figure~\ref{plate3} shows the geometry and loading conditions adopted for the simulations. For the sake of simplicity, as in~\citet{seiler2019}, a simplified test set-up is considered, where forces are applied on the left-most edges of the specimen as uniformly distributed loads. The load is discretized in 100 time steps per cycle.

Consider the parameters shown in table~\ref{tab:cts} and the AT-2 damage model (equation~\eqref{w0}). For the numerical simulations, an unstructured mesh of bilinear quadrilateral elements is employed, with local mesh refinement along the direction of the initial crack with a characteristic element size $h_\mathrm{min}=0.13$ mm.%, such that $\ell_\mathrm{d}/h_\mathrm{min}=1.95$. 

\begin{table}[!h]
\centering
\small
\caption{Material parameters for the compact-tension specimen under cyclic loading.}
\small \begin{tabular}{cccccccccccc}
\hline
$K$       & $\nu$   & $w_0$     & $\ell_\mathrm{d}$ & $\gamma_0$ & $k$  & $\sigma^\mathrm{p}$ & $H^\mathrm{iso}$ & $H^\mathrm{kin}$ & $\ell_{\mathrm{p}}$ & $\beta$ \\ 
{[}MPa{]} & {[}-{]} & {[}MPa{]} & {[}mm{]}   & {[}MPa{]}  & {[}-{]}         & {[}MPa{]}             & {[}MPa{]}          & {[}MPa{]}          & {[}mm{]}     & {[}-{]} \\ \hline
175000  & 0.3    &  45.23       &         0.25       & 500         & 0.4   & 125 / 480 / 720                   & 0                  & 36500                & 0.75                   & 0.1    \\ \hline  
\end{tabular}
\label{tab:cts}
\end{table}

In this example, we focus on the influence of ductility in the fatigue fracture process. Consequently, the simulation is computed with different values for the plastic yield strength $\sigma^\mathrm{p}$. Due to the loading conditions, forces are applied either until the initial crack length is increased to achieve a total length $a\approx 21.5$ mm (figure~\ref{plate3}), or until the maximum allowable force is reached at a given time step due to the damage/fatigue-induced cyclic softening response. 

\begin{figure}[!ht]
\centering
\includegraphics[scale=0.67, trim={0cm 0cm 0cm 0cm},clip]{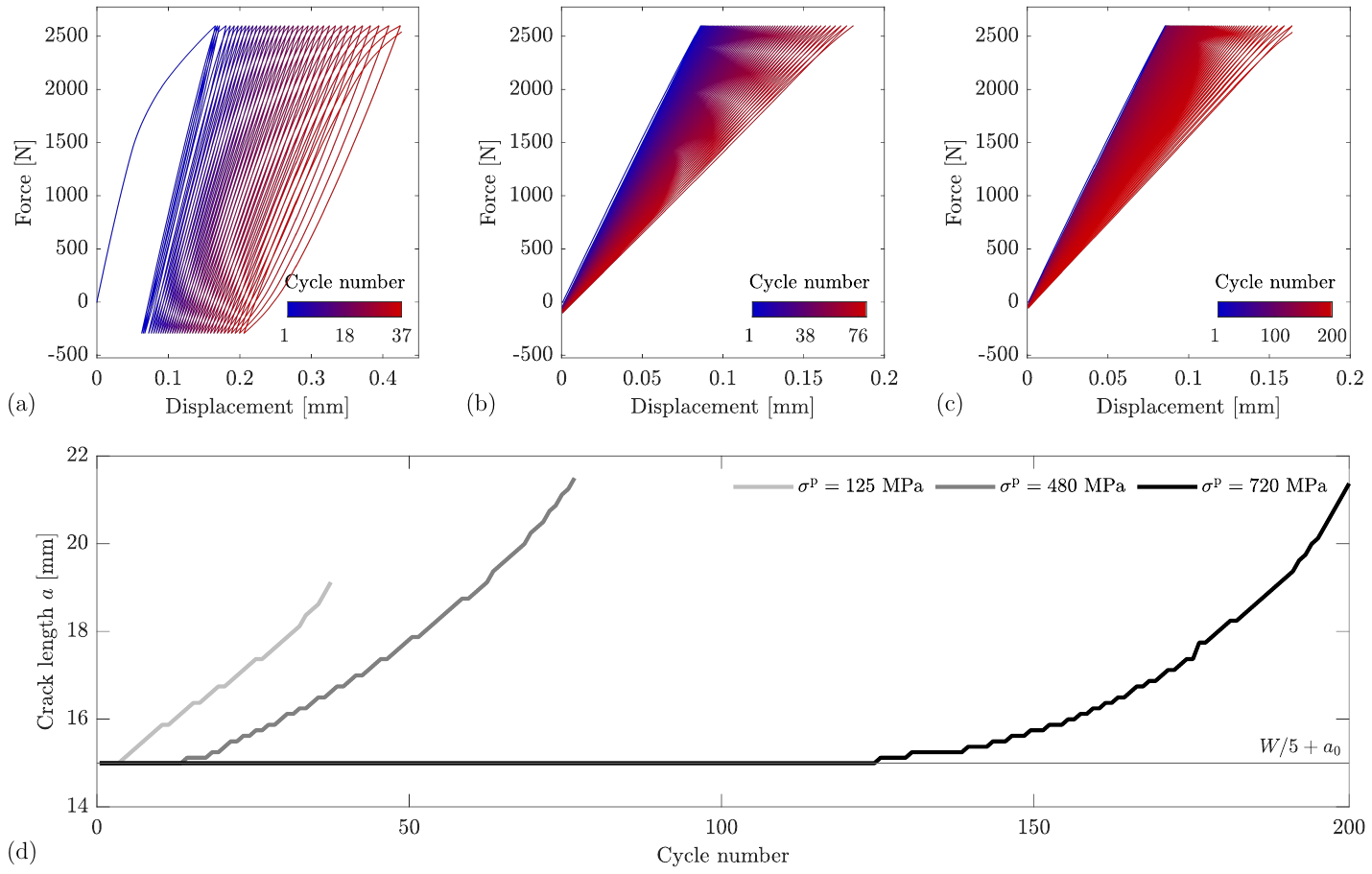}
\caption{Force-displacement curves at the top border (displacements at location A; figure~\ref{plate3}) for (a) $\sigma^\mathrm{p}=125$ MPa, (b) $\sigma^\mathrm{p}=480$ MPa and (c) $\sigma^\mathrm{p}=720$ MPa. The cyclic crack growth for the three cases is shown in (d).}
\label{cts1}
\end{figure}

The resulting force-displacement curves are shown in figure~\ref{cts1}, along with the cyclic crack-growth curves. Figure~\ref{cts1}(a) shows a ductile response, exhibiting considerable ratcheting effects. As observed in the homogeneous responses in section~\ref{sec:hom}, the ratcheting rate progressively grows as damage evolves and, in particular, as the crack propagates. The simulation ends after 37 cycles, where the maximum allowable force is reached, and the crack has achieved a total length $a=19.12$ mm. A markedly less ductile response is observed as $\sigma^\mathrm{p}$ increases, with the plastic loops exhibiting very limited ratcheting. The low accumulation of plastic strains thus delays the crack propagation process and allows the specimen to withstand much more cycles prior to attaining a maximum allowable force. For $\sigma^\mathrm{p}=480$ MPa, a total length $a=21.49$ mm is achieved after 76 cycles, while for $\sigma^\mathrm{p}=720$ MPa, a total length $a=21.37$ mm is achieved after 200 cycles.

\begin{figure}[!ht]
\centering
\includegraphics[scale=0.7, trim={0cm 0cm 0cm 0cm},clip]{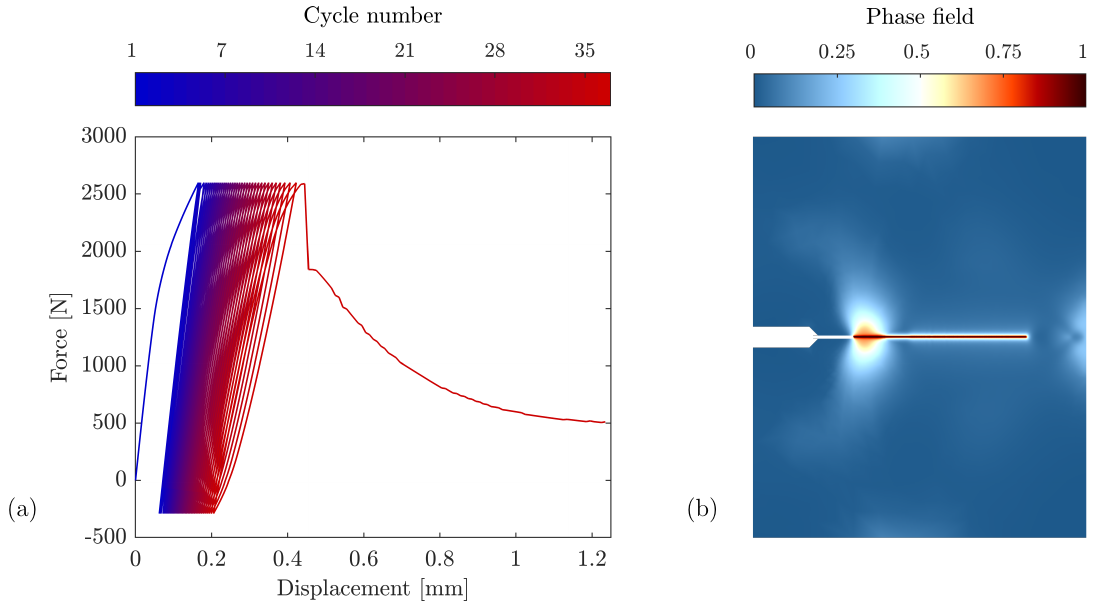}
\caption{(a) Force-displacement curve  measured at the top border (displacements at location A; figure~\ref{plate3}) and  (b) phase-field  contour plot for $\sigma^\mathrm{p}=125$ MPa. The initial cyclic force-loading stage is followed by monotonically imposed displacements.}
\label{cts3}
\end{figure}

Figure~\ref{cts3} shows the result of continuing the simulation after the cyclic force-loading stage for $\sigma^\mathrm{p}=125$ MPa, where monotonic displacements are imposed after the maximum allowable force has been reached. The contour plot shows the phase-field variable with a relatively large damage region around the crack tip during the cyclic loading process. This is followed first by an instant of brutal crack propagation and then by a stable propagation stage, where a narrow crack is observed to grow as displacements are monotonically~imposed.

\begin{figure}[!ht]
\centering
\includegraphics[scale=0.7, trim={0cm 0cm 0cm 0cm},clip]{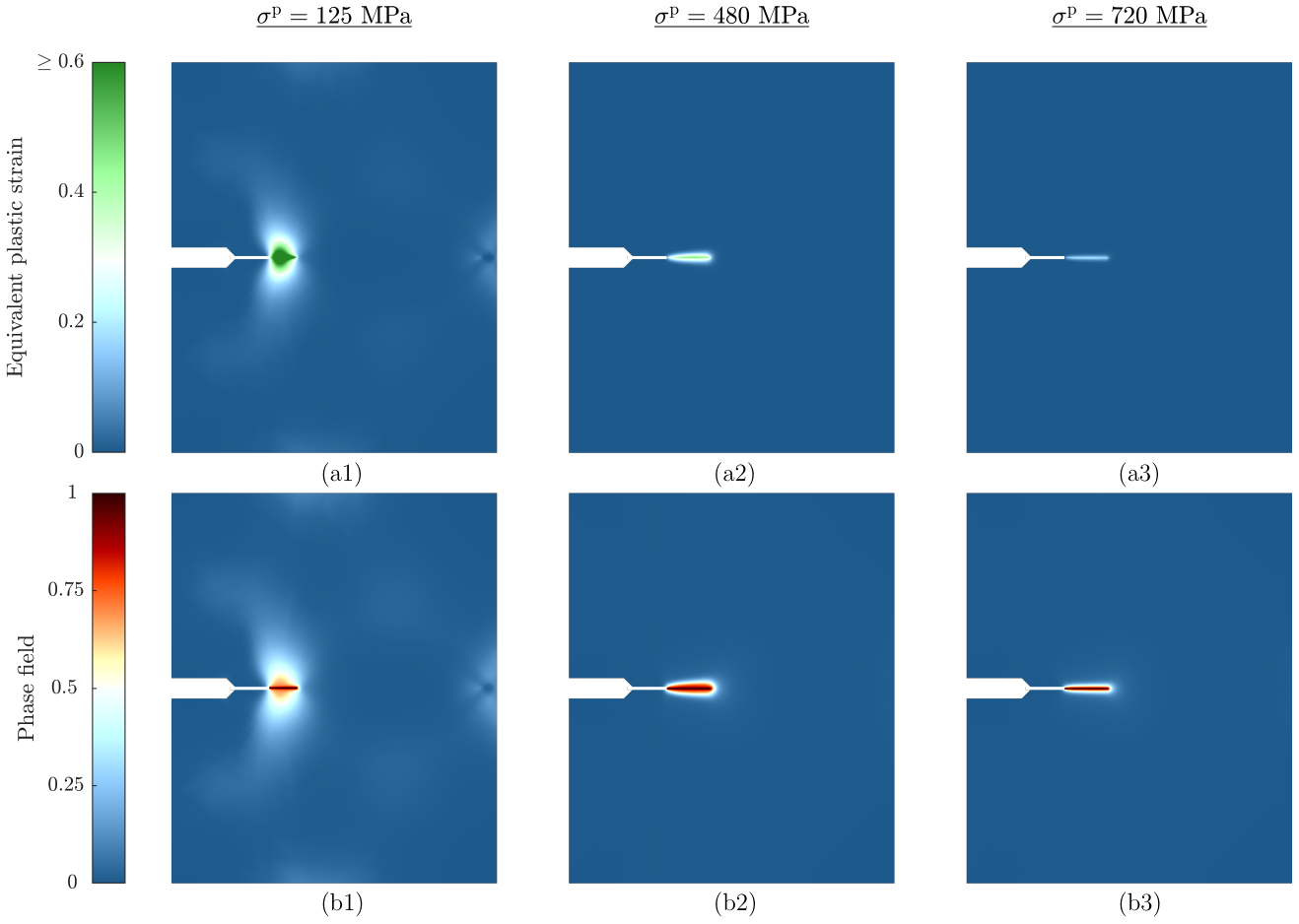}
\caption{Equivalent plastic strains (a1)-(a3) and phase field (b1)-(b3) for different yield strengths.}
\label{cts2}
\end{figure}

Figure~\ref{cts2} shows contour plots for the equivalent plastic strains and the phase-field/damage variable for the different yield strengths. As the response becomes less ductile, a markedly smaller plastic region is observed, resulting in a smaller damage zone around the crack tip.

\FloatBarrier

\section{Conclusions}

In the context of the energetic formulation, we have presented a coupled gradient-enhanced plasticity-damage model that embeds, in a unified way, characteristic features of low- and high-cycle fatigue. The proposed variational model is able to account for cyclic failure under both force and displacement loading by combining a phase-field description of fatigue cracks with cyclic plasticity, including multi-surface kinematic hardening, isotropic hardening/softening and ratcheting. The multi-field governing equations are derived from the principles of the energetic formulation, leading to a robust numerical implementation based on an alternate minimization scheme.

The results of numerical simulations indicate that several material responses can be captured by the cyclic plasticity model, including cyclic hardening and cyclic softening effects, as well as ratcheting under force loading and stress relaxation under displacement loading. Once damage coupling is introduced, the plastic energy accumulation entails a plastic fatigue effect, associated with the low-cycle fatigue regime. The model is further enriched by a fatigue degradation function that is driven by free energy accumulation and degrades the damage resisting force. In the absence of plastic strains, this feature accounts for  brittle fracture processes under high-cycle fatigue. The combination of elastic and plastic fatigue mechanisms conceived in an energetic framework allows for a physically sound description of a broad range of behaviors. These results are evidenced in the study of homogeneous uniaxial responses. Moreover, the results of numerical simulations in a 2D setting highlight the ability of the model to objectively  describe fatigue-induced ductile fracture, including the initiation, propagation and merging of ductile cracks.

The present study lays the groundwork for future theoretical developments and a broad range of applications. In addition to the most evident applications in civil and mechanical engineering, the modeling of cyclic inelastic behavior is relevant in a variety of fields ranging from geomechanics~\citep{suiker2003,niemunis2005,franccois2010} to electromechanics, e.g., to account for plastic deformations in electrode materials~\citep{brassart2013,peigney2018}, and biomechanics~\citep{martin2015,dong2020}, where plastic deformations are observed in biological tissues~\citep{Preziosi2010,Preziosi2011,Sciume2013}.  Concerning the proposed model, quantitative comparisons with experimental results are crucial to calibrate the various mechanisms included in the model and to assess their predictive ability. On the other hand, from a numerical perspective, the development of efficient techniques are of major interest to handle the computational cost of cyclic loading.

%The model includes two fatigue effects. The first aims to describe the main features of low-cycle fatigue by coupling the phase-field approach to fracture to a cyclic plasticity model that includes multi-surface kinematic hardening with gradient-enhanced isotropic hardening/softening and a ratcheting variable. The second effect is able to describe high-cycle fatigue, driven by elastic energy accumulation. The interplay between both fatigue mechanisms allows to capture a wide range of responses under both force control and displacement control, as shown in several numerical simulations. 
%highlight the capability of the proposed model to capture a wide range of responses.

\section*{Acknowledgements}

R. Alessi acknowledges the Italian Ministry of Education, University and Research (MIUR) under the PRIN~2017 20177TTP3S grant.

\clearpage

\small
\bibliography{bibliog}

\end{document}